%% file: main.tex
\documentclass[reqno]{amsart}
\usepackage{amsfonts}
\usepackage{amsmath}
\usepackage{amsthm, amscd}
\usepackage{mathtools}
\usepackage{amssymb}
\usepackage{mathrsfs}
\usepackage{graphicx}
\usepackage{caption}
\usepackage{subcaption}
\usepackage{float}
\usepackage{xypic}
\usepackage[abs]{overpic}
\usepackage[alphabetic]{amsrefs}
\usepackage{etex}
\usepackage{tikz-cd}

\usepackage{hyperref}
\hypersetup{
	colorlinks,
	linkcolor=[rgb]{0,0,0.7},
	urlcolor=[rgb]{0,0,0.4},
	citecolor=[rgb]{0.4,0.1,0}
}

\allowdisplaybreaks

\theoremstyle{plain}\newtheorem{Theorem}{Theorem}[section]
\theoremstyle{plain}\newtheorem{Corollary}[Theorem]{Corollary}
\theoremstyle{plain}\newtheorem{Lemma}[Theorem]{Lemma}
\theoremstyle{plain}\newtheorem{Definition}[Theorem]{Definition}
\theoremstyle{plain}\newtheorem{Proposition}[Theorem]{Proposition}
\theoremstyle{plain}\newtheorem{Assumption}[Theorem]{Assumption}
\theoremstyle{plain}
\theoremstyle{plain}\newtheorem*{Theorem*}{Theorem}
\theoremstyle{plain}

\theoremstyle{remark}\newtheorem{remark}[Theorem]{Remark}
\theoremstyle{remark}\newtheorem{Example}[Theorem]{Example}
\theoremstyle{remark}\newtheorem*{Notation*}{Notation}

\theoremstyle{plain}
\makeatletter
\newtheorem*{rep@theorem}{\rep@title}
\newcommand{\newreptheorem}[2]{
\newenvironment{rep#1}[1]{
 \def\rep@title{#2 \ref{##1}}
 \begin{rep@theorem}}
 {\end{rep@theorem}}}
\makeatother
\newreptheorem{Theorem}{Theorem}
\newreptheorem{Proposition}{Proposition}
\newreptheorem{Corollary}{Corollary}

\numberwithin{equation}{section}

\usepackage{lipsum}

\DeclareMathOperator{\Imm}{Im}

\DeclareMathOperator{\SO}{SO}

\DeclareMathOperator{\pt}{pt}
\DeclareMathOperator{\id}{id}

\DeclareMathOperator{\lk}{lk}

\DeclareMathOperator{\Diff}{Diff}
\DeclareMathOperator{\Homeo}{Homeo}
\DeclareMathOperator{\Emb}{Emb}
\DeclareMathOperator{\Map}{Map}
\DeclareMathOperator{\inte}{int}
\DeclareMathOperator{\Co}{Co}
\DeclareMathOperator{\Cone}{Cone}
\DeclareMathOperator{\sk}{sk}

\newcommand{\bQ}{\mathbb{Q}}
\newcommand{\bR}{\mathbb{R}}

\newcommand{\bZ}{\mathbb{Z}}

\newcommand{\cS}{\mathcal{S}}

\author{Jianfeng Lin}
\address{Yau Mathematical Sciences Center, Tsinghua University, Beijing 100084, China}
\email{linjian5477@mail.tsinghua.edu.cn}
\author{Yi Xie}
\address{Beijing International Center for Mathematical Research, Peking University, Beijing 100871, China}
\email{yixie@pku.edu.cn}
\author{Boyu Zhang}
\address{Department of Mathematics, The University of Maryland at College Park, Maryland 20742, USA}
\email{bzh@umd.edu}

\title{On the mapping class groups of 4-manifolds with 1-handles}

\begin{document}

\begin{abstract}
We develop a framework that generalizes Budney--Gabai's $W_3$ invariant on $\pi_0\textrm{Diff}(S^1\times D^3,\partial)$ to 4-manifolds with 1-handles. As applications, we show that if $M=(S^1\times D^3)\natural \hat M$ where $\hat M$ either has the form $I\times Y$  or is a punctured aspherical manifold, then the center of the mapping class group of $M$ is of infinite rank.
\end{abstract}

\maketitle

\section{Introduction}
\input{introduction}

\section{The scanning map $\mathcal{S}$}
\input{map_to_emb}

\section{The mapping space}
\label{sec_mapping_space}
\input{mapping_space}

\section{Homotopy groups of the configuration spaces}
\label{sec_homotopy_conf_space}
\input{homotopy_conf_space}

\section{Spectral sequence from a fibration tower}\label{sec_spec_seq}
\input{spectral_sequence}

\section{Proof of Theorem \ref{Thm-main}}
\label{secpi20}
\input{pi_2_hat_M_zero}

\section{The case $\hat M = I\times Y$}
\label{secIxY}
\input{I_times_Y}

\appendix
\section{The topology of simplicial compactifications}
\input{appendix}

\bibliographystyle{amsalpha}
\bibliography{references}

\end{document}

%% file: introduction.tex

Let $M$ be a compact smooth manifold with boundary. Define $\Diff(M,\partial )$ to be the group of diffeomorphisms of $M$ that are equal to the identity near $\partial M$. Endow $\Diff(M,\partial )$ with the $C^\infty$ topology. We use $\pi_i\Diff(M,\partial )$ to denote the $i$th homotopy group of $\Diff(M,\partial)$ with the base point $\id_M$. 
If $M\hookrightarrow M'$ is a codimension-zero smooth embedding, then it induces a homomorphism from $\Diff(M,\partial )$ to $\Diff(M', \partial )$ via an extension by the identity map. 

The following result was proved by Budney--Gabai and Watanabe independently. 
\begin{Theorem}[\cite{BG2019,Watanabe2020}]\label{Thm-BG}
The mapping class group $\pi_0 \Diff(S^1\times D^3,\partial)$ is an abelian group of infinite rank. 
\end{Theorem}

 In this paper, we study the image of $\pi_0 \Diff(S^1\times D^3,\partial)$ to $\pi_0 \Diff(M,\partial)$ induced by an embedding of $S^1\times D^3$ into $M$, where $M$ is a $4$--manifold. Our main results generalize Theorem \ref{Thm-BG} to a class of 4-manifolds with 1-handles. 

\begin{Theorem}\label{Thm-main}
	Suppose $\hat M$ is a smooth connected compact $4$--manifold with boundary, let $p_0\in \partial \hat M$ be a fixed point. Let $M = (S^1\times D^3)\natural \hat M$ be the boundary connected sum of $S^1\times D^3$ with $\hat M$ at $p_0$, and view $S^1\times D^3$ as a codimension-0 submanifold of $M$. Assume the following conditions hold simultaneously:
    \begin{enumerate}
    \item $\pi_{2}^\bQ(\hat{M}) =0$. 
    \item\label{item_main_pi1_pi3_not intersect} Let $b_1,b_2$ be distinct points in the connected component of $\partial \hat M$ containing $p_0$. For every pair of elements $\alpha\in \pi_{1}(\hat{M},b_1)$ and $\beta\in \pi_{k}(\hat{M},b_2)$ with $k=3$ or $4$, there is a map $\gamma:(D^{1},\partial D^{1})\to (\hat{M},b_1)$ that represents $\alpha$, and a map $\mu: (D^{k},\partial D^{k})\to (\hat{M},b_2)$ that represents \emph{a nonzero multiple} of $\beta$, such that \[\gamma(D^1)\cap \mu(D^k)=\emptyset.\]
        \end{enumerate}
    Then the image of  $\pi_0 \Diff(S^1\times D^3,\partial )$ in $\pi_0 \Diff(M,\partial)$ is of infinite rank.
\end{Theorem}
Here and in the following, we use $\pi_i^\bQ$ $(i\ge 2)$ to denote the tensor product of the $i^{th}$ homotopy group with $\bQ$. 

\begin{remark}
    Note that the image of $\pi_0 \Diff(S^1\times D^3,\partial )$ in $\pi_0 \Diff(M,\partial)$ is contained in the center of $\pi_0 \Diff(M,\partial)$. This is because if $f$ is in this image and $g$ is in $\pi_0 \Diff(M,\partial)$, one can isotope $f$ so that its support is contained in a collar neighborhood of $\partial M$ that is disjoint from the support of $g$. Therefore, Theorem \ref{Thm-main} implies that the center of $\pi_0 \Diff(M,\partial)$ is an abelian group of infinite rank. 
\end{remark}

We give a few examples such that the conditions of Theorem \ref{Thm-main} hold.
\begin{Example}
  If $\hat M$ is aspherical, then the conditions of Theorem \ref{Thm-main} hold.
\end{Example}

\begin{Example}
    Suppose $\pi_2^\bQ(\hat M) = 0$ and $\pi_1(\partial \hat M,p_0)\to \pi_1(\hat M,p_0)$ is surjective. In this case, Condition (2)  of Theorem \ref{Thm-main} holds because we may take the representative of $\alpha$ to be contained in $\partial \hat M$, and take the representative of $\beta$ to be in the interior of $\hat M$ except for the base point. 
\end{Example}

\begin{Example}
    Suppose $X$ is a compact, connected, aspherical 4-manifold with $\partial M\neq \emptyset$ and $\hat M$ is the closure of $X\setminus (\cup_{1\le i \le m} B_i)$ in $X$, where $\{B_i\}_{1\le i\le m}$ is a collection of disjoint embedded closed balls in the interior of $X$. Suppose $p_0\in \partial X$. Then the universal cover of $\hat M$ is homotopy equivalent to the wedge sum of a family of $S^3$'s, so $\pi_2^\bQ(\hat M)\cong \pi_4^\bQ(\hat M) =0$. For each $\partial B_i$, let $q_i\in \partial B_i$ be a fixed point, let $\eta_i\in \pi_3(\partial B_i,q_i)$ be a generator, let $\gamma_i$ be a fixed arc from $b_2$ to $q_i$, and let $\eta_i'\in \pi_3(\hat M,b_2)$ be represented by the $\gamma_i$--action on $\eta_i$. Then $\pi_3(\hat M,b_2)$ is generated by $\{\eta_i'\}$ and their actions under $\pi_1(\hat M)$, so Condition (2)  of Theorem \ref{Thm-main} holds.

    In particular, if $X=D^4$, then $M$ is $S^1\times D^3$ with finitely many interior balls removed. 
\end{Example}

In order to state the next theorem, 
let $\Homeo(M,\partial)$ be the group of homeomorphisms of $M$ that are equal to the identity near $\partial M$, and endow $\Homeo(M,\partial)$ with the compact open topology. 

\begin{Theorem} 
\label{Thm-main-2}
Suppose $Y$ is a connected compact $3$--manifold with a non-empty boundary. Let $I=[0,1]$. Let $\hat{M}=I\times Y$ and let  $M = (S^1\times D^3)\natural \hat M$. Then the image of $\pi_0 \Diff(S^1\times D^3,\partial )$ in $\pi_0 \Homeo(M,\partial)$ is of infinite rank.   
\end{Theorem}

\begin{remark}
Under the assumptions of Theorem \ref{Thm-main-2}, we have $M\cong I\times ((S^1\times D^2)\natural Y)$.  For $f,g\in  \Homeo(M,\partial)$, one can always isotope them so that the support of $f$ is contained in $(0,1/2)\times ((S^1\times D^2)\natural Y)$ and the support of $g$ is contained in $(1/2,1)\times ((S^1\times D^2)\natural Y)$. As a result, $\pi_0 \Homeo(M,\partial)$ and $\pi_0\Diff(M,\partial)$ are both abelian. Therefore, Theorem \ref{Thm-main-2} implies that $\pi_0 \Homeo(M,\partial)$ and $\pi_0\Diff(M,\partial)$ are abelian groups of infinite rank. 
\end{remark}

Budney and Gabai’s proof of Theorem \ref{Thm-BG} involves ideas from embedding calculus developed by
Weiss \cite{Weiss1}, Goodwillie--Weiss \cite{GW}, Sinha \cite{Sinha-cpt,Sinha-cosimplicial} and Goodwillie--Klein \cite{GK}. Watanabe’s proof, on the other hand, depends on the configuration space integral introduced by Kontsevich. Both proofs make use of the configuration space $C_n(M)$ of ordered $n$ points in a manifold $M$ and their Fulton--MacPherson compactifications. It is believed that this may not be a coincidence. 
In fact, Randal--Williams \cite{RW} suggests that the embedding calculus and the configuration space integrals are closely related and in some cases equivalent.

The argument of Budney and Gabai can be summarized as follows.
In \cite{Sinha-cosimplicial}, a Taylor tower of topological spaces was constructed to describe the homotopy type of the embedding space \( \Emb(I, M) \) (see Section \ref{subsec_embedding_space} for the definition of \( \Emb(I, M) \)). Using the third stage \( T_3 \Emb(I, S^1 \times D^3) \) of this tower, Budney and Gabai defined an invariant \( W_3 \) on \( \pi_2^{\mathbb{Q}} \Emb(I, S^1 \times D^3) \), which takes values in a quotient group of \( \pi_5^{\mathbb{Q}}(C_3(S^1 \times D^3)) \). They also defined a scanning map \( \cS: \pi_0 \Diff(S^1 \times D^3, \partial) \to \pi_2 \Emb(I, S^1 \times D^3) \) (see Section \ref{subsec_S}). Let $i_\bQ$ be the canonical map from $\pi_2 \Emb(I, S^1 \times D^3)$ to $\pi_2 ^\bQ \Emb(I, S^1 \times D^3)$.  The proof of Theorem \ref{Thm-BG} is finished by showing that the image of the composition map $W_3\circ i_\bQ\circ \cS$ is of infinite rank.

In this paper, we generalize the constructions of Budney and Gabai to $4$--manifolds with $1$--handles. We give an interpretation of the Budney--Gabai's $W_3$ map using the Bousfield--Kan spectral sequence of a fibration tower (see Section \ref{sec_spec_seq}). This allows us to generalize the map $W_3$ to arbitrary manifolds $M$. When $M$ has the form $(S^1\times D^3)\natural \hat M$, the images of  diffeomorphisms supported in $S^1\times D^3$ under the map $W_3\circ i_\bQ\circ \mathcal{S}$ take values in a quotient group of $\pi_5^\bQ C_3 (M)$. 
Theorem \ref{Thm-main} is proved by a careful analysis of the spectral sequence.  Theorem \ref{Thm-main-2} is proved by generalizing an argument in \cite{BG2023} using similar ideas.

Throughout this paper, all manifolds are assumed to be smooth.

\subsection*{Acknowledgments} We would like to express our sincere gratitude to Budney, Gabai, and Watanabe for helping us understand their work. 
J.L. is partially supported by
NSFC Grant 12271281. 
Y.X. is partially supported by National Key R\&D Program of China 2020YFA0712800, NSFC 12071005 and NSFC 12341105.
B.Z. is partially supported by NSF grant
DMS-2405271 and a travel grant from the Simons Foundation.

%% file: map_to_emb.tex
In this section, we define the scanning map $\mathcal{S}: \pi_0\Diff(M,\partial)\to \pi_2 \Emb(I,M)$ when $M$ has the form $(S^1\times D^3)\natural \hat M$. By definition, the map $\mathcal{S}$ is a priori a set-theoretic map and may not be a homomorphism. 

\subsection{The embedding space}
\label{subsec_embedding_space}
We start by defining the space $\Emb(I, M)$ and set up some notations along the way. 

Let $M$ be a compact $4$--manifold with boundary. Let $p_0, p_1$ be distinct fixed points on $\partial M$. Fix two unit vectors $v_0\in T_{p_0}M$, $v_1\in T_{p_1}M$ such that $v_0$ and $v_1$ are transverse to $T\partial M$, the vector $v_0$ points inward to $M$, and the vector $v_1$ points outward from $M$. Let $I=[0,1]$. Let $\Emb(I,M)$ denote the space of smooth embeddings $i:I\to M$ such that $i(0)=p_0$, $i(1) = p_1$, $i'(0)=v_0$, $i'(1) = v_1$, and that the image of $(0,1)\subset I$ under $i$ is contained in the interior of $M$. Endow $\Emb(I,M)$ with the $C^\infty$--topology. It is clear that the homeomorphism type of the space $\Emb(I,M)$ does not depend on the choice of $p_1,p_2, v_1, v_2$. 

For later reference, we define another embedding space  $ \Emb^* (I,M)$ to be the space of all smooth embeddings $i:I\to M$ such that $i(0)=p_0$, $i(1)=p_1$, the vectors $i'(0), i'(1)$ are non-zero and transverse to $\partial M$, the vector $i'(0)$ points inward to $M$, the vector $i'(1)$ points outward from $M$, and the image of $(0,1)\subset I$ in $M$ is contained in the interior of $M$. Endow $ \Emb^* (I,M)$ with the $C^\infty$--topology.  
For $p\in \partial M$, let $T_p^+ M$ be the subset of $T_pM$ consisting of non-zero vectors transverse to $T_p\partial M$ and pointing outward from $M$, and let $T_p^- M$ be the subset of $T_pM$ consisting of non-zero vectors transverse to $T_p \partial M$ and pointing inward to $M$.
Then
$$
\Emb(I,M) \to \Emb^*(I,M) \to T_{p_0}^-M \times T_{p_1}^+M 
$$
is a Serre fibration. Hence the embedding of $\Emb(I,M)$ into $\Emb^*(I,M)$ is a weak homotopy equivalence. 

\subsection{The map $\cS$}
\label{subsec_S}
Now assume $M$ has the form $(S^1\times D^3) \natural \hat M$ where $\hat M$ is a connected compact $4$--manifold with boundary. We define a set-theoretic map $\mathcal{S}$ from
 $\pi_0\Diff(M,\partial)$ to $\pi_2(\Emb(I,M))$.

View $S^1\times D^3$ as a codimension--$0$ submanifold of $M$. Note that $\partial (S^1\times D^3)\setminus \partial M$ is an open $3$--ball.
Without loss of generality, we may assume there exists $t_0\in S^1$ such that $\{t_0\}\times \partial D^3$ as a subset of $S^1\times D^3$ is contained in $\partial M\cap \partial (S^1\times D^3)$. Fix such a $t_0$, and let $p_0$, $p_1$ be fixed points on $\{t_0\}\times \partial D^3$. 
Let unit vectors $v_0\in T_{p_0}M$ and $v_1\in T_{p_1}M$ be transverse to $T\partial M$ such that $v_0$ points inward to $M$ and $v_1$ points outward from $M$. Define $\Emb(I,M)$ and $\Emb^*(I,M)$ with respect to $p_0,p_1,v_0,v_1$. 

Let $\inte(-)$ denote the interior of a manifold. 
There exists a smooth family of embedded arcs in $\Emb^*(I,M)$, parametrized by an open $2$--disk, such that the images of $\inte(I)$ sweep through $\{t_0\}\times \inte(D^3)$.  Choose such a family and denote it by 
\[
\iota:\inte(D^2)\to \Emb^*(I,M).
\]
We choose $\iota$ so that there exists a unique point $\xi_0\in\inte(D^2)$ such that $\iota(\xi_0)\in \Emb(I,M)$.  We also assume that $p_0,p_1,v_0,v_1,\iota$ are chosen so that $\iota(\xi_0)$ is given by a linear map from $I$ into $\{t_0\}\times D^3$. Fix $p_0,p_1,v_0,v_1,\iota$ from now.

For each $f\in \Diff(M,\partial)$, let $f_*\iota:\inte(D^2)\to \Emb^*(I,M)$ denote the map such that $(f_*\iota)(\xi) = f\circ (\iota(\xi))$ for all $\xi\in \inte(D^2)$. Since $f$ is equal to the identity near $\partial M$, the maps $\iota$ and $f_*\iota$ are the same outside a compact subset of $\inte(D^2)$. 

Identify $D^2$ with a unit disk in $\bR^2$, and let $D^2_r$ denote the concentric disk with radius $r$. Consider the following maps on two copies of $\inte(D^2)$:
\begin{align*}
    \iota: \inte(D^2) &\to \Emb^*(I,M),\\
    f_*\iota: \inte(D^2) &\to \Emb^*(I,M).
\end{align*}
Then for $\epsilon>0$ sufficiently small, the two maps agree on $\partial D^2_{1-\epsilon}$ and hence can be glued to define a map from $S^2$ to $\Emb^*(I,M)$. Let $\xi_0$ from the first copy of $\inte(D^2)$ be the base point, then the glued map defines an element in $\pi_2(\Emb^*(I,M), \iota(\xi_0))$. It is straightforward to verify that when $\epsilon$ is sufficiently small, this element in $\pi_2(\Emb^*(I,M), \iota(\xi_0))$ is independent of the choice of $\epsilon$, and it also remains invariant under a homotopy of $f$ within $\Diff(M,\partial)$.  

As a result, the above construction defined a set-theoretic map from $\pi_0(\Diff(M,\partial))$ to $\pi_2(\Emb^*(I,M), \iota(\xi_0))$. Since $\pi_2(\Emb^*(I,M), \iota(\xi_0))$ is canonically isomorphic to $\pi_2(\Emb(I,M), \iota(\xi_0))$, we obtained a map from $\pi_0(\Diff(M,\partial))$ to $\pi_2(\Emb(I,M), \iota(\xi_0))$. We denote this map by $\cS$. For notational convenience, if $f\in \Diff(M,\partial)$ and $[f]$ is its homotopy class in $\Diff(M,\partial)$, we will also use $\cS(f)$ to denote $\cS([f])\in \pi_2(\Emb(I,M), \iota(\xi_0))$.

\begin{Lemma}
\label{lem_hom_cS_S1xD3}
	Assume $M = (S^1\times D^3) \natural \hat M$ and $f_1,f_2\in \Diff(M,\partial M)$ are supported in $S^1\times D^3$. Then $\cS(f_1\circ f_2) = \cS(f_1)+\cS(f_2)$. 
\end{Lemma}

\begin{proof}
After isotopy, we may assume that there are two disjoint balls $B_1,B_2$ in the interior of $D^3$ such that $f_1$, $f_2$ are supported in $S^1\times B_1$ and $S^1\times B_2$ respectively. Recall that $\xi_0\in \inte (D^2)$ is the base point. After further isotopies of $f_1$ and $f_2$ if necessary, we may assume that $\iota(\xi_0)$, which is a straight line segment in $\{t_0\}\times D^3$, is disjoint from both $\{t_0\}\times B_1$ and $\{t_0\}\times B_2$. We may also assume that the disk $D^2$ can be decomposed into a union of two half-disks $D^2 = S\cup T$ such that $S$ and $T$ have disjoint interiors, $\xi_0\in S\cap T$, the arcs in $\iota(T)$ are disjoint from $\{t_0\}\times B_2$, and the arcs in $\iota(S)$ are disjoint from $\{t_0\}\times B_1$.

Let $\mu_1: S^2\to \Emb^*(I,M)$ be defined by gluing $\iota$ and $(f_2)_*\iota$, let $\mu_2: S^2\to \Emb^*(I,M)$ be defined by gluing $(f_1)_*\iota$ and $(f_1\circ f_2)_*\iota$. We only need to show that $\mu_1$ and $\mu_2$ are homotopic relative to the base point $\xi_0$. 

Let $S_{1-\epsilon} = S\cap D^2_{1-\epsilon}$, $T_{1-\epsilon} = T\cap D^2_{1-\epsilon}$. Then $\mu_1$, $\mu_2$ are defined by maps on the doubles of $S_{1-\epsilon} \cup T_{1-\epsilon}$. For notational convenience, we will denote the two copies of $S_{1-\epsilon}$, $T_{1-\epsilon}$ by $S_{1-\epsilon}^{(i)}$, $T_{1-\epsilon}^{(i)}$ for $i=1,2$.

Let $\hat S = S_{1-\epsilon}^{(1)} \cup S_{1-\epsilon}^{(2)}$. Then $\hat S$ is a closed disk in the double of $S_{1-\epsilon} \cup T_{1-\epsilon}$. 
Note that by definition, $\mu_1$ and $\mu_2$ only differ on a compact subset of $\inte(\hat S)$. Let $\tau$ be the reflection on $\hat S$ that exchanges the two copies of $S_{1-\epsilon}$. Then both $\mu_1$ and $\mu_2$ are invariant under the reflection $\tau$, and $\mu_1=\mu_2$ near the boundary of $\hat S$. Therefore, $\mu_1|_{\hat S}$ and $\mu_2|_{\hat S}$ are homotopic relative to the boundary. Since the base point is in the complement of $\hat S$, the desired result is proved.
\end{proof}

Budney--Gabai proved Theorem \ref{Thm-BG} by showing that the image of $\pi_0\Diff(S^1\times D^3,\partial )$ in $\pi_2\Emb(I,S^1\times D^3)$ under the map $\cS$ is of infinite rank. In the following, we consider the homomorphism from $\pi_2\Emb(I,S^1\times D^3)$ to $\pi_2\Emb(I,M)$ induced by the embedding of $S^1\times D^3$ in $M$, and study the image of $\pi_0\Diff(S^1\times D^3,\partial )$ in $\pi_2\Emb(I,M)$.

%% file: mapping_space.tex
In this section, we define a family of maps from $\Emb(I,M)$ to certain mapping spaces (see Definition \ref{defn_Map_n}). These mapping spaces are the terms in the Taylor tower of $\Emb(I,M)$ 
as constructed in \cite{Sinha-cosimplicial}. Although this construction is closely related to embedding calculus, it does not logically depend on any results from embedding calculus, and we will give a self-contained exposition.

\subsection{Compactifications of configuration spaces}
In this subsection, we construct a sequence of mapping spaces from the configuration spaces of points on $M$ as approximations of $\Emb(I,M)$. We follow \cite{Sinha-cosimplicial} on the notations for configuration spaces and their compactifications. All the notations are also reviewed in the Appendix.

Given a smooth manifold $M$,
we use $C_n(\inte(M))$ to denote the configuration space of (ordered) $n$ points in $\inte(M)$. 
Let $C_n\langle M,\partial \rangle$ denote the boundary version of the simplicial compactification of $C_n(\inte(M))$ as defined in \cite[Definition 4.10]{Sinha-cosimplicial}\footnote{Our notation $C_n\langle M,\partial \rangle$ corresponds to the notation $C_n\langle [M,\partial] \rangle$ in \cite{Sinha-cosimplicial}}. See Definition \ref{def_compactify_w_boundary} in the appendix for a review of the definition of $C_n\langle M,\partial \rangle$.

Following \cite[Definition 4.9]{Sinha-cosimplicial}, let $C_n'(\inte(M))$ and $C_n'\langle M,\partial\rangle$ be defined as the pullbacks of the bundle $(STM)^n\to M^n$ under the canonical maps
$C_n(\inte(M))\to M^n$ and $C_n\langle M,\partial \rangle \to M^n$ respectively.
Then $C_n'(\inte(M))$ is a bundle over $C_n(\inte(M))$ and $C_n'\langle M,\partial \rangle$ is a bundle over $C_n\langle M, \partial\rangle$, both with fiber $(S^{3})^n$. 

By Proposition \ref{prop_homotopy_equiv_Q}, the embedding of $C_n'(\inte(M))$ into $C_n'\langle M,\partial \rangle$ is a homotopy equivalence. 

By \cite[Corollary 4.22]{Sinha-cosimplicial}, the sequence of spaces $\{C_n'\langle M, \partial \rangle\}_{n\ge 1}$ has a cosimplicial structure. We briefly review the definition of the coface maps and the codegeneracy maps. For $i=1,\dots,n$, the codegeneracy map $\sigma^i_n: C_n'\langle M, \partial  \rangle \to C_{n-1}'\langle M, \partial \rangle$ is defined by forgetting the $i^{th}$ point in the configuration. For $i=1,\dots,n$, the coface map $\delta^i_n: C_n'\langle M\rangle \to C_{n+1}'\langle M \rangle$ is defined by doubling the $i^{th}$ point in the direction of the corresponding unit tangent vector. 
To describe the maps $\delta^0_n$ and $\delta^{n+1}_n$, let $p_0,p_1,v_0,v_1$ be as in the definition of $\Emb(I,M)$. A point in $C_n\langle M,\partial \rangle$ is given by a configuration of $(n+2)$ points $x_0,\dots,x_{n+1}$, where $x_0=p_0$ and $x_{n+1}=p_{1}$.
Then $\delta^0_n: C_n'\langle M\rangle \to C_{n+1}'\langle M \rangle$ is defined to be the map that doubles the point $x_0$ along the direction of $v_0$, and $\delta^{n+1}_n: C_n'\langle M\rangle \to C_{n+1}'\langle M \rangle$ is defined to be the map that doubles the point $x_{n+1}$ along the direction of $v_1$.
The maps $\delta^i_n$ are cofibrations.

Recall that $I=[0,1]$.
Let $\widetilde C_n'\langle I,\partial \rangle$ be the connected component of $C_n'\langle I,\partial \rangle$ such that the $n+2$ points are given in increasing order and all the tangent vectors are positive. The space $\widetilde C_n'\langle I,\partial\rangle$ is diffeomorphic to the $n$ dimensional simplex.

\subsection{The mapping spaces}
\begin{Definition}
   Let $X$ be a topological space. 
A \emph{$\Delta^n$ structure} on $X$ is defined to be the choice of a subset $X_S$ of $X$ for each $S\subset \{0,\dots, n\}$, such that $X_S\cap X_{S'} = X_{S\cap S'}$ for all $S,S'\subset \{0,\dots, n\}$, and that $X_{\{0,\dots,n\}} = X$, $X_{\emptyset} = \emptyset$.
\end{Definition}

The $n$--dimensional simplex has a $\Delta^n$ structure, where the subspace associated with $S\subset \{0,\dots, n\}$ is defined to be the convex hull of the vertexes indexed by $S$. 

The space $X=C_n'\langle M,\partial \rangle$ has a $\Delta^n$ structure defined as follows. For each $S\subset \{0,\dots, n\}$, the inclusion of $S$ in $\{0,\dots, n\}$ defines a map from $C_{(\# S)-1}'\langle M,\partial \rangle$ to $C_n'\langle M,\partial \rangle$ by the cosimplicial structure, and we define $X_S$ to be the image of this map. By definition, $X_{\{0,\dots,n\}} = X$ and  $X_{\emptyset} = \emptyset$. It is also straightforward to verify that $X_S\cap X_{S'} = X_{S\cap S'}$ for all $S,S'$. Note that by definition, each $X_{S}$ is the image of a composition of coface maps, so the inclusions of $X_S$ in $X$ are cofibrations. 

\begin{Definition}
\label{defn_Map_n}
Let $\Map_n(M)$ denote the space of continuous maps from $\widetilde C_n'\langle I,\partial \rangle$ to $C_n'\langle M,\partial \rangle$ that preserves the $\Delta^n$ structures, endowed with the compact open topology.
\end{Definition}

We have a canonical evaluation map
\begin{align}
\Psi_n: \Emb(I,M) &\to \Map_n(M) \label{eqn_def_Psi_n}\\
            f &\mapsto f_*, \nonumber
\end{align}
where $f_*:  \widetilde C_n'\langle I,\partial\rangle \to C_n'\langle M,\partial \rangle$ is the map on configuration spaces induced by the embedding $f$.

\begin{Definition}
    Suppose $\{X_S\}_{S\subset \{0,\dots,n\}}$ is a $\Delta^n$ structure on $X$. For $m\le n$, define the \emph{$m$--skeleton} of $X$ to be the union of $X_S$ for all $S$ such that $\#S \le m+1$. Denote the $m$--skeleton of $X$ by $\sk^{m}X$. 
\end{Definition}

\begin{Definition}
Suppose $\{X_S\}_{S\subset \{0,\dots,n\}}$ and $\{Y_S\}_{S\subset \{0,\dots,n\}}$ are two spaces with $\Delta^n$ structures. We say that a map $f:\sk^{m}X\to \sk^{m}Y$ \emph{preserves the $\Delta^n$ structure on skeletons}, if $f(X_S)\subset Y_S$ for all $S$ such that $\# S\le m+1$. 
\end{Definition}

Let $\Psi_n\circ\iota (\xi_0)$ be the image of the base point of $\Emb(I,M)$ in $\Map_n(M)$. 

\begin{Definition}
     For $m\le n$, let $\Map_{n,m}(M)$ denote the space of maps from $\sk^{m}\widetilde C_n'\langle[I,\partial ]\rangle$ to $\sk^{m}C_n'\langle[M,\partial ]\rangle$ that preserve the $\Delta^n$ structure on skeletons and that are homotopic to the restriction of $\Psi\circ \iota(\xi_0)$ to $\sk^{m}\widetilde C_n'\langle[I,\partial ]\rangle$. Endow $\Map_{n,m}(M)$ with the compact open topology.
\end{Definition}

The space $\Map_{n,m}(M)$ is path connected, and the restriction map 
$$
\Map_{n,m}(M) \to \Map_{n,m-1}(M)
$$
is a fibration. Note that by definition, $\Map_n(M) = \Map_{n,n}(M)$.

\begin{remark}
By \cite[Theorem 5.4, Lemma 5.12]{Sinha-cpt}, the map $\Psi_3: \Emb(I,M) \to \Map_3(M)$ induces an isomorphism on $\pi_i$ for all $i\le 2$.
This result suggests that we may use $\Map_3(M)$ to study the image of $\pi_0\Diff(M,\partial)$ under the map $\cS$. 
In fact, we have the following commutative diagram

\begin{equation}\label{diagram_S}
\begin{tikzcd}
\pi_0\Diff(S^1\times D^3, \partial) \arrow[r, "\cS"] \arrow[d,"E"] &  \pi_2\Emb(I, S^1\times D^3) \arrow[d] 
\arrow[r,"(\Psi_{3})_\ast"] & \pi_2\Map_3(S^1\times D^3) \arrow[d]\\
\pi_0\Diff(M, \partial) \arrow[r, "\cS"] & \pi_2\Emb(I, M) \arrow[r,"(\Psi_{3})_\ast"] & \pi_2\Map_3(M)
\end{tikzcd}
\end{equation}
where the first vertical map $E$ is defined by extending the diffeomorphisms on $S^1 \times D^3$ by the identity outside of 
$S^1 \times D^3$. It is proved in \cite{BG2019} that the image of $(\Psi_{3})_\ast\circ \cS$ from $\pi_0\Diff(S^1\times D^3, \partial)$ is of infinite rank.  
In Sections \ref{secpi20} and \ref{secIxY}, we will study the rightmost vertical map and show that 
the image of $(\Psi_{3})_\ast\circ \cS \circ E$ is also of infinite rank.
\end{remark}

\subsection{Homotopy groups with base points in a simply connected set}
\label{subsec_base_point_simply_conn}
For later reference, we introduce a notion of homotopy groups with base points in a simply connected subset. Suppose $X$ is a topological space and $B$ is a non-empty simply connected subset of $X$. Suppose $b_1,b_2\in B$ are distinct points and let $\gamma$ be an arc from $b_1$ to $b_2$ in $B$, then $\pi_i(X,b_1)$ and $\pi_i(X,b_2)$ are isomorphic via $\gamma$. Since $B$ is simply connected, the isomorphism does not depend on the choice of $\gamma$. Hence the homotopy groups $\pi_i(X,b)$ are canonically isomorphic for all $b\in B$. When the set $B$ is clear from the context, we will denote the canonically identified homotopy groups of $X$ with base points in $B$ by $\pi_i(X)$. 

Suppose $B$ is a simply connected subset in $X$ and $B'$ is a simply connected subset in $X'$. Suppose $f:X\to X'$ is a continuous map such that $f(B)\subset B'$. Then $f$ induces a homomorphism from $\pi_i(X)$ to $\pi_i(X')$, where the base point of $\pi_i(X)$ is taken in $B$, and the base point of $\pi_i(X')$ is taken in $B'$.

Suppose $M=(S^1\times D^3)\natural \hat M$ as above.
We construct a simply connected subset of $C_n'\langle M,\partial \rangle$ as the set of base points. Recall that $\iota(\xi_0)$ denotes the base point of $\Emb(I,M)$, and it is given by a linear embedding of $I$ into $\{t_0\}\times D^3$. 
For each $n$, let $B_n$ be the image of $\Psi_n(\iota(\xi_0))$, which is a subset of $C_n'\langle M,\partial \rangle$. 
Then $B_n$ is homeomorphic to an $n$--dimensional simplex. Geometrically, it consists of configurations $(x_1,\dots,x_n)$ such that all the points are contained in the image of $\iota(\xi_0)$ (which is a linear arc in $\{t_0\}\times D^3$) in order and all the directional vectors are parallel to $\iota(\xi_0)$ and point at the positive direction. For the rest of the paper, we take the set $B_n$ to be the set of base points for $C_n'\langle M,\partial \rangle$. We use $\pi_i C_n'\langle M,\partial\rangle$ to denote the $i^{th}$ homotopy group of $C_n'\langle M,\partial\rangle$ with base points in $B_n$. 

Let $\delta^i_n$, $\sigma^i_n$ be the coface and the codegeneracy maps, then $\delta^i_n(B_{n-1}) \subset B_n$, $\sigma^i_n(B_{n+1}) = B_n$. So  $\sigma^i_n$ induces a homomorphism $\pi_i C_n'\langle M,\partial\rangle \to \pi_i C_{n-1}'\langle M,\partial \rangle$, and $\delta^i_n$ induces a homomorphism  $\pi_i C_n'\langle M,\partial\rangle \to \pi_i C_{n+1}'\langle M,\partial \rangle$. 

We will also frequently consider the homotopy groups of $C_n \langle M,\partial \rangle$, $C_n '(\inte(M))$, $C_n (\inte(M))$, and $M$.
Define the set of base points of $C_n \langle M,\partial \rangle$ to be the projection image of $B_n$. Define the set of base points of $C_n '(\inte(M))$ to be the intersection of $B_n$ with $C_n '(\inte(M))$, and define the set of base points of $C_n (\inte(M))$ to be the projection image of the set of base points of  $C_n '(\inte(M))$.
The set of base points of $M$ is defined to be the image of $\iota(\xi_0)$. 

%% file: homotopy_conf_space.tex

This section studies the homotopy groups of $C_n'\langle M,\partial\rangle$. The base points of the homotopy groups are given in Section \ref{subsec_base_point_simply_conn}. 

\subsection{General setups}
We start by setting up some notations and general results. All the results and constructions in this subsection only depend on the assumption that $M$ is a connected compact $4$--manifold with a non-empty boundary.

Since $M$ has a non-empty boundary, there is a section $M\to STM$, which defines a lifting $\mathfrak{s}: C_{n}\langle M,\partial \rangle \to C'_{n}\langle M,\partial \rangle$. The map $\mathfrak{s}$ can be chosen so that it takes the set of base points into the set of base points. 
 Since the projection $C'_{n}\langle M,\partial \rangle \to C_{n}\langle M,\partial \rangle$ is a fiber bundle, we have:
\begin{Lemma}
\label{lem_pi_i_Cn'_split}
    $\pi_iC'_{n}\langle M,\partial \rangle \cong \pi_iC_{n}\langle M,\partial \rangle  \oplus \pi_i((S^3)^{n}).$ \qed 
\end{Lemma}

\begin{remark}
\label{rmk_def_s*}
From now on, we will fix a choice of $\mathfrak{s}$ and identify $\pi_i C_n\langle M,\partial \rangle $ with a subspace of $\pi_i C'_n\langle M,\partial \rangle $ via the map $\mathfrak{s}_*$. Similarly, we identify $\pi_i C_n(\inte(M))$ with a subspace of $\pi_i C_n'(\inte(M))$ via the map $(\mathfrak{s}|_{C_n(\inte(M)})_*$.
\end{remark}

\begin{Definition}
For $1\le i \neq j \le n$, Let $w_{ij}$ be the element in $\pi_3C_n(\inte(M))$ defined by letting the $i^{th}$ point rotate around the $j^{th}$ point in a small coordinate neighborhood. The sign of $w_{ij}$ is chosen so that in a coordinate chart, the vector from the $i^{th}$ point to the $j^{th}$ point defines a map from $S^3$ to $\bR^4\setminus\{0\}$ with degree $1$, where the orientation of the local coordinate chart is given by the orientation of $M$. 
\end{Definition}
We also abuse notation and use $w_{ij}$ to denote its image in $\pi_3C_n\langle M,\partial \rangle$ under the inclusion map, and its image in $\pi_3C_n'\langle M,\partial \rangle$ under the map $\mathfrak{s}_*$. 

\begin{Definition}
Let $\tau_i:C_1'\langle M,\partial \rangle \to C_n'\langle M,\partial \rangle$ be the map in the cosimplicial structure induced by the map from $\{0,1\}$ to $\{0,\dots,n\}$ which takes $(0,1)$ to $(i-1,i)$. In other words, $\tau_i$ is a composition of coface maps characterized by the following properties: (1) the $i^{th}$ point in $\tau_i(x)$ is equal to the unique point in $x\in C_1'\langle M,\partial \rangle $, (2) for $j<i$, the $j^{th}$ point is equal to the fixed starting point on the boundary of $M$, (3) for $j>i$, the $j^{th}$ point is equal to the fixed endpoint on the boundary of $M$.
\end{Definition}

\begin{Lemma}
    The maps $\tau_i$ induce injections on all homotopy groups. 
\end{Lemma}

\begin{proof}
    Consider the forgetful map $p_i: C_n'\langle M,\partial \rangle\to C_1'\langle M,\partial \rangle$ that forgets every point except for the $i^{th}$ point. Then $p_i\circ \tau_i=\id$, so $\tau_i$ induces injections on all homotopy groups.
\end{proof}

\begin{Lemma}
\label{lem_pi1Cn_dir_sum}
    $\pi_1C_n'\langle M,\partial\rangle = \oplus_{i=1}^n (\tau_i)_*\pi_1C_1'\langle M,\partial\rangle.$
\end{Lemma}
\begin{remark}
On the right-hand side, $ (\tau_i)_*\pi_1C_1'\langle M,\partial\rangle $ are subspaces of $\pi_1C_n'\langle M,\partial\rangle $, and the direct sum is internal. 
\end{remark}
\begin{proof}
We have a chain of isomorphisms
\begin{align*}
(\pi_1C_1'\langle M,\partial\rangle)^n & \to (\pi_1C_1\langle M,\partial\rangle)^n \leftarrow (\pi_1C_1(\inte(M))^n \cong \pi_1(\inte(M)^n)\leftarrow
\\
& \leftarrow \pi_1C_n(\inte(M))\to \pi_1C_n\langle M,\partial \rangle \leftarrow \pi_1C_n'\langle M,\partial \rangle.
\end{align*}
 The first arrow is induced by the projection from $C_1'\langle M,\partial \rangle$ to $C_1\langle M,\partial\rangle$; it induces an isomorphism on $\pi_1$ because the map is a fiber bundle with fiber $S^3$. The second arrow is given by the inclusion of $\inte(M)$ in $C_1\langle M,\partial\rangle$, which is a homotopy equivalence by Proposition \ref{prop_homotopy_equiv_Q}. The third relation is an obvious isomorphism. The fourth arrow is given by the canonical map from $C_n(\inte(M))$ to $M^n$ that evaluates the positions of points; it induces an isomorphism on $\pi_1$ because the map is an open embedding and the complement has codimension at least $4=\dim M$. The next arrow is given by the inclusion of $C_n(\inte(M))$ in $C_n\langle M,\partial\rangle$, which is a homotopy equivalence by Proposition \ref{prop_homotopy_equiv_Q}. The final arrow is induced by the projection map, which is a fiber bundle with fiber $(S^3)^n$ and hence induces an isomorphism on $\pi_1$.

It then follows from a straightforward diagram chasing that the composition of the isomorphisms is equal to $\oplus (\tau_i)_*$.
\end{proof}

\begin{Definition}
Suppose $\alpha\in \pi_1 C_1'\langle M,\partial \rangle \cong \pi_1(M)$. We use $t_i^\alpha$ to denote the action of $(\tau_i)_*(\alpha)\in \pi_1C_n'\langle M,\partial \rangle$ on $\pi_iC_n'\langle M,\partial \rangle$. 
\end{Definition}

By Lemma \ref{lem_pi1Cn_dir_sum}, for $i\neq j$ and $\alpha,\beta\in \pi_1(M)$, the actions $t_i^\alpha$ and $t_j^\beta$ are commutative.

For later reference, we record the following property of $\pi_1$ actions on Whitehead products.

\begin{Lemma}
\label{lem_whitehead_under_pi1}
Suppose $X$ is a path connected space and $x_0\in X$ is the base point. Suppose $\gamma$ is a loop based at $x_0$. For $f\in \pi_*(X,x_0)$, we use $f^\gamma$ to denote the action of $\gamma$ on $f$. Then the Whitehead product operation on $\pi_*(X,x_0)$ satisfies $[f,g]^\gamma = [f^\gamma, g^\gamma]$.
\end{Lemma}

\begin{proof}
    Suppose $f$, $g$ are represented by the maps $\hat f:S^i\to X$ and $\hat g:S^j\to X$. By definition, the maps $\hat f$ and $\hat g$ map the base points to $x_0$. Let $\hat f^\gamma$ and $\hat g^\gamma$ be the action of $\gamma$ on $\hat f$ and $\hat g$ as maps. Then $\hat f$ can be freely homotoped to $\hat f^{\gamma}$ such that the trajectory of the base point is $\gamma$, and the similar result holds for $\hat g$ and $\hat g^\gamma$. As a consequence, the map that defines $[f,g]$ can be freely homotoped to a map that defines $[f^\gamma, g^\gamma]$ such that the trajectory of the base point is $\gamma$. Hence the result is proved. 
\end{proof}

\subsection{Properties of $w_{ij}$}
In this subsection, we establish several basic properties of $w_{ij}\in \pi_3C_n(\inte(M))$. All the results in this subsection only depend on the assumption that $M$ is a connected compact $4$--manifold with a non-empty boundary. 

\begin{Lemma}
\label{lem_w_flip_ij}
    Suppose $i\neq j$. Then $w_{ij} = w_{ji}$.
\end{Lemma}
\begin{proof}
    This is because the antipodal map on $S^3$ is homotopic to the identity. 
\end{proof}

\begin{Lemma}
\label{lem_w_ti_action_jk}
Suppose $\alpha\in \pi_1(M)$, and $i,j,k$ are distinct. Then 
    $t_k^\alpha w_{ij} =  w_{ij}$.
\end{Lemma}
\begin{proof}
We may assume $\alpha$
is disjoint from the coordinate
neighborhood where the $i^{th}$ point rotates around the $j^{th}$ point.
Then the lemma is clear.
\end{proof}

\begin{Lemma}
\label{lem_flip_t_action_w}
Suppose $\alpha\in \pi_1(M)$, and $i\neq j$. Then 
    $t_i^\alpha w_{ij} = t_j^{\alpha^{-1}} w_{ij}$. 
\end{Lemma}

\begin{proof}
    One can directly construct a homotopy between the two sides by letting the ball in which the $i^{th}$ point rotates around the $j^{th}$ point move along $\alpha$. 
\end{proof}

\begin{remark}
Let $\beta\in \pi_1(M)$. 
  Applying the action of $t_j^\beta$ to Lemma \ref{lem_flip_t_action_w}, we have
  \begin{equation}
  \label{eqn_ti_tj_wij}
     t_j^\beta  t_i^\alpha w_{ij} =  t_j^{\beta} (t_j^{\alpha^{-1}} w_{ij}) = t_j^{\beta \alpha^{-1}} w_{ij} = t_{i}^{\alpha \beta^{-1}} w_{ij}. 
  \end{equation} 
\end{remark}

\begin{Lemma}
\label{lem_whitehead_zero_three_w}
Suppose $i,j,k$ are distinct. Then
    $[w_{ij}+w_{ik},w_{jk}] = 0$. 
\end{Lemma}

\begin{proof}
    Define a map $f:S^3\times S^3\to C_n'\langle M,\partial \rangle$ as follows: the first factor of $S^3$ parametrizes a rotation of the $j^{th}$ point around the $k^{th}$ point such that their trajectories are included in a closed ball $D$, and the second factor parametrizes a rotation of the $i^{th}$ point around $D$. Then the restriction of $f$ to $ S^3\times \pt$ has the homotopy class $w_{jk}$, and the restriction of $f$ to $\pt\times S^3$ has the homotopy class $w_{ij} + w_{ik}$. Hence the result is proved.
\end{proof}

By Lemmas \ref{lem_whitehead_zero_three_w}, \ref{lem_w_ti_action_jk}, and \ref{lem_whitehead_under_pi1}, we have that for $i,j,k$ distinct and $\alpha,\beta\in \pi_1(M)$, 
\[
[t_i^{\alpha} w_{ik} + t_i^{\alpha}t_j^{\beta} w_{ij} , t_j^\beta w_{jk}]  = 0.
\]
Hence by \eqref{eqn_ti_tj_wij},
\begin{equation}
\label{eqn_linear_rel_wh_prod_0}
[t_i^{\alpha}w_{ik} + t_i^{\alpha \beta^{-1}}w_{ij}, t_j^\beta w_{jk}] =0.
\end{equation}
This implies 
\begin{equation}\label{eqn_linear_rel_wh_prod}
\begin{split}
[t^{\alpha}_{1}w_{13},t^{\beta}_{2}w_{23}]=-[t^{\alpha\beta^{-1}}_{1}w_{12},t^{\beta}_{2}w_{23}],\\
[t^{\alpha}_{1}w_{12},t^{\beta}_{1}w_{13}]=-[t^{\alpha}_{1}w_{12},t^{\alpha^{-1}\beta}_{2}w_{23}].
\end{split}
\end{equation}
The first equation is obtained by taking $(i,j,k)= (1,2,3)$ in \eqref{eqn_linear_rel_wh_prod_0}, and the second equation is obtained by taking $(i,j,k) = (3,2,1)$ and replacing $(\alpha,\beta)$  with $(\beta^{-1},\alpha^{-1})$ in \eqref{eqn_linear_rel_wh_prod_0}.

\begin{remark}
Recall that we also use $w_{ij}$ to denote its images in $\pi_3C_n\langle M,\partial \rangle$ and $\pi_3C_n'\langle M,\partial \rangle$.  Since these are defined to be the images of $w_{ij}\in \pi_3C_n(\inte(M))$ under inclusion maps, the above results extend immediately to $\pi_3C_3\langle M,\partial \rangle$ and  $\pi_3C_3'\langle M,\partial \rangle$.
\end{remark}

We will also need the following formulas, which are straightforward to verify from the definitions.
\begin{equation}
\label{eqn_delta^i_2_on_w12}
\begin{split}
(\delta^0_2)_*(t^{\alpha}_{1}w_{12}) &= t^{\alpha}_{2}w_{23} , \\
(\delta^1_2)_*(t^{\alpha}_{1} w_{12}) &= t^{\alpha}_{1}w_{13}+t^{\alpha}_{2}w_{23},\\ (\delta^2_2)_*(t^{\alpha}_{1}w_{12}) &= t^{\alpha}_{1}w_{12}+t^{\alpha}_{1}w_{13},\\
(\delta^3_2)_*(t^{\alpha}_{1}w_{12})&=t^{\alpha}_{1}w_{12}.
\end{split}
\end{equation}

We will also need the following property of $w_{12}$. 

\begin{Lemma}
\label{lem_w12_ima_delta11_wh_zero}
    Consider the coface map $\delta_1^1:C_1'\langle M,\partial \rangle\to C_2'\langle M,\partial \rangle$. The Whitehead product of $w_{12}$ with the image of 
    \[
    (\delta_1^1)_*:\pi_3 C_1'\langle M,\partial \rangle\to \pi_3 C_2'\langle M,\partial \rangle
    \]
    is zero in $\pi_5^\bQ  C_2'\langle M,\partial \rangle$. 
\end{Lemma}

\begin{proof}
Let $pr:C_1'\langle M,\partial \rangle \to M$ be the projection. 
Suppose $f:S^3\to C_1'\langle M,\partial \rangle$ is a map representing an element in $\pi_3C_1'\langle M,\partial \rangle$. Without loss of generality, assume the image of $f$ is contained in $C_1'(\inte(M))$, which is equal to $STM$. Let $E = (pr\circ f)^* (TM)$ be the pullback of $TM$ by $pr\circ f$, and let $f_E: E\to TM$ be the corresponding bundle map. 

The map $f$ assigns each point $x\in S^3$ a point $pr(f(x))$ in $M$ and a unit tangent vector at $pr(f(x))$, so it defines a section $s:S^3\to E$.  Since $\pi_2(\SO(4))=0$, the bundle $E$ is trivial as an Euclidean vector bundle. 
Let $\varphi: E \to S^3\times \bR^4$ be a trivialization, let $p:S^3\times \bR^4\to \bR^4$ be the projection map. Then $p\circ \varphi\circ s$ is a map from $S^3$ to $S^3$. 

Note that the principal bundle $\SO(3)\hookrightarrow \SO(4)\to \SO(4)/\SO(3) \cong S^3$ is trivial because $\pi_2(\SO(3))= 0$, so the map $p\circ \varphi\circ s$ lifts to $\SO(4)$. As a result, there exists a trivialization $\varphi': E \to S^3\times \bR^4$ such that $\varphi'\circ s$ is a constant section.
Hence the map $s: S^3\to E$ extends to a smooth embedding $\hat s: S^3\times S^3 \to E$ whose image is the unit sphere bundle.  

View the composition $f_E\circ \hat s$ as a map to $C_1'(\inte(M))$, then  $\delta_1^1\circ f_E\circ \hat s$ gives a map from $S^3\times S^3$ to $C_2'\langle M,\partial \rangle$ whose restriction to $S^3\times \pt$ gives $w_{12}$ and whose restriction to $\pt \times S^3$ is $f$. 
\end{proof}

\subsection{Homotopy groups of configuration spaces}
\label{subsec_htpy_conf_space}
In this subsection, we compute $\pi_i^\bQ C_n'\langle M,\partial \rangle$ for $i=3,4,5$ and $n=2,3$ in terms of the homotopy groups of $\hat M$.  These are the computations that we need in order to study $\pi_2\Emb(I,M)$. 
Results in this subsection depend on the assumption that $M=(S^1\times D^3)\natural \hat M$ and $\pi_2^\bQ (\hat M) = 0$.

\subsubsection{The general setup}
By Proposition \ref{prop_homotopy_equiv_Q}, we know that $\pi_i^\bQ C_n'\langle M,\partial \rangle$ is canonically isomorphic to $\pi_i^\bQ C_n'(\inte(M))$. The following lemma allows us to reduce the computation to a computation of the homotopy groups of punctured $M$.

\begin{Lemma}
\label{lem_hmty_grp_Cn'_split}
Consider the fiber bundle   
\[
ST(\inte(M)\setminus\{p_1,\dots,p_n\}) \to C_{n+1}'(\inte(M)) \stackrel{\pi}{\to} C_n'(\inte(M))
\]
where $\pi$ is the forgetful map that removes the last point in the configuration and $p_1,\dots,p_n$ are distinct points in $\inte(M)$. There exists a map $s:  C_n'(\inte(M)) \to C_{n+1}'(\inte(M))$ such that $\pi\circ s \simeq \id$. 
\end{Lemma}

\begin{proof}
    We define the map $s$ directly. Let $N(\partial M)$ be a collar neighborhood of $M$. Let $f:M\to M\setminus N(\partial M)$ be a homeomorphism that is isotopic to the identity when regarded as a map to $M$. Then $f$ induces a homeomorphism $f_*$ from $C_n'(\inte(M))$ to $C_n'(\inte(M\setminus N(\partial M)))$. Fix a point $p\in N(\partial M)$ and fix a unit tangent vector $v$ at $p$. Define $s$ by taking the first $n$ points in $s(x)$ to be given by $f_*(x)$, and taking the last point in $s(x)$ to be given by $(p,v)$. Then the map $s$ satisfies the desired properties. 
\end{proof}

As a result, we have 
\begin{equation}
\label{eqn_hmty_grp_Cn'_split}
\pi_i C_{n+1}'(\inte(M))  \cong \pi_i C_{n}'(\inte(M)) \oplus \pi_i \,ST(\inte(M)\setminus\{p_1,\dots,p_n\}).
\end{equation}
Since $M$ is a connected manifold with a non-empty boundary, the unit sphere bundle $ST(M)$ has a section, and hence 
\begin{equation}
\label{eqn_hmty_grp_ST_split}
  \pi_i ST(\inte(M)\setminus\{p_1,\dots,p_n\}) \cong \pi_i (S^3)\oplus \pi_i(\inte(M)\setminus\{p_1,\dots,p_n\}).
\end{equation}

Equations \eqref{eqn_hmty_grp_Cn'_split} and \eqref{eqn_hmty_grp_ST_split} allow us to reduce the computations of the homotopy groups of $\pi_i C_{n}'(\inte(M))$ to the computations of homotopy groups of $\inte(M)\setminus\{p_1,\dots,p_n\}$. Note that $\inte(M)\setminus\{p_1,\dots,p_n\}$ is homotopy equivalent to $M\vee (\vee_n S^3)$.

\subsubsection{Computation of $\pi_3^\bQ C_n'\langle M,\partial\rangle$ for $n=2,3$}
In order to compute $\pi_3^\bQ C_n'\langle M,\partial\rangle$ for $n=2,3$ using \eqref{eqn_hmty_grp_Cn'_split} and \eqref{eqn_hmty_grp_ST_split}, we need to first compute the homotopy groups of $M\vee S^3$ and $M\vee(\vee_2 S^3)$.

We start by recalling the following result in homotopy theory, which is an immediate consequence of the excision theorem in $\bQ$--coefficients. 

\begin{Lemma}
\label{lem_wedge_X_Y_pi3Q}
    Suppose $X$, $Y$ are simply connected and $\pi_2^\bQ(X) \cong \pi_2^\bQ(Y)=0$. Let $i_X, i_Y$ be the embeddings of $X$ and $Y$ in $X\vee Y$. Suppose $k=3$ or $4$. Then $(i_X)_*$, $(i_Y)_*$ are both injective on $\pi_k^\bQ$, and we have 
    \[
    \pi_k^\bQ(X\vee Y) = (i_X)_*\pi_k^\bQ(X)\oplus (i_Y)_*\pi_k^\bQ(Y),
    \]
    where the right-hand side is an internal direct sum. \qed
\end{Lemma}

We apply Lemma \ref{lem_wedge_X_Y_pi3Q} to compute $\pi_3(M\vee (\vee_n S^3))$.

\begin{Definition}
\label{defn_notation_M_vee_n_S3}
    Let $i_M$ be the inclusion of $M$ in $M\vee (\vee_n S^3)$, and let $i_{S^3}^{(k)}$ be the inclusion of $S^3$ in $M\vee (\vee_n S^3)$ as the $k^{th}$ component in the wedge sum. 
Fix a generator of $\pi_3(S^3)$, and 
let $\eta_k\in \pi_3^\bQ(M\vee (\vee_n S^3))$ denote its image via $i_{S^3}^{(k)}$. Note that there is a canonical isomorphism between $\pi_1(M)$ and $\pi_1(M\vee (\vee_n S^3))$ induced by the inclusion map. For $\alpha\in \pi_1(M)\cong \pi_1(M\vee (\vee_n S^3))$, let $t^\alpha\eta_k$ denote the action of $\alpha$ on $\eta_k$.
\end{Definition}

\begin{Lemma}
\label{lem_pi3Q_M_wedge_nS3}
We have $\eta_k\neq 0 \in \pi_3^\bQ(M\vee (\vee_n S^3))$ for all $k$, and $(i_M)_*$ is injective on $\pi_3^\bQ$. We also have 
\[
\pi_3^\bQ (M\vee (\vee_n S^3)) = (i_M)_*(\pi_3^\bQ M) \, \bigoplus \Big(\bigoplus_{\stackrel{\alpha\in \pi_1(M)}{k=1,\dots, n}} \bQ\cdot (t^\alpha \eta_k)\Big),
\]
where the right-hand side is an internal direct sum.
\end{Lemma}

\begin{proof}
Let $\widetilde{M}$ be the universal cover of $M$. 
Then the universal cover of $M\vee (\vee_n S^3)$ is the wedge sum of $\widetilde{M}$ and a collection of $S^3$'s indexed by the set $\pi_1(M)\times \{1,\dots,n\}$. Therefore the result follows from Lemma \ref{lem_wedge_X_Y_pi3Q}. 
\end{proof}

As a result, we have the following propositions. 

\begin{Proposition}
\label{prop_C2'_pi3}
The following properties hold for $\pi_3^\bQ C_2'\langle M,\partial\rangle$.
\begin{enumerate}
    \item         Each element of the form $t_1^\alpha w_{12}$ $(\alpha\in \pi_1(M))$ is non-zero in $\pi_3^\bQ C_2'\langle M,\partial\rangle$.
    \item We have 
    \begin{equation}
    \label{eqn_C2'_pi3}
    \pi_3^\bQ C_2'\langle M,\partial\rangle = \Big(\bigoplus_{\alpha \in \pi_1(M)} \bQ\cdot (t_1^\alpha w_{12}) \Big) \bigoplus\Big(\bigoplus_{i=1,2}(\tau_i)_* (\pi_3^\bQ C_1'\langle M,\partial\rangle)\Big),
    \end{equation}
    where the right-hand side is an internal direct sum. 
\end{enumerate}
\end{Proposition}

\begin{proof}
    Equations \eqref{eqn_hmty_grp_Cn'_split}, \eqref{eqn_hmty_grp_ST_split}, and Lemma \ref{lem_pi3Q_M_wedge_nS3} yield a direct sum decomposition of $\pi_3^\bQ C_2'\langle M,\partial\rangle$ which is abstractly isomorphic to the right-hand side of \eqref{eqn_C2'_pi3}. By a straightforward diagram chasing, each component of the decomposition is equal to the right-hand side of \eqref{eqn_C2'_pi3} as a subspace of $\pi_3^\bQ C_2'\langle M,\partial\rangle$.
\end{proof} 

The same argument gives the following description of $\pi_3^\bQ C_3'\langle M,\partial \rangle$.

\begin{Proposition}
    \label{prop_C3'_pi3}
    The following properties hold for $\pi_3^\bQ C_3'\langle M,\partial\rangle$.
\begin{enumerate}
    \item         Each element of the form $t_i^\alpha w_{ij}$ for $1\le i<j\le 3$, $\alpha\in \pi_1(M)$ is non-zero in $\pi_3^\bQ C_3'\langle M,\partial\rangle$.
    \item We have \[
    \pi_3^\bQ C_3'\langle M,\partial\rangle = \Big(\bigoplus_{\alpha \in \pi_1(M),i<j} \bQ\cdot (t_i^\alpha w_{ij}) \Big) \bigoplus\Big(\bigoplus_{i=1,2,3}(\tau_i)_* (\pi_3^\bQ C_1'\langle M,\partial\rangle)\Big),
    \]
        where the right-hand side is an internal direct sum. 
\end{enumerate}
\end{Proposition}

\begin{proof}
   The result follows from the same argument as Proposition \ref{prop_C2'_pi3}. The only difference is that we need to apply  \eqref{eqn_hmty_grp_Cn'_split} two times for $n=1$ and $2$ to obtain the desired decomposition of $\pi_3^\bQ C_3'\langle M,\partial\rangle $. 
\end{proof}

\subsubsection{Computation of $\pi_4^\bQ C_n'\langle M,\partial\rangle$}
Note that $\pi_4^\bQ(S^3)=0$. Hence the same argument as before yields the following result. 
\begin{Proposition}
        \label{prop_Cn'_pi4}
We have 
\[
    \pi_4^\bQ C_n'\langle M,\partial\rangle = \bigoplus_{1\le i \le n}(\tau_i)_* (\pi_4^\bQ C_1'\langle M,\partial\rangle),
    \]
        where the right-hand side is an internal direct sum. 
\end{Proposition}

\subsubsection{Computation of $\pi_5^\bQ C_n'\langle M,\partial\rangle$ for $n=2,3$}
To compute $\pi_5^\bQ C_n'\langle M,\partial \rangle$, we recall the following result from homotopy theory.

\begin{Lemma}
\label{lem_wedge_homotopy_group_whitehead}
Suppose $X$ is $n$--connected and $Y$ is $m$--connected. Let $i_X,i_Y$ be the inclusion of $X$ and $Y$ in $X\vee Y$. Then 
\[\pi_{n+m+1}(X\vee Y) = (i_X)_*\pi_{n+m+1}(X) \oplus (i_Y)_*\pi_{n+m+1}(Y) \oplus (\pi_{n+1}(X)\otimes \pi_{m+1}(Y)).\] \qed  
\end{Lemma}
Here, $\pi_{n+1}(X)\otimes \pi_{m+1}(Y)$ maps injectively into $\pi_{n+m+1}(X\vee Y)$ by the Whitehead product and is identified with its image, so the right-hand side is an internal direct sum.

We apply Lemma \ref{lem_wedge_homotopy_group_whitehead} to compute $\pi_5^\bQ(M\vee(\vee_n S^3))$. The following lemma uses the notation given by Definition \ref{defn_notation_M_vee_n_S3}.

\begin{Lemma}
\label{lem_S3_vee_M_pi5}
    Fix a full order ``$<$'' on $\pi_1(M)$. Then 
    \begin{enumerate}
        \item For each $\alpha\in \pi_1(M)$ and $k\in\{1,\dots,n\}$, the map $\pi_3^\bQ(M)\to \pi_5^\bQ(M\vee(\vee_n S^3))$ defined by taking the Whitehead product with $t^\alpha\eta_k$ is injective.
        \item For each $\alpha,\beta\in \pi_1(M)$ and $k,l\in \{1,\dots,n\}$ with $(\alpha,k)\neq (\beta,l)$, we have $\big[t^\alpha \eta_k, t^\beta \eta_l\big] \neq 0 \in \pi_5^\bQ(M\vee(\vee_n S^3))$. 
        \item  We have 
    \begin{multline*}
            \pi_5^\bQ(M\vee(\vee_n S^3)) = (i_M)_*\pi_5^\bQ(M)\bigoplus \Big(\bigoplus_{\stackrel{\alpha\in\pi_1(M)}{1\le k \le n}} [t^\alpha\eta_k,(i_M)_*\pi_3^\bQ(M)] \Big)\\
            \bigoplus \Big(\bigoplus_{\stackrel{\alpha,\beta\in\pi_1(M)}{1\le k<l\le n}}\bQ\cdot [t^\alpha\eta_k, t^\beta\eta_l]\Big)
            \bigoplus \Big(\bigoplus_{\stackrel{\alpha<\beta}{1\le k \le n}}\bQ\cdot [t^\alpha\eta_k, t^\beta\eta_k]\Big),
    \end{multline*}
    where the right-hand side is an internal direct sum.
    \end{enumerate}
\end{Lemma}

\begin{proof}
    Note that $\pi_5^\bQ(S^3)=0$. If $\pi_2(M)=0$, then the result follows from applying Lemma \ref{lem_wedge_homotopy_group_whitehead} to the universal cover of $M\vee(\vee_n S^3)$. In general, we know that $\pi_2^\bQ(M)=0$, so we may take the localization in $\bQ$ on the universal cover of $M\vee(\vee_n S^3)$ and apply Lemma \ref{lem_wedge_homotopy_group_whitehead}. 
\end{proof}

As a consequence, we have the following results.
\begin{Proposition}
\label{prop_C2'_pi5}
The homotopy group $\pi_5^\bQ C_2'\langle M,\partial\rangle$ satisfies the following properties.
\begin{enumerate}
    \item  For each $t_1^\alpha w_{12}$ with $\alpha\in \pi_1(M)$, the map
   \begin{align*}
        \pi_3^\bQ C_1\langle M,\partial \rangle &   \to  \pi_5^\bQ C_2'\langle M,\partial \rangle \\
        x & \mapsto [t_1^\alpha w_{12}, (\tau_{2})_* (\mathfrak{s})_*(x)]
    \end{align*}
is injective, where $\mathfrak{s}_*$ is the embedding of $\pi_3^\bQ C_1\langle M,\partial \rangle$ in $\pi_3^\bQ C_1'\langle M,\partial \rangle$ given by Remark \ref{rmk_def_s*}.
\item For each $\alpha\neq \beta\in \pi_1(M)$, we have $[t_1^\alpha w_{12}, t_1^\beta w_{12}]\neq 0$ in $\pi_5^\bQ C_2'\langle M,\partial\rangle$.
\item Fix a full order ``$<$'' on $\pi_1(M)$. Then $  \pi_5^\bQ C_2'\langle M,\partial\rangle$ is the direct sum of the following subspaces:
\begin{enumerate}
    \item $(\tau_i)_* (\pi_5^\bQ C_1'\langle M,\partial\rangle)$ for $i=1,2$.
    \item $[t_1^\alpha w_{12},(\tau_2)_* (\mathfrak{s})_*(\pi_3^\bQ C_1\langle M,\partial\rangle)] $ for $\alpha\in \pi_1(M)$.
    \item $\bQ\cdot [t_1^\alpha w_{12}, t_1^\beta w_{12}]$ for $\alpha<\beta$. 
\end{enumerate}
\end{enumerate}
\end{Proposition}

\begin{proof}
    The proof follows from the same argument as Proposition \ref{prop_C2'_pi3}. The only difference is that we use Lemma \ref{lem_S3_vee_M_pi5} to find $\pi_5^\bQ(M\vee S^3)$. 
\end{proof}

\begin{Proposition}
\label{prop_C3'_pi5}
The homotopy group $\pi_5^\bQ C_3'\langle M,\partial\rangle$ satisfies the following properties.
\begin{enumerate}
    \item  For each $t_i^\alpha w_{ij}$ with $\alpha\in \pi_1(M)$ and $1\le i<j\le 3$, the map
   \begin{align*}
        \pi_3C_1^\bQ \langle M,\partial \rangle &   \to  \pi_5^\bQ C_3'\langle M,\partial \rangle \\
        x & \mapsto [t_i^\alpha w_{ij}, (\tau_{j})_* (\mathfrak{s})_* (x)]
    \end{align*}
is injective, where $\mathfrak{s}_*$ is the embedding of $\pi_3^\bQ C_1 \langle M,\partial \rangle$ in $\pi_3^\bQ C_1'\langle M,\partial \rangle$ given by Remark \ref{rmk_def_s*}.
\item If $\alpha\neq \beta$ and $1\le i<j\le 3$, we have $[t_i^\alpha w_{ij}, t_i^{\beta} w_{ij}]\neq 0$ in $\pi_5^\bQ C_3'\langle M,\partial\rangle$.
\item If $\alpha,\beta\in \pi_1(M)$, then $[t_1^\alpha w_{13}, t_2^\beta w_{23}]\neq 0$ in $\pi_5^\bQ C_3'\langle M,\partial\rangle$.
\item Fix a full order ``$<$'' on $\pi_1(M)$. Then $  \pi_5^\bQ C_3'\langle M,\partial\rangle$ is the direct sum of the following subspaces:
\begin{enumerate}
    \item $(\tau_i)_* (\pi_5^\bQ C_1'\langle M,\partial\rangle)$ for $i=1,2,3$.
    \item $[t_i^\alpha w_{ij},(\tau_j)_*(\mathfrak{s})_*(\pi_3^\bQ C_1\langle M,\partial\rangle)] $ for $\alpha\in \pi_1(M)$ and $1\le i<j\le 3$.
    \item $\bQ\cdot [t_i^\alpha w_{ij}, t_i^\beta w_{ij}]$ for $\alpha<\beta$  and $1\le i<j\le 3$.
    \item $\bQ\cdot [t_1^\alpha w_{13}, t_2^\beta w_{23}]$ for $\alpha,\beta\in \pi_1(M)$.
\end{enumerate}
\end{enumerate}
\end{Proposition}

\begin{proof}
    The result follows from the same argument as before. 
\end{proof}

Note that by \eqref{eqn_linear_rel_wh_prod}, Part (3) of Proposition \ref{prop_C3'_pi5} implies that $[t_i^{\alpha} w_{ik} ,t_j^{\beta} w_{jk}]\neq 0$ as long as $i,j,k$ are distinct. 
By \eqref{eqn_linear_rel_wh_prod}, if  $[t_1^\alpha w_{13}, t_2^\beta w_{23}]$ are linearly independent for all $\alpha,\beta\in \pi_1(M)$, then  $[t_1^\alpha w_{12}, t_2^\beta w_{23}]$ are linearly independent for all $\alpha,\beta\in \pi_1(M)$, and the two sets span the same linear space.
Hence we have the following corollary.

\begin{Corollary}
\label{cor_dir_sum_decom_pi5_C3'}
    Fix a full order ``$<$'' on $\pi_1(M)$. Then $  \pi_5^\bQ C_3'\langle M,\partial\rangle$ is the direct sum of the following subspaces:
\begin{enumerate}
    \item $(\tau_i)_* (\pi_5^\bQ C_1'\langle M,\partial\rangle)$ for $i=1,2,3$.
    \item $[t_i^\alpha w_{ij},(\tau_j)_*(\mathfrak{s})_*(\pi_3^\bQ C_1 \langle M,\partial\rangle)] $ for $\alpha\in \pi_1(M)$ and $1\le i<j\le 3$.
    \item $\bQ\cdot [t_i^\alpha w_{ij}, t_i^\beta w_{ij}]$ for $\alpha<\beta$  and $1\le i<j\le 3$.
    \item $\bQ\cdot [t_1^\alpha w_{12}, t_2^\beta w_{23}]$ for $\alpha,\beta\in \pi_1(M)$.
\end{enumerate}
\end{Corollary}

\subsection{The diagonal action}
\label{subsec_dia_act}
For the rest of Section \ref{sec_homotopy_conf_space}, we assume that $M=(S^1\times D^3)\natural \hat M$ and $\hat{M}$ satisfies the assumption of Theorem \ref{Thm-main}, which is restated as follows. 
\begin{Assumption}\label{assump pi_2=0} 
$\hat{M}$ is compact, connected, has a non-empty boundary, and satisfies the following conditions:
    \begin{enumerate}
    \item $\pi_{2}^\bQ(\hat{M})=0$. 
    \item Let $b_1,b_2$ be distinct points in the interior of $\hat M$. For every pair of elements $\alpha\in \pi_{1}(\hat{M},b_1)$ and $\beta\in \pi_{k}(\hat{M},b_2)$ with $k=3$ or $4$, there is a map $\gamma:(D^{1},\partial D^{1})\to (\hat{M},b_1)$ that represents $\alpha$, and a map $\mu: (D^{k},\partial D^{k})\to (\hat{M},b_2)$ that represents \emph{a nonzero multiple} of $\beta$, such that \[\gamma(D^1)\cap \mu(D^k)=\emptyset.\]
        \end{enumerate}
\end{Assumption}

Recall from Section \ref{subsec_base_point_simply_conn} that the set of base points of $M$ is given by a straight line in $\{t_0\}\times D^3$.  
Let $\{t_1\}\times D^3$ be a $D^3$--slice in $S^1\times D^3$ such that $t_0\neq t_1$ and that its boundary is disjoint from $\hat M$. Let $M^\vee$ denote the complement of a small tubular neighborhood of $\{t_1\}\times D^3$ in $M$. Then $M^\vee$ deformation retracts to $\hat M$. Define the set of base points on $M^\vee$ to be the same as the set of base points on $M$. Then the embedding of $M^\vee$ to $M$ induces homomorphisms on the homotopy groups with respect to the chosen base point sets. 

\begin{remark}
    Since $\hat M$ is a deformation retract of $M^\vee$, all of their homotopy groups are canonically isomorphic. In the following, we will frequently consider the map $\pi_i (M^\vee) \to \pi_i (M)$. It is therefore more convenient to use $M^\vee$ instead of $\hat M$ because the set of base points of $M^\vee$ can be chosen to be equal to the set of base points of $M$.
\end{remark}

\begin{remark}
By the Seifert-van Kampen theorem, the maps $\pi_1(M^\vee)\to \pi_1(M)$ and $\pi_1(S^1\times D^3)\to \pi_1(M)$ induced by inclusions are both injective. For notational convenience, we will often use $\pi_1(M^\vee)$ and $\pi_1(S^1\times D^3)$ to denote their images in $\pi_1(M)$.
\end{remark}

\begin{Definition}
\label{def_diag_action}
    For $\alpha\in \pi_1(M)$, let $\mathfrak{t}^\alpha$ denote the action by $\prod_{i=1}^n t_i^{\alpha}$ on $\pi_*C_n'\langle M,\partial\rangle$. 
\end{Definition}

We call $\mathfrak{t}^\alpha$ the \emph{diagonal action} by $\alpha$. 

The following two lemmas are consequences of Assumption \ref{assump pi_2=0}.

\begin{Lemma}
\label{lem_diagonal_map_commutes_with_cosimplicial} 
For $k=3,4$ and $n=1,2$, the diagonal action $\mathfrak{t}^\alpha$ commutes with all maps of the form
\[
\pi^{\bQ}_{k} C_n'\langle M,\partial\rangle \to \pi^{\bQ}_{k} C_{n+1}'\langle M, \partial \rangle
\]
and 
\[
\pi^{\bQ}_{k} C_{n+1}'\langle M,\partial\rangle \to \pi^{\bQ}_{k} C_{n}'\langle M, \partial \rangle
\]
induced by the coface and codegeneracy maps.
\end{Lemma}

\begin{proof}
    The fact that $\mathfrak{t}^\alpha$ commutes with all $\sigma_n^i$ and all $\delta^i_n:C_n'\langle M,\partial\rangle \to C_{n+1}'\langle M, \partial \rangle$ for $i=1,\dots,n$ is clear from the definition. We only need to show that $\mathfrak{t}^\alpha$ commutes with  $\delta^0_n$ and $\delta^{n+1}_n$.  
 We focus on the map $\delta^0_n$, the case for $\delta^{n+1}_n$ is similar.

Suppose $\beta\in \pi_k^\bQ C_n'\langle M,\partial\rangle$ with $k=3$ or $4$, $n=1$ or $2$. We show that 
\begin{equation}
\label{eqn_talpha_deta0n_comm}
\delta_n^0\circ \mathfrak{t}^\alpha(\beta) = \mathfrak{t}^\alpha\circ \delta_n^0(\beta).
\end{equation}
 By Propositions \ref{prop_C2'_pi3}, \ref{prop_C3'_pi3}, \ref{prop_Cn'_pi4}, we may assume without loss of generality that either  $\beta$ is in the image of $(\tau_i)_*$ for some $i$, or $n=3$ and $\beta$ has the form $t_i^\eta w_{ij}$ for $i<j$. 

 If $\beta=t_i^\eta w_{ij}$, Equation \eqref{eqn_talpha_deta0n_comm} follows immediately from Equation \eqref{eqn_delta^i_2_on_w12} and Lemmas \ref{lem_w_ti_action_jk}, \ref{lem_flip_t_action_w}. 

 From now, assume $\beta=(\tau_i)_*(\hat \beta)$ for some $1\le i\le n$ and $\hat \beta\in \pi_k^\bQ C_1'\langle M,\partial \rangle$. 
 Since $\pi_k^\bQ C_1'\langle M,\partial \rangle$ is generated by $\pi_k^\bQ C_1'\langle M^\vee ,\partial \rangle$ and its images under $\pi_1(M)$ actions, we may further assume without loss of generality that $\hat \beta\in t_1^{\alpha'} \pi_k^\bQ C_1'\langle M^\vee ,\partial \rangle$ for some $\alpha'\in \pi_1(M)$.  Since $\pi_{1}(M)=\pi_{1}(M^\vee)\ast \pi_{1}(S^{1}\times D^3)$, we only need to consider the cases when $\alpha\in \pi_{1}(M^\vee)$ or $\pi_{1}(S^1\times D^3)$.
 
 Let $q$ denote the projection from $C_1'\langle M,\partial \rangle$ to $M$. Let $b_1, b_2$ be two distinct points in the set of base points of $C'(\inte(M))$. By Assumption \ref{assump pi_2=0} and the assumptions on $\alpha$ and $\hat \beta$, there exist a positive integer $r$, a representative $\mu$ of $r\cdot q_*(\hat \beta)$, and a representative $\gamma$ of $\alpha$, such that the base point of $\gamma$ is $b_1$ and the base point of $\mu$ is $b_2$, and that
 \begin{equation}
 \label{eq: interior don't intersect}
\gamma(D^{1})\cap \mu(D^{k})=\emptyset.    
 \end{equation}

 We show that $\delta_n^0\circ \mathfrak{t}^\alpha(r\beta) = \mathfrak{t}^\alpha\circ \delta_n^0(r\beta).$ Since we are working in $\bQ$ coefficients, this implies the desired equation \eqref{eqn_talpha_deta0n_comm}.

 Let $\hat \mu: S^k \to  C_n'\langle M,\partial \rangle$ be a representative of $r\beta$ such that the trajectory of $x_i$ is given by $\mu$, and the trajectories of $x_j$ are constant for all $j\neq i$.  
 The corresponding representatives of $\delta^0_n\circ \mathfrak{t}^{\alpha}[r\beta]$ and $ \mathfrak{t}^{\alpha}\circ \delta^0_n[r\beta]$, which we denote as $\eta_{1}, \eta_{2}: S^{k}\to C'_{n+1}\langle M,\partial \rangle$ respectively, are only different on the trajectory of the point $x_1$. The trajectory of $x_1$ in  $\eta_{1}$ is the constant map $c:S^{k}\to STM$ with value $(p_0,v_0)$, while the trajectory of $x_1$ in $\eta_{2}$ is given by the $\gamma$--action on $c$, denoted by $\gamma\cdot c$. Here, $\gamma:D^1\to M$ is the aforementioned loop that represents $\alpha$. There is a natural homotopy $h:[0,1]\times S^{k}\to STM$ from  $c$ to $\gamma\cdot c$ whose image is contained in the image of $\gamma$. 
By \eqref{eq: interior don't intersect}, there is a homotopy between $\eta_{1}$ and $\eta_2$ which fixes the trajectories of $x_{2},\cdots, x_{n+1}$ and restricts to $h$ on $x_{1}$.
\end{proof}

The next lemma is proved by a similar argument.
\begin{Lemma}
\label{lem_ti_tauj}
    Suppose $i,j$ are distinct indices in $\{1,\dots,n\}$, and $k=3$ or $4$. Let $(\tau_j)_*$ denote the homomorphism $(\tau_j)_*: \pi_k^\bQ C_1'\langle M,\partial \rangle \to \pi_k^\bQ C_n'\langle M,\partial \rangle$. Let $\alpha \in \pi_1(M)$. Then $t_i^\alpha$ acts trivially on $\Imm(\tau_j)_*$.
\end{Lemma}

\begin{proof}
    Similar to the proof of Lemma \ref{lem_diagonal_map_commutes_with_cosimplicial}, we only need to show that if $\alpha\in \pi_{1}(M^\vee)$ or $\pi_{1}(S^1\times D^3)$ and if $\hat \beta\in t_1^{\alpha'} \pi_k^\bQ C_1'\langle M^\vee ,\partial \rangle$ for some $\alpha'\in \pi_1(M)$, then
    \begin{equation}
    \label{eqn_ti_tauj}
    t_i^\alpha \big((\tau_j)_*(\hat\beta)\big)=(\tau_j)_*(\hat\beta).
    \end{equation}
    Let $q$ denote the projection from $C_1'\langle M,\partial \rangle$ to $M$. 
    By Assumption \ref{assump pi_2=0}, there exist a positive integer $r$, a representative $\mu$ of $r\cdot q_*(\hat \beta)$, and a representative $\gamma$ of $\alpha$, such that \eqref{eq: interior don't intersect} holds. By replacing $\hat \beta$ with $r\hat \beta$, we may assume $r=1$ without loss of generality. Using the above representatives, the two sides of \eqref{eqn_ti_tauj} are only different in the trajectory of $x_i$, and the trajectory of $x_i$ can be homotoped from the left-hand side to the right-hand side without intersecting the trajectory of the other points. Hence the result is proved. 
\end{proof}

\begin{Lemma}
\label{lem_wij_prof_tau_k_zero}
    Suppose $i,j,k$ are distinct indices in $\{1,\dots,n\}$. Let $m=3$ or $4$. Let $(\tau_k)_*$ denote the homomorphism $(\tau_k)_*: \pi_m^\bQ C_1'\langle M,\partial \rangle \to \pi_m^\bQ C_n'\langle M,\partial \rangle$. Then $[t_i^\alpha w_{ij}, \Imm((\tau_k)_*)]=0$ in $\pi_{m+2}^\bQ C_n'\langle M,\partial \rangle$.
\end{Lemma}
\begin{proof}
We only need to prove 
\begin{equation}
\label{eqn_wij_tauj_whitehead_0}
[w_{ij}, \Imm((\tau_k)_*)]=0
\end{equation}
The general result then follows by applying $t_i^\alpha$ on both sides of \eqref{eqn_wij_tauj_whitehead_0} and invoking Lemma \ref{lem_ti_tauj}. 

Let $\beta\in \pi^{\bQ}_{m}C_1'\langle M,\partial \rangle$. Since $M$ has a nonempty boundary, we can find a nonzero $r$ and a representative $\mu: S^m\to M$ of $rq_{*}(\beta)\in \pi_{m} M$ such that $\mu$ is not surjective. Thus we may assume that the trajectory of the $i^{th}$ and $j^{th}$ points in  $w_{ij}$ is contained in the complement of the image of $\mu$ and thus disjoint from the trajectory of the $k^{th}$ point in $(\tau_k)_*(r\beta)$. Hence one can construct a map from $S^3\times S^3$ extending the map from $S^3\vee S^3$ defined by the wedge sum of $w_{ij}$ and $(\tau_k)_*(r\beta)$. 
\end{proof}

\subsection{Cancellations between $(\delta_1^i)_*$}
Note that there is a retraction from $M$ to $S^1\times D^3$ that maps $\hat M$ to an closed ball in $\partial (S^1\times D^3)$. This induces a retraction from $\pi_1(M)$ to $\pi_1(S^1\times D^3)$. We introduce the following definition. 
\begin{Definition}
\begin{enumerate}
    \item Define $N_0^{(2)}\subset \pi_3^\bQ C_2'\langle M,\partial\rangle$ to be the subspace generated by elements of the form $t_1^\alpha w_{12}$ such that $\alpha$ is in the kernel of the retraction from $\pi_1(M)$ to $\pi_1(S^1\times D^3)$.
    \item Define $N_0^{(3)}\subset \pi_3^\bQ C_3'\langle M,\partial\rangle$ to be the subspace generated by elements of the form $t_i^\alpha w_{ij}$ such that $1\le i<j\le 3$ and $\alpha$ is in the kernel of the retraction from $\pi_1(M)$ to $\pi_1(S^1\times D^3)$.
\end{enumerate}
\end{Definition}

We have the following properties.
\begin{Lemma}\label{lem: N invariant under diagonal}
    The spaces $N_0^{(2)}$ and $N_0^{(3)}$ are invariant under the diagonal actions. 
\end{Lemma}
\begin{proof}
    By Lemma \ref{lem_flip_t_action_w}, we have $\mathfrak{t}^\alpha t_i^\beta w_{ij} = t_i^{\alpha\beta\alpha^{-1}} w_{ij}.$ Since the kernel of the retraction from $\pi_1(M)$ to $\pi_1(S^1\times D^3)$ is a normal subgroup of $\pi_1(M)$, the result is proved. 
\end{proof}
\begin{Lemma}
\label{lem_N0(2)_to_N0}
    For all $i$, we have $(\delta^i_2)_*(N_0^{(2)}) \subset N_0^{(3)}$.
\end{Lemma}

\begin{proof}
    The result follows from a direct computation using \eqref{eqn_delta^i_2_on_w12}. 
\end{proof}

The following lemma establishes a cancellation property among $(\delta_1^i)_*$'s on $\pi_3^\bQ$. 
\begin{Lemma}
\label{lem_delta^i_1_cancel_mod_N0(2)}
    Consider the coface maps $\delta^i_1:C_1'\langle M,\partial \rangle \to C_2'\langle M,\partial \rangle$. The image of 
\[
(\delta^0_1)_*-(\delta^1_1)_*+(\delta^2_1)_*: \pi_3^\bQ C_1'\langle M,\partial \rangle \to \pi_3^\bQ C_2'\langle M,\partial \rangle
\]
is contained in $N_0^{(2)}$.
\end{Lemma}
\begin{proof}
By a straightforward diagram chasing using \eqref{eqn_hmty_grp_Cn'_split} and \eqref{eqn_hmty_grp_ST_split}, we know that the image of $\pi_3^\bQ C_2'\langle M^\vee,\partial\rangle$ in $\pi_3^\bQ C_2'\langle M,\partial\rangle$ is the direct sum of the following spaces:
    \begin{enumerate}
        \item $\bQ\cdot (t_1^\alpha w_{12})$, where $\alpha\in \pi_1(M^\vee)$. 
        \item $(\tau_i)_* \pi_3^\bQ  C_1'\langle M^\vee,\partial\rangle$. 
    \end{enumerate}

    Suppose $\mu\in\pi_3^\bQ C_1'\langle M^{\vee},\partial \rangle$. Then $((\delta^0_1)_*-(\delta^1_1)_*+(\delta^2_1)_*)(\mu)$ is in the image of $\pi_3^\bQ C_2'\langle M^\vee,\partial \rangle \to \pi_3^\bQ C_2'\langle M,\partial \rangle $ induced by the inclusion, and the image of $((\delta^0_1)_*-(\delta^1_1)_*+(\delta^2_1)_*)(\mu)$ under $\sigma_2^1$ and $\sigma_2^2$ are both zero. Therefore, it must be contained in the space generated by $\bQ\cdot (t_1^\alpha w_{12})$ with $\alpha\in \pi_1(M^\vee)$, which is a subspace of $N_0^{(2)}$. 

    Since $\pi_2^\bQ(\hat M)=0$, we know that $\pi^{\bQ}_3C_1'\langle M,\partial \rangle $ is generated by $\pi^{\bQ}_3C_1'\langle M^\vee,\partial \rangle $ and its $\pi_1(M)$--actions. Since $(\delta^0_1)_*-(\delta^1_1)_*+(\delta^2_1)_*$ commutes with the diagonal actions and $N_0^{(2)}$ is invariant under diagonal actions, we conclude that the image of $(\delta^0_1)_*-(\delta^1_1)_*+(\delta^2_1)_*$ is contained in $N_0^{(2)}$.
\end{proof}

We use Lemma \ref{lem_delta^i_1_cancel_mod_N0(2)} to prove the following result, which will be used later. 
\begin{Lemma}
\label{lem_wh_tau1_tau2_compare}
Suppose $x\in \pi_3^\bQ C_1'\langle M,\partial \rangle$. Then the following relation holds in $\pi_5^\bQ C_2'\langle M,\partial \rangle$:
\[
[t_1^\alpha w_{12}, (\tau_1)_* x + (\tau_2)_* (t_1^{\alpha^{-1}}x)]\in [t_1^\alpha w_{12}, t_1^\alpha N_0^{(2)}].
\]
\end{Lemma}

\begin{proof}
Let $x' = t_1^{\alpha^{-1}}x$. By Lemma \ref{lem_delta^i_1_cancel_mod_N0(2)}, 
\[
[w_{12}, (\tau_1)_* (x') + (\tau_2)_* (x') - (\delta_1^1)_*(x')] \in [w_{12}, N_0^{(2)}]. 
\]
By Lemma \ref{lem_w12_ima_delta11_wh_zero}, $[w_{12}, (\delta_1^1)_*(x')]=0$. So $[w_{12}, (\tau_1)_* (x') + (\tau_2)_* (x')] \in [w_{12}, N_0^{(2)}]$. Applying $t_1^\alpha$ to both sides of the equation, we have 
\[
[t_1^\alpha w_{12}, t_1^\alpha (\tau_1)_* (x') + t_1^\alpha (\tau_2)_* (x')]\in [t_1^\alpha w_{12}, t_1^\alpha N_0^{(2)}].
\]
Note that $t_1^\alpha (\tau_1)_* (x') = (\tau_1)_* (x)$,  $t_1^\alpha (\tau_2)_* (x') = (\tau_2)_* (x')$. So the desired result is proved. 
\end{proof}

\begin{Corollary}
\label{cor_wh_tau1_tau2_compare}
Suppose $x\in \pi_3^\bQ C_1'\langle M,\partial \rangle$ and $\alpha$ is in the kernel of the retraction from $\pi_1(M)$ to $\pi_1(S^1\times D^3)$. Then 
\[
[t_1^\alpha w_{12}, (\tau_1)_* x + (\tau_2)_* (t_1^{\alpha^{-1}}x)]\in [t_1^\alpha w_{12}, N_0^{(2)}].
\]
\end{Corollary}
\begin{proof}
This follows from Lemma \ref{lem_wh_tau1_tau2_compare} and the observation that $t_1^\alpha N_0^{(2)}= N_0^{(2)}$ when $\alpha$ is in the kernel of the retraction from $\pi_1(M)$ to $\pi_1(S^1\times D^3)$.
\end{proof}

\subsection{Images of the coface maps}

\begin{Definition}
\label{def_N_in_pi5_C3'}
    Define $N\subset \pi_5^\bQ C_3'\langle M,\partial\rangle$ to be the subspace generated by the following:
\begin{enumerate}
    \item $(\tau_i)_* (\pi_5^\bQ C_1'\langle M,\partial\rangle)$ for $i=1,2,3$.
    \item $[t_i^\alpha w_{ij},(\tau_j)_*(\mathfrak{s})_*(\pi_3^\bQ C_1\langle M,\partial\rangle)] $ for $\alpha\in \pi_1(M)$ and $1\le i<j\le 3$.
    \item \label{item_N_wij_component}
    $\bQ\cdot [t_i^\alpha w_{ij}, t_i^\beta w_{ij}]$ for $(i,j)=(1,2)$ or $(2,3)$.  
    \item \label{item_N_w13_component}
    $\bQ\cdot [t_1^\alpha w_{13}, t_1^\beta w_{13}]$ where $\alpha$ or $\beta$ is in $\pi_1(M)\setminus \pi_1(S^1\times D^3)$.   
    \item 
        \label{item_nontrivial_summand_N}
        $\bQ\cdot [t_1^\alpha w_{12}, t_2^\beta w_{23}]$ where $\alpha$ or $\beta$ is in $(\pi_1(M)\setminus \pi_1(S^1\times D^3))\cup\{1\}$.
    \item 
        \label{item_S1xD3_summand_N}
    The images of the 4 composition maps 
    \[
    \pi_5C_2'\langle S^1\times D^3,\partial\rangle \xrightarrow{(\delta^i_2)_*} \pi_5C_3'\langle S^1\times D^3,\partial\rangle \to \pi_5C_3'\langle M,\partial\rangle,
    \]
    where the second map is induced by the inclusion. 
\end{enumerate}
\end{Definition}

The main result of this subsection is the following. 
\begin{Proposition}
\label{prop_coface_C2_C3_in_N}
    Suppose $\hat M$ satisfies Assumption \ref{assump pi_2=0} and $\delta_2^i: C_2'\langle M,\partial \rangle \to C_3'\langle M,\partial \rangle$ is a coface map. Then the image of $(\delta_2^i)_*$ on $\pi_5^\bQ$ is contained in $N$. 
\end{Proposition}

The rest of Section \ref{sec_homotopy_conf_space} is devoted to the proof of Proposition \ref{prop_coface_C2_C3_in_N}.

\begin{Lemma}
\label{lem_N0_prod_in_N}
    $[N_0^{(3)},\pi_3^\bQ C_3'\langle M,\partial\rangle]\subset N$. 
\end{Lemma}
\begin{proof}
    By Propositions \ref{prop_C3'_pi3} and \ref{prop_C3'_pi5},  Equation \eqref{eqn_linear_rel_wh_prod}, and Lemma \ref{lem_wij_prof_tau_k_zero}, we only need to show that if $\alpha$ is in the kernel of the retraction from $\pi_1(M)$ to $\pi_1(S^1\times D^3)$ and $i<j$, then 
    \begin{equation}
    \label{eqn_ti_wij_taui}
    [t_{i}^\alpha w_{ij}, (\tau_i)_*(\pi_3^\bQ C_1'\langle M,\partial \rangle )]\subset N.
    \end{equation}

Suppose $(i,j)=(1,2)$ and $x\in \pi_3^\bQ C_1'(\inte(M))\cong \pi_3^\bQ C_1'\langle M,\partial \rangle$. 
Pushing forward Corollary \ref{cor_wh_tau1_tau2_compare} by the map $\delta_2^3$, we have
    \begin{equation}
    \label{eqn_t1_w12_tau1}
    [t_{1}^\alpha w_{12}, (\tau_1)_*(x)]\in [t_{1}^\alpha w_{12}, N_0^{(3)}].
    \end{equation}
Since we assume $x\in \pi_3^\bQ C_1'(\inte(M))$, and $\pi_5^\bQ C_3'(\inte(M))$ is canonically isomorphic to $\pi_5^\bQ C_3'\langle M,\partial \rangle$, relation \eqref{eqn_t1_w12_tau1} also holds in $\pi_5^\bQ C_3'(\inte(M))$. 
Note that $C_3'(\inte(M))$ admits an involution that switches the second and the third point in the configuration. Applying this involution to \eqref{eqn_t1_w12_tau1}, we have
\[
    [t_{1}^\alpha w_{13}, (\tau_1)_*(x)]\in [t_{1}^\alpha w_{13}, N_0^{(3)}].
\]
Finally, pushing forward Corollary \ref{cor_wh_tau1_tau2_compare} by the map $\delta_2^0$ yields
\[
    [t_{2}^\alpha w_{23}, (\tau_2)_*(x)]\in [t_{2}^\alpha w_{23}, N_0^{(3)}].
\]
Since $[t_i^\alpha w_{ij}, N_0^{(3)}]\subset N$ when $\alpha$ is in the kernel of the retraction from $\pi_1(M)$ to $\pi_1(S^1\times D^3)$, \eqref{eqn_ti_wij_taui} is proved. 
\end{proof}

\begin{remark}
\label{rmk_identify_piiCnMvee_within_piiCnM}
Since there exists a retraction from $M$ to $\hat M$, the inclusion map $M^\vee \to M$ induces injections on all homotopy groups. By a straightforward diagram chasing using \eqref{eqn_hmty_grp_Cn'_split} and \eqref{eqn_hmty_grp_ST_split}, the maps $\pi_i C_n'\langle M^\vee, \partial \rangle \to \pi_i C_n'\langle M, \partial \rangle $ are all injections. In the following, we will use $\pi_iC_n'\langle M^\vee, \partial \rangle$ to denote its image in $\pi_i C_n'\langle M, \partial \rangle$ for notational convenience.
\end{remark}

Now we prove Proposition \ref{prop_coface_C2_C3_in_N}. To better organize the proof, we divide the result into 4 parts.

\begin{Lemma}
\label{lem_im_delta_2^0_in_N}
   Consider the coface map $\delta_2^0:C_2'\langle M,\partial\rangle \to C_3'\langle M,\partial\rangle $. The image of $(\delta_2^0)_*$ on $\pi_5^\bQ$ is contained in $N$.
\end{Lemma}
\begin{proof}
The map $\delta_2^0$ inserts a first point into the configuration at the boundary. Therefore, the image of $t_1^\alpha w_{12}$ is $t_2^\alpha w_{23}$ (see also \eqref{eqn_delta^i_2_on_w12}), and for $i=1,2$, the image of 
\[
(\tau_i)_*(\pi^{\bQ}_kC_1'\langle M,\partial \rangle ) \subset \pi^{\bQ}_k C_2'\langle M,\partial \rangle
\]
under $(\delta^0_{2})_{*}$ is 
\[
(\tau_{i+1})_*(\pi_kC_1'\langle M,\partial \rangle ) \subset \pi^{\bQ}_k C_3'\langle M,\partial \rangle.
\]
Hence the result follows immediately from Proposition \ref{prop_C2'_pi5}.
\end{proof}

\begin{Lemma}
\label{lem_im_delta_2^3_in_N}
   Consider the coface map $\delta_2^3:C_2'\langle M,\partial\rangle \to C_3'\langle M,\partial\rangle $. The image of $(\delta_2^3)_*$ on $\pi_5^\bQ$ is contained in $N$.
\end{Lemma}
\begin{proof}
    This is similar to the previous case. The image of $t_1^\alpha w_{12}$ under the map  $(\delta_2^3)_*$ is $t_1^\alpha w_{12}$, and for $i=1,2$, the image of 
$(\tau_i)_*(\pi^{\bQ}_kC_1'\langle M,\partial \rangle ) \subset \pi^{\bQ}_k C_2'\langle M,\partial \rangle$
is 
$(\tau_{i})_*(\pi_kC_1'\langle M,\partial \rangle ) \subset \pi_k C_3'\langle M,\partial \rangle.$
Hence the result follows from Proposition \ref{prop_C2'_pi5}.
\end{proof}

\begin{Lemma}
\label{lem_im_delta_2^1_in_N}
   Consider the coface map $\delta_2^1:C_2'\langle M,\partial\rangle \to C_3'\langle M,\partial\rangle $. The image of $(\delta_2^1)_*$ on $\pi_5^\bQ$ is contained in $N$.
\end{Lemma}
\begin{proof}
    $\delta_2^1$ is the map that doubles the first point in the configuration, so the image of $t_1^\alpha w_{12}$ is $t_1^\alpha w_{13} + t_2^\alpha w_{23}$. 

    The composition of $\tau_1: \pi^{\bQ}_5 C_1'\langle M,\partial \rangle\to C_2'\langle M,\partial \rangle$ with $\delta_2^1$ is equal to $\delta_2^3\circ \delta_1^1$. So the image of $(\tau_1)_*(\pi_5^{\bQ} C_1'\langle M,\partial \rangle)\subset  \pi^{\bQ}_5 C_2'\langle M,\partial \rangle$ under $(\delta_2^1)_*$ is contained in the image of $(\delta_2^3)_*$, which is contained in $N$ by Lemma \ref{lem_im_delta_2^3_in_N}.

    The image of 
\[(\tau_2)_*(\pi^{\bQ}_5 C_1'\langle M,\partial \rangle ) \subset \pi^{\bQ}_5 C_2'\langle M,\partial \rangle\]
under $(\delta_2^1)_*$ is 
\[(\tau_{3})_*(\pi^{\bQ}_5C_1'\langle M,\partial \rangle ) \subset \pi^{\bQ}_5 C_2'\langle M,\partial \rangle,\] which is contained in $N$.

The image of $[t_1^\alpha w_{12},(\tau_2)_*(\pi_3C_1'\langle M,\partial\rangle)] $ under $(\delta_2^1)_*$ equals 
\[
[t_1^\alpha w_{13}+ t_2^\alpha w_{23}, (\tau_{3})_*(\pi_3C_1'\langle M,\partial \rangle )],
\]
which is contained in $N$.

The image of $[t_1^\alpha w_{12}, t_1^\beta w_{12}]$ under $(\delta_2^1)_*$ equals 
\[
[t_1^\alpha w_{13}+ t_2^\alpha w_{23},t_1^\beta w_{13}+ t_2^\beta w_{23}].
\]
By \eqref{eqn_linear_rel_wh_prod}, 
    \[
    [t_1^{\alpha} w_{13},  t_2^\beta w_{23}] =- [t_1^{\alpha \beta^{-1}} w_{12}, t_2^\beta w_{23}],\quad 
    [t_2^\alpha w_{23}, t_1^\beta w_{13}] =   [t_1^{\beta \alpha^{-1}}w_{12}, t_2^\alpha w_{23}]. 
    \]
    So 
    \begin{align}
    & [t_1^\alpha w_{13}+ t_2^\alpha w_{23},t_1^\beta w_{13}+ t_2^\beta w_{23}] 
    \nonumber
    \\
    = \, & 
    [t_1^\alpha w_{13}, t_1^\beta w_{13}]  + [t_2^\alpha w_{23}, t_2^\beta w_{23}] 
    - [t_1^{\alpha \beta^{-1}} w_{12}, t_2^\beta w_{23}]+ [t_1^{\beta \alpha^{-1}}w_{12}, t_2^\alpha w_{23}]. 
    \label{eqn_simplify_(w13+w23)_product}
    \end{align}
    If either $\alpha \notin \pi_1(S^1\times D^3)$ or $\beta\notin \pi_1(S^1\times D^3)$, then the right-hand side of \eqref{eqn_simplify_(w13+w23)_product} is contained in the subspace spanned by \eqref{item_N_wij_component}, \eqref{item_N_w13_component}, \eqref{item_nontrivial_summand_N} in the definition of $N$ (Definition \ref{def_N_in_pi5_C3'}). If both $\alpha$ and $\beta$ are contained in  $\pi_1(S^1\times D^3)$, then \eqref{eqn_simplify_(w13+w23)_product} is contained in subspace \eqref{item_S1xD3_summand_N} in the definition of $N$. 
\end{proof}

\begin{Lemma}
\label{lem_im_delta_2^2_in_N}
   Consider the coface map $\delta_2^2:C_2'\langle M,\partial\rangle \to C_3'\langle M,\partial\rangle $. The image of $(\delta_2^2)_*$ on $\pi_5^\bQ$ is contained in $N$.
\end{Lemma}

\begin{proof}
    $\delta_2^2$ is the map that doubles the second point in the configuration, so the image of $t_1^\alpha w_{12}$ is $t_1^\alpha w_{12} + t_1^\alpha w_{13}$.

    The image of 
$(\tau_1)_*(\pi^{\bQ}_5 C_1'\langle M,\partial \rangle ) \subset \pi^{\bQ}_5 C_2'\langle M,\partial \rangle$
under $(\delta_2^2)_*$ is 
\[(\tau_{1})_*(\pi^{\bQ}_5 C_1'\langle M,\partial \rangle ) \subset \pi^{\bQ}_5 C_3'\langle M,\partial \rangle,\] which is contained in $N$.

    The composition of $(\tau_2)_*: \pi_5^\bQ C_1'\langle M,\partial \rangle\to \pi_5^\bQ C_2'\langle M,\partial \rangle$ with $\delta_2^2$ is equal to $\delta_2^0\circ \delta_1^1$. So the image of $(\tau_2)_*(\pi_5^\bQ C_1'\langle M,\partial \rangle)\subset  \pi_5^\bQ C_2'\langle M,\partial \rangle$ under $(\delta_2^2)_*$ is contained in the image of $(\delta_2^0)_*$, which is contained in $N$ by Lemma \ref{lem_im_delta_2^0_in_N}.

The image of $[t_1^\alpha w_{12},(\tau_2)_*(\pi_3C_1'\langle M,\partial\rangle)] $ under $(\delta_2^2)_*$ equals 
\[
[t_1^\alpha w_{13}+ t_2^\alpha w_{23}, (\delta_2^0\circ \delta_1^1)_*(\pi_3C_1'\langle M,\partial \rangle )].
\]
By Lemmas \ref{lem_N0(2)_to_N0} and \ref{lem_delta^i_1_cancel_mod_N0(2)}, we have 
\[
 (\delta_2^0\circ \delta_1^1)_*(x) \equiv (\tau_2)_*(x) + (\tau_3)_*(x)\quad \mod N_0^{(3)}
\]
for every $x\in \pi_3C_1'\langle M,\partial \rangle $. Hence by Lemma \ref{lem_N0_prod_in_N}, we know that 
\[
[t_1^\alpha w_{13}+ t_2^\alpha w_{23}, (\delta_2^0\circ \delta_1^1)_*(x)]\equiv [t_1^\alpha w_{13}+ t_2^\alpha w_{23},  (\tau_2)_*(x) + (\tau_3)_*(x)]\quad \mod N
\]
for every $x\in \pi_3C_1'\langle M,\partial \rangle $. By Lemma \ref{lem_wij_prof_tau_k_zero}, $[t_i^\alpha w_{ij}, \Imm((\tau_k)_*)] = 0$ if $i,j,k$ are distinct, so 
\[
[t_1^\alpha w_{13}+ t_2^\alpha w_{23},  (\tau_2)_*(x) + (\tau_3)_*(x)] \equiv [t_2^\alpha w_{23},  (\tau_2)_*(x)] \quad \mod N.
\]
Pushing forward Lemma \ref{lem_wh_tau1_tau2_compare} by $\delta_2^0$, we know that 
\[
[t_2^\alpha w_{23},  (\tau_2)_*(x)] \equiv \big[t_2^\alpha w_{23}, t_2^\alpha \big((\delta_2^0)_* N_0^{(2)}\big)\big] \quad \mod N.
\]
The space $\big[t_2^\alpha w_{23}, t_2^\alpha \big((\delta_2^0)_* N_0^{(2)}\big)\big]$ is generated by elements of the form
$[t_2^\alpha w_{23}, t_2^{\alpha\beta} w_{23}]$, where $\beta$ is in the kernel of the retraction from $\pi_1(M)$ to $\pi_1(S^1\times D^3)$. So it is contained in subspace \eqref{item_N_wij_component} in the definition of $N$. 

The image of $[t_1^\alpha w_{12}, t_1^\beta w_{12}]$ under $(\delta_2^1)_*$ equals 
\[
[t_1^\alpha w_{12}+ t_1^\alpha w_{13},t_1^\beta w_{12}+ t_1^\beta w_{13}].
\]
By \eqref{eqn_linear_rel_wh_prod}, 
    \[
    [t_1^{\alpha} w_{12},  t_2^\beta w_{13}] =- [t_1^{\alpha} w_{12}, t_2^{\alpha^{-1}\beta} w_{23}],\quad 
    [t_1^\alpha w_{13}, t_1^\beta w_{12}] = [t_1^\beta w_{12}, t_2^{\beta^{-1}\alpha}w_{23}].  
    \]
    So 
    \begin{align}
    & [t_1^\alpha w_{12}+ t_1^\alpha w_{13},t_1^\beta w_{12}+ t_2^\beta w_{13}] 
    \nonumber
    \\
    = \, & 
    [t_1^\alpha w_{12}, t_1^\beta w_{12}]  + [t_1^\alpha w_{13}, t_1^\beta w_{13}] 
    -[t_1^{\alpha} w_{12}, t_2^{\alpha^{-1}\beta} w_{23}]+ [t_1^\beta w_{12}, t_2^{\beta^{-1}\alpha}w_{23}]. 
    \label{eqn_simplify_(w12+w13)_product}
    \end{align}
    If $\alpha \notin \pi_1(S^1\times D^3)$ or $\beta\notin \pi_1(S^1\times D^3)$, then the right-hand side of \eqref{eqn_simplify_(w12+w13)_product} is contained in the subspace spanned by \eqref{item_N_wij_component}, \eqref{item_N_w13_component}, \eqref{item_nontrivial_summand_N} in the definition of $N$. If both $\alpha$ and $\beta$ are contained in  $\pi_1(S^1\times D^3)$, then \eqref{eqn_simplify_(w12+w13)_product} is contained in subspace \eqref{item_S1xD3_summand_N} in the definition of $N$. 
\end{proof}

In conclusion, Proposition \ref{prop_coface_C2_C3_in_N} is proved. 

\begin{remark}
    The quotient of $\pi_5^\bQ(C_3'\langle M,\partial M\rangle)$ by $N$ is isomorphic to the quotient of $W$ from \cite[Proposition 3.4]{BG2019} by the subspace generated by $t_1^p t_3^q[w_{12}, w_{23}]$ with $p=0$ or $q=0$. 
\end{remark}

%% file: spectral_sequence.tex

In this section, we use a spectral sequence from fibration towers to generalize Budney--Gabai's W3 invariant to $M=(S^1\times D^3)\natural \hat M$.
\subsection{Spectral sequence from a fibration tower}
\label{subsec_ss_fib_tower}
We first review a general construction of spectral sequences on homotopy groups from fibration towers.
This is a version of the Bousfield--Kan spectral sequence. 
Consider a sequence of maps $f_i:X_i\to X_{i-1}$ that are fibrations. Assume further that each $X_i$ is path connected. Assume that $X_i=\pt$ for $i\le 0$ and that $f_i=\id$ for $i$ sufficiently large. Take a base point $b_i$ in each $X_i$ such that $f_{i}(b_{i})=b_{i-1}$. Let $F_i$ be the path component of $f_i^{-1}(b_{i-1})$ that contains $b_i$.  Assume further that $\pi_1(F_i,b_i)$ are all abelian.

Let $\tilde \pi_1(F_i)$ be the kernel of $\pi_1(F_i)\to \pi_1(X_i)$, and let $\tilde \pi_1^\bQ(F_i) = \tilde \pi_1(F_i) \otimes \bQ$.  For each $i$, there is a long exact sequence:
\begin{equation}\label{eq: exact sequence from fibration}
\cdots\to \pi_n^{\bQ}(X_i)\to  \pi_{n}^{\bQ}(X_{i-1}) \to \pi_{n-1}^{\bQ}(F_i)\to \dots\to \pi_{2}^{\bQ}(X_i)\to \pi_{2}^{\bQ}(X_{i-1}) \to \tilde\pi_1^{\bQ}(F_i)\to 0.
\end{equation}
For notational convenience, let $\tilde\pi_n(F_i)=\pi_n(F_i)$ for $n\ge 2$, and let $\tilde \pi_n^\bQ(F_i) = \tilde \pi_n(F_i) \otimes \bQ$. Define $X_\infty$ to be $X_i$ for $i$ sufficiently large. By summing up the above exact sequences, one obtains an exact couple:
$$
\to \oplus_{i\ge 1,j} \tilde\pi_i^{\bQ}(F_j) \to  \oplus_{i\ge 2,j} \pi_i^{\bQ}(X_j) \to  \oplus_{i\ge 2,j} \pi_i^{\bQ}(X_j) \to, 
$$
which induces a spectral sequence. See e.g., \cite{BottTu}. The first page of this spectral sequence is 
\[E^{1}_{*,*}=\oplus_{i,j }E^{1}_{i,j}=\oplus_{i\ge 1,j}\tilde \pi_i^{\bQ}(F_j),\] 
and the differential on the $k$-th page has the form 
\[d_{i,j}^{k}:E^{k}_{i,j}\to E^{k}_{i-1,j+k}.\]
The spectral sequence converges to a graded associated group of $\bigoplus_{i\ge 2} \pi_i^{\bQ}(X_\infty)$. More precisely, we have $E^{\infty}_{1,j} = 0$, and for all $i\ge 2$,
\[
E^{\infty}_{i,j}= \ker(\pi_{i}^{\bQ}(X_{\infty})\to \pi_{i}^{\bQ}(X_{j-1}))/\ker(\pi_{i}^{\bQ}(X_{\infty})\to \pi_{i}^{\bQ}(X_{j})). 
\]

Now we consider the special case when $f_i=\id$ for all $i\ge 4$. In this case, we have  
$$E^{\infty}_{2,3}= E^{3}_{2,3}= \ker ((f_3)_*:\pi_{2}^{\bQ}(X_3)\to \pi_{2}^{\bQ}(X_2)),$$ 
and
\[
E^{\infty}_{2,3}=(\tilde\pi_{2}^{\bQ}(F_3)/\operatorname{Im}(d^{1}_{3,2}))/\operatorname{Im}d^{2}_{3,1} = (\pi_{2}^{\bQ}(F_3)/\operatorname{Im}(d^{1}_{3,2}))/\operatorname{Im}d^{2}_{3,1}.
\]
The differential $d^{1}_{3,2}: \pi_{3}^{\bQ}(F_{2})\to \pi_{2}^{\bQ}(F_{3})$ is given by the composition 
\begin{equation}\label{eq: d^{1}_{3,2}}
\pi_{3}^{\bQ}(F_2) \to \pi_{3}^{\bQ}(X_2) \to \pi_{2}^{\bQ}(F_3), 
\end{equation}
where the first map is induced by the inclusion and the second map is in the long exact sequence of the fibration $f_3:X_3\to X_2$. 

The term $E^{2}_{3,1}$, which is the domain of the differential $d^{2}_{3,1}$, is the kernel of the composition
$$
\pi_{3}^{\bQ}(F_1) \to \pi_{3}^{\bQ}(X_1) \to \pi_{2}^{\bQ}(F_2).
$$
And the differential \begin{equation}\label{eq: d^{2}_{3,1}}
d^{2}_{3,1}: E^{2}_{3,1}\to E^{2}_{2,3}=\pi_{2}^{\bQ}(F_3)/\operatorname{Im}d^{1}_{3,2}    
\end{equation} is given by the composition
$$E^2_{3,1}\xrightarrow{\text{inclusion}}\pi_{3}^{\bQ}(F_1) \to \pi_{3}^{\bQ}(X_1) \xrightarrow{\text{lift}} \pi_{3}^{\bQ}(X_2) \to \pi_{2}^{\bQ}(F_3)\xrightarrow{\text{projection}}\pi_{2}^{\bQ}(F_3)/\operatorname{Im}d^{1}_{3,2},
$$
where all the unlabeled arrows are from (\ref{eq: exact sequence from fibration}).

\subsection{The spectral sequence on mapping spaces}
\label{sec_ss_map_space}
We apply the spectral sequence in Section \ref{subsec_ss_fib_tower} to the case 
\[
X_{i}=\begin{cases} \pt &\text{for }i<0,\\
\Map_{3,i}(M) &\text{for } 0\leq i\leq 2,\\
\Map_{3}(M) &\text{for } i\ge 3.
\end{cases}
\]
The maps $f_i:X_{i}\to X_{i-1}$ are restrictions of the domain. Let the base point $b_{3}$ be given by the image of $\iota(x_0)$. Let $b_{2},b_{1}$ be the restriction of $b_{3}$ to the respective skeletons. 

\begin{Proposition}\label{prop: d1} (1) For all $0\leq i\leq 3$ and $1\leq n$, we have a natural isomorphism 
\begin{equation}\label{eq: homotopy group of Fi}
\pi_{n}^{\bQ}(F_{i})\cong \oplus_{\sigma}\pi_{n+i}^{\bQ}(C'_{i}\langle M,\partial \rangle),    
\end{equation}
where the direct sum is taken over all injective monotone maps $\sigma$ from $\{0,\dots,i\}$ to $\{0,1,2,3\}$.

(2) Under the isomorphism (\ref{eq: homotopy group of Fi}), the differential 
\[
d^{1}_{n,i}:\pi_{n}^{\bQ}(F_{i})\to \pi_{n-1}^{\bQ}(F_{i+1})
\]
can be expressed as 
\[
\oplus_{\sigma,\tau}d_{\sigma,\tau}: \oplus_{\sigma}\pi_{n+i}^{\bQ}(C'_{i}\langle M,\partial \rangle)\to \oplus_{\tau}\pi_{n+i}^{\bQ}(C'_{i+1}\langle M,\partial \rangle).
\]
Here $d_{\sigma,\tau}$ is defined as follows: If $\Imm(\sigma)\subset \Imm(\tau)$, then
\[
 d_{\sigma,\tau}: C'_{i}\langle M,\partial \rangle\to C'_{i+1}\langle M,\partial \rangle,
\]
is the coface map induced by the unique map $\delta:\{0,\dots,i\}\to \{0,\dots,i+1\}$ such that $\sigma=\tau\circ \delta$. If $\Imm(\sigma)\not\subset \Imm(\tau)$, then $d_{\sigma,\tau}=0$.
\end{Proposition}
\begin{proof} (1) Each $\sigma$ corresponds to an $i$-dimensional face of $\Delta^{3}$, denoted by $\Delta^{\sigma}$. Let $F_{\sigma}$ be the space of $\Delta^i$ structure-preserving maps from $\Delta^{\sigma}$ to $C'_{i}\langle M,\partial\rangle$ which equals $b_{3}$ when restricted to $\partial \Delta^{\sigma}$. Then by definition, we have 
\[
F_{i}=\prod_{\sigma}F_{\sigma}.
\]
On the other hand, each $F_{\sigma}$ is homotopy equivalent to $\Omega^{i}C'_{i}\langle M,\partial\rangle$. So we have a natural isomorphism 
\[
\pi_{n}^{\bQ}(F_{i})\cong \oplus_{\sigma}\pi_{n}^{\bQ}(F_{\sigma})\cong \pi_{n+i}^{\bQ}(C'_{i}\langle M,\partial\rangle).
\]
(2) This follows from a careful unwinding of definitions. We pull back the fibration $F_{i+1}\to X_{i+1}\to X_{i}$ along the natural inclusion $F_{\sigma}\to F_{i}\to X_{i}$. Then use the naturality of the boundary map in the long exact sequence. 
\end{proof}

\subsection{The domain of $d^2_{3,1}$}
\begin{Lemma}\label{lem_face_relation} 
Suppose for some $\alpha_{0},\alpha_{1},\alpha_{2}\in \pi_{4}^{\bQ}(C'_{1}\langle M,\partial\rangle)$, we have 
\begin{equation}\label{eq: face relation}
(\delta_2^0)_*(\alpha_{0})-(\delta_2^1)_*(\alpha_{1})+(\delta_2^2)_*(\alpha_{2})=0\in \pi_{4}^{\bQ}(C'_{2}\langle M,\partial \rangle).    
\end{equation}
Then $\alpha_{0}=\alpha_{1}=\alpha_{2}$.    
\end{Lemma}
\begin{proof}
Consider the codegeneracy map $\sigma_2^1$, which is the map that forgets the first point in the configuration. We have 
\[
0=\sigma_2^1\big((\delta_2^0)_*(\alpha_{0})-(\delta_2^1)_*(\alpha_{1})+(\delta_2^2)_*(\alpha_{2})\big) = \alpha_0 - \alpha_1,
\]
so $\alpha_0= \alpha_1$. Similarly, by applying $\sigma_2^2$ to \eqref{eq: face relation}, we have $\alpha_1=\alpha_2$.
\end{proof}

\begin{Lemma}\label{lem_kerd^1_in_diagonal}
Under the isomorphism 
\[
\pi_{3}^{\bQ}(F_{1})\cong \oplus_{\sigma}\pi_{4}^{\bQ}(C'_{1}\langle M,\partial\rangle)
\]
provided by Proposition \ref{prop: d1},  the kernel of 
\[
d^{1}_{3,1}:\pi_{3}^{\bQ}(F_{1})\to \pi_{3}^{\bQ}(F_{2})
\]
is contained in the diagonal.
\end{Lemma}
\begin{proof}
This follows directly from Proposition \ref{prop: d1} and Lemma \ref{lem_face_relation}.
\end{proof}

By Lemma \ref{lem_kerd^1_in_diagonal}, the domain $E^2_{3,1}$ of $d^2_{3,1}$ is contained in the diagonal. In the following, we will regard $E^2_{3,1}$ as a subspace of $\pi_{4}^{\bQ}(C_1'\langle M,\partial\rangle)$, identified via the diagonal map.

\subsection{Description of the map $d_{3,1}^2$}
\label{subsec_d_31^2_description}
Let 
$\zeta\in E^2_{3,1}\subset \pi_{4}^{\bQ}(C_1'\langle M,\partial\rangle).$ Suppose $\zeta$ is in the image of the homotopy group in $\bZ$--coefficients. In other words, suppose $\zeta$ is represented by a map from $S^4$ to $C_1'\langle M,\partial\rangle$.
We give an explicit description of $d_{3,1}^2(\zeta)$, which is an element in $E_{2,3}^2\cong \pi^{\bQ}_5(C_3'\langle M,\partial \rangle)/\Imm(d_{3,2}^1)$. 

Identify $D^5$ with $D^3\times D^2$, where the first factor $D^3$ is viewed as the cone over the boundary of a tetrahedron $T$. Label the vertices of $T$ by $v_0,\dots, v_3$. For $0\le i<j\le 3$, let $e_{ij}$ be the edge that connects $v_i$ and $v_j$. For each subset $A$ of $\partial T$, let $\Cone(A)$ denote the cone over $A$ as a subset of the cone over $T$. Then $\Cone(e_{ij})$ is homeomorphic to a $2$--ball. 

Let $K = (\cup_{i,j}\Cone(e_{ij})\times D^2 )\cup (\partial (D^3\times D^2))$. Then $K$ is a subset of $D^3\times D^2$. Let $\partial^- K = \partial (D^3\times D^2)$, $\partial ^+ K = \cup_{i} \textrm{Cone}(v_i)\times D^2$, $\partial^* K = \partial^- K \cup \partial^+ K$. Then the space $K\setminus \partial^*K$ is a disjoint union of $6$ open $4$--balls. For $0\le i<j\le 3$, let $K_{ij} = \Cone(e_{ij})\times D^2$. Then $K = \partial^- K\cup (\cup_{i,j} K_{ij})$ and each $K_{ij}$ is homeomorphic to a closed $4$--ball. Let $\partial^- K_{ij} = K_{ij}\cap \partial^- K$, $\partial^+ K_{ij} = K_{ij}\cap \partial^+ K$, $\partial^* K_{ij} = K_{ij}\cap \partial^* K$.

Let $\mu_0(\zeta)$ be a map from $K$ to $ C_3'\langle M \rangle $, such that 
\begin{enumerate}
\item The restriction of $\mu_0(\zeta)$ to $\partial^* K$ is given by a constant map to a base point. 
\item The restriction of $\mu_0(\zeta)$ to $K_{ij}$ represents the image of $\zeta$ under the composition of coface maps $\delta_2^j\circ \delta_1^i: C_1'\langle M,\partial\rangle\to C_3'\langle M,\partial \rangle $. 
\end{enumerate}
The map $\mu_0(\zeta)$ is uniquely defined up to homotopy relative to $\partial^* K$. 

Since $\zeta \in E^2_{3,1} = \ker d^1_{3,1}$, we know that the map $\mu_0(\zeta)$ can be extended to a map $\mu(\zeta):D^3\times D^2\to C_3'\langle M,\partial\rangle$. The extended map $\mu(\zeta)$ is not unique, but the homotopy class of $\mu(\zeta)$ can only change by adding elements in the image of $d^1_{3,2}$. As a consequence, the map $\mu(\zeta)$ gives a uniquely defined element in $E^2_{2,3}$, which is the range of $d^2_{3,1}$. 

It follows from a straightforward diagram chasing that the map $\mu(\zeta)$ represents the element $d^2_{3,1}(\zeta)$.

%% file: pi_2_hat_M_zero.tex

This section proves Theorem \ref{Thm-main}. Assume that $M=(S^1\times D^3)\natural \hat M$ and   $\hat{M}$ satisfies Assumption \ref{assump pi_2=0}. Let $E_{i,j}^k,d_{i,j}^k$ be as in Section \ref{sec_ss_map_space}. We prove Theorem \ref{Thm-main} by studying the images of the differentials $d^{2}_{3,1}$ and $d^{1}_{3,2}$. Recall the definition of the diagonal $\pi_1(M)$--action from Definition \ref{def_diag_action}. By Lemma \ref{lem_diagonal_map_commutes_with_cosimplicial}, the domain and range of $d^{2}_{3,1}$ both admit induced diagonal $\pi_1(M)$--actions.

\begin{Lemma}
\label{lem_d_31^2_equivariant}
  Suppose $\hat M$ satisfies Assumption \ref{assump pi_2=0}. Then  $d_{3,1}^2$ is equivariant with respect to the diagonal action by $\pi_1(M)$. 
\end{Lemma}

\begin{proof}
    We continue to use the notation from Section \ref{subsec_d_31^2_description}, and introduce some additional notation as follows: Suppose $X$ is a compact smooth $m$--dimensional manifold with corners and $A\subset \partial X$ is a finite union of $(m-1)$--dimensional faces of $X$. Also, assume that $X$ is endowed with a Riemannian metric. Let $N_\epsilon(A)$ denote the open neighborhood of $A$ with radius $\epsilon$, and let $d_A$ denote the distance function to $A$. Then for $\epsilon$ sufficiently small, $X\setminus N_\epsilon(A)$ is diffeomorphic to $X$. Let $\varphi(X,A,\epsilon) : X \to X\setminus N_\epsilon(A)$ be a diffeomorphism that is defined when $\epsilon$ is sufficiently small.  We also require that the composition $X \xrightarrow{\varphi(X,A,\epsilon) } X\setminus N_\epsilon(A) \hookrightarrow X$ varies smoothly on $\epsilon$. 
This condition can be achieved when $\epsilon>0$ is sufficiently small. 

Now we describe $\mu_0(\mathfrak{t}^{\alpha} \zeta)$, where $\alpha\in \pi_1(C_1'\langle M,\partial \rangle)$. Let $\gamma$ be a loop in $C_3'\langle M,\partial \rangle$ that represents the diagonal element corresponding to $\alpha$. 

Note that to determine the homotopy class of $\mu(\mathfrak{t}^{\alpha} \zeta)$ relative to $\partial (D^3\times D^2)$ and modulo the image of $d_{3,2}^1$, one only needs to know the homotopy class of $\mu_0(\mathfrak{t}^{\alpha} \zeta)$ relative to $\partial^- K$. Therefore, we only describe $\mu_0(\mathfrak{t}^{\alpha} \zeta)$ up to homotopy relative to $\partial^- K$.

Choose $\epsilon>0$ sufficiently small. 
Let $N_\epsilon(\partial^*K_{ij})$ denote the open $\epsilon$ neighborhood of $\partial^*K_{ij}$ in $K_{ij}$, and define $N_\epsilon(\partial^\pm K_{ij})$ similarly.

By Lemma \ref{lem_diagonal_map_commutes_with_cosimplicial}, the restriction of $\mu_0(\mathfrak{t}^{\alpha} \zeta)$ to $K_{ij}$ is given by $(\delta_2^j\circ \delta_1^i)_*(\mathfrak{t}^{\alpha} \zeta)= \mathfrak{t}^{\alpha} (\delta_2^j\circ \delta_1^i)_*(\zeta)$. 
After a homotopy relative to $\partial^* K$, we may assume that $\mu_0(\mathfrak{t}^{\alpha} \zeta)$ equals $\mu_0(\zeta) \circ \varphi(K_{ij}, \partial^*K_{ij}, \epsilon)^{-1}$ on $K_{ij}\setminus N_\epsilon(\partial^*K_{ij})$, and equals $\gamma \circ d_{\partial^*K_{ij}}$ on $N_\epsilon(\partial^*K_{ij})$, where $\gamma$ is a loop representing $\alpha$. 

By Lemma \ref{lem_free_hmtp_to_pi_1_action} below, we know that the map $\mu_0(\mathfrak{t}^{\alpha}\zeta )$ can be homotoped relative to $\partial^- K$, such that after the homotopy, its restriction to $K_{ij}\setminus N_\epsilon(\partial^-K_{ij})$ is given by $(\delta_2^j\circ \delta_1^i)_*(\zeta)$, and its restriction to $N_\epsilon(\partial^-K_{ij})$ is given by $\gamma \circ d_{\partial^-K_{ij}}$.  As a result, the map $\mu_0(\mathfrak{t}^\alpha \zeta)$ extends to a map on $D^3\times D^2$ that represents the element $\mathfrak{t}^\alpha [\mu(\zeta)]$ in $\pi_5 (C_3'\langle M\rangle)$. Hence the lemma is proved.
\end{proof}

Now we state and prove the aforementioned Lemma \ref{lem_free_hmtp_to_pi_1_action}.

\begin{Lemma}
\label{lem_free_hmtp_to_pi_1_action}
    Suppose $X$ is a topological space, $A$ is a closed subspace of $X$, and suppose that there exists an open neighborhood $U$ of $A$ with a homeomorphism $h:U \to A\times [0,1)$ such that $h(A) = A\times \{0\}$. Then there exists a homeomorphism $\varphi: X\to X\cup_{A\sim A\times\{0\}} A\times[0,1]$ such that $\varphi(a) = (a,1)$ for all $a\in A$. Suppose $(Y,y_0)$ is a topological space with a base point. For each map $f:(X,A)\to (Y,y_0)$ and a loop $\gamma:[0,1]\to Y$ based as $y_0$, we define $f^\gamma$ by
    \[
    f^\gamma(x) = \begin{cases}
        f\circ \varphi(x) \quad\quad \text{ if } \varphi(x)\in X, \\
        \gamma \circ p \circ \varphi(x) \quad \text{ if } \varphi(x)\in A\times [0,1],
    \end{cases}
    \]
    where $p: A\times [0,1]\to [0,1]$ is the projection map. If there exists a map 
    \[
    H: X\times [0,1]\to Y
    \]
    such that $H(x,1) = g(x)$, $H(x,0) = f(x)$ for all $x\in X$, and $H(a,t) = \gamma(t)$ for all $a\in A$, then $g$ is homotopic to $f^\gamma$ relative to $A$.
\end{Lemma}

\begin{proof}
    Consider the map $\iota: X\cup_{A\sim A\times\{0\}}  A\times[0,1]\to X\times [0,1]$ defined by $\iota(x) = (x,0)$ if $x\in X$ and $\iota(a,t) = (a,t)$ for $(a,t)\in A\times [0,1]$. Let $\tau:X\to X\times [0,1]$ be defined by $\tau(x) = (x,1)$. 
    Then $\iota\circ \varphi^{-1}$ is homotopic to $\tau$ relative to $A$. The desired result is proved by composing $H$ with the homotopy from $\tau$ to $\iota\circ \varphi^{-1}$.
\end{proof}

Recall that the space $M^\vee$ and its base points were defined in Section \ref{subsec_dia_act}. Each $\alpha\in \pi_1(M)$ acts on $\pi_4(M)$, and $\alpha \cdot \pi_4(M^\vee)$ is a subgroup of $\pi_4(M)$. The subgroup $\alpha \cdot \pi_4(M^\vee)$ only depends on the image of $\alpha$ in the left coset of $\pi_1(M^\vee)$ in $\pi_1(M)$.
\begin{Lemma}
Suppose $\pi_2^\bQ(\hat M) = 0$. Then $\pi^\bQ_4(M) = \oplus_{[\alpha]} \alpha \cdot \pi^\bQ_4(M^\vee)$, where $[\alpha]$ takes values in the left coset of $\pi_1(M^\vee)$ in $\pi_1(M)$, and the right-hand side is an internal direct sum. 
\end{Lemma}

\begin{proof}
Let $\widetilde{M^\vee}$ be the universal cover of $M^\vee$. 
    The universal cover of $M$ is homotopy equivalent to the wedge sum of a collection of copies of $\widetilde{M^\vee}$, where each copy of $\widetilde{M^\vee}$ in the sum corresponds to an element in the left coset of $\pi_1(M^\vee)$ in $\pi_1(M)$. Since $\pi_1(\widetilde{M^\vee}) \cong \pi_2^\bQ (\widetilde{M^\vee}) = 0$, the $\pi_4^\bQ$ of the wedge sum of copies of  $\widetilde{M^\vee}$ equals the direct sum of $\pi_4^\bQ$ of each copy. So the result is proved.
\end{proof}

Recall that by Remark \ref{rmk_identify_piiCnMvee_within_piiCnM}, we identify  $\pi_iC_n'\langle M^\vee, \partial \rangle$ with its image in $\pi_i C_n'\langle M, \partial \rangle$.
By Proposition \ref{prop_Cn'_pi4} and a straightforward diagram chasing, we have the following result. 
\begin{Lemma}
    \[\pi^\bQ_4 C_1'\langle M,\partial \rangle = \oplus_{[\alpha]}\mathfrak{t}^\alpha \cdot \pi^\bQ_4 C_1'\langle M^\vee,\partial \rangle,\]
    \[
    \pi^\bQ_4C_2'\langle M,\partial \rangle = \oplus_{[\alpha]}\mathfrak{t}^\alpha \cdot \pi^\bQ_4 C_2'\langle M^\vee,\partial\rangle,
    \]
    where $[\alpha]$ takes values in the left coset of $\pi_1(M^\vee)$ in $\pi_1(M)$, and the direct sums on the right-hand side are internal direct sums.  \qed
\end{Lemma}

Recall that $E_{3,1}^2$ is the domain of $d_{3,1}^2$. Consider the spectral sequence from the fibration tower of mapping spaces defined by $M^\vee$, and let $(d_{3,1}^2)^\vee$ be its second-page differential corresponding to $d_{3,1}^2$, and let $(E_{3,1}^2)^\vee$ be the domain of $(d_{3,1}^2)^\vee$. Then by Lemma \ref{lem_diagonal_map_commutes_with_cosimplicial} and Proposition \ref{prop: d1}, we have
\begin{Lemma}
\label{lem_E_31^2_decomp_vee}
    $E_{3,1}^2  =  \oplus_{[\alpha]}\mathfrak{t}^\alpha\cdot(E_{3,1}^2)^\vee,$ where $[\alpha]$ takes values in the left coset of $\pi_1(M^\vee)$ in $\pi_1(M)$,  and the direct sums on the right-hand side is an internal direct sum. \qed
\end{Lemma}
Note that the range of $d_{3,1}^2$ is the quotient space of $\pi_5C_3'\langle M,\partial \rangle$ by the images of $(\delta_2^i)_*$ for $i=0,1,2,3$.  Recall that $N\subset \pi_5C_3'\langle M,\partial \rangle$ was defined in Definition \ref{def_N_in_pi5_C3'}, and recall that by Proposition \ref{prop_coface_C2_C3_in_N}, the images of $(\delta_2^i)_*$  are all contained in $N$. In the following, we identify $(E_{3,1}^2)^\vee$ with its image in $E_{3,1}^2$ by the homomorphism induced by $M^\vee\hookrightarrow M$.

\begin{Corollary}
The image of $d_{3,1}^2$ is contained in $N/\Imm(d_{3,2}^1)$.
\end{Corollary}
\begin{proof}
    By Lemma \ref{lem_d_31^2_equivariant}, $d_{3,1}^2$ is equivariant with respect to the diagonal actions. So, by Lemma \ref{lem_E_31^2_decomp_vee}, the image of $d_{3,1}^2$ is generated by $d_{3,1}^2((E_{3,1}^2)^\vee)$ and the diagonal actions. By the naturality of the spectral sequence with respect to embeddings of manifolds, $d_{3,1}^2((E_{3,1}^2)^\vee)$ is contained in the quotient image of $\pi_5^\bQ C_3'\langle M^\vee,\partial \rangle$ in $E_{2,3}^2\cong \pi_5^\bQ C_3'\langle M,\partial \rangle/\Imm(d_{3,2}^1)$. 
    
    Note that $\mathfrak{t}^{\alpha}t_{i}^{\beta}w_{ij}=t^{\alpha\beta\alpha^{-1}}_{i}w_{ij}.$
    When $\alpha\in \pi_{1}(M)$ and $\beta\in \pi_{1}(M^\vee)$, we have $\alpha\beta\alpha^{-1}\in (\pi_{1}(M)\setminus \pi_{1}(S^1\times D^3) )\cup \{1\}$. Using this, it is straightforward to check that the images of $\pi_5^\bQ C_3'\langle M^\vee,\partial \rangle$ in $\pi_5^\bQ C_3'\langle M,\partial \rangle$  under diagonal actions by $\pi_1(M)$ are always contained in $N$. Hence the result is proved. 
\end{proof}

\begin{proof}[Proof of Theorem \ref{Thm-main}]
By the spectral sequence, the cokernel of $d_{3,1}^2$, which is denoted by $E_{2,3}^3$, is a subgroup of $\pi_2^\bQ \Map_3(M)$. By the previous results, $E_{2,3}^3$ maps to $\pi_2^\bQ C_3'\langle M,\partial \rangle/N$ by a quotient map. 

The only non-zero term on the third page of the spectral sequence associated with $S^1\times D^3$ is on the grading of $E_{2,3}^3$, which is a quotient space of $\pi_5^\bQ C_3'\langle S^1\times D^3,\partial \rangle$. Therefore, the image of $\pi_2^\bQ\Map_3(S^1\times D^3)$ to $\pi_2 ^\bQ\Map_3(M)$ is contained in $E_{2,3}^3$. Budney--Gabai \cite{BG2019} constructed an infinite collection of diffeomorphisms on $S^1\times D^3$ and computed their images in the quotient space of $\pi_5^\bQ C_3'\langle S^1\times D^3,\partial \rangle$. Using their computation, it is straightforward to verify that the images of these elements in $\pi_2^\bQ C_3'\langle M,\partial \rangle /N$ generate a space of infinite rank. So the theorem is proved.
\end{proof}

%% file: I_times_Y.tex

In this section, we study the case when $M=(S^1\times D^3)\natural \hat M$ and $\hat M = I\times Y$, where $Y$ is a compact, connected $3$--manifold with a non-empty boundary. The main result of this section is the proof of Theorem \ref{Thm-main-2}. We will prove the theorem in two steps. First, we show that the image of $\pi_0\Diff(S^1\times D^3,\partial)$ in $\pi_0\Diff(M,\partial)$ has infinite rank. Then, we study the canonical homomorphism from $\pi_0\Diff(M,\partial)$ to $\pi_0\Homeo(M,\partial)$, and show that the image of $\pi_0\Diff(S^1\times D^3,\partial)$ is still of infinite rank after being mapped to  $\pi_0\Homeo(M,\partial)$.

\subsection{Notation and conventions}
Now we study the case when $\hat M = I\times Y$. In this case, we identify $M = S^1\times D^3 \natural \hat M$ with $I\times (S^1\times D^2\natural Y)$. Write $M_0 = S^1\times D^2\natural Y$, then $M= I\times M_0$. 

Let $\widetilde{M_0}$ be the boundary connected sum of $\mathbb{R}\times D^2$ with countably many copies of $Y$. 
Then $\widetilde{M_0}$ is a normal $\bZ$--covering of $M_0$ that unwinds the $S^1$ factor of $S^1\times D^2\subset M_0$. Let $\widetilde{M}= I\times \widetilde{M_0}$, then $\widetilde{M}$ is a normal $\bZ$--covering of $M$.
For $x_1\neq x_2\in \widetilde M$, we say that $x_1$ is \emph{above} $x_2$, if the projections of $x_1$ and $x_2$ to 
$\widetilde M_0$ are the same, and the projection of $x_1$ to $I=[0,1]$ is greater than the projection of $x_2$ to $I$. 
The deck transformation group is $\bZ$, and we denote the action of  $\alpha\in\bZ$ on $\widetilde{M}$ by $t^\alpha$. 
We define
$$
\tilde C_n(M):=\{(x_1,\cdots,x_n)\in \widetilde M^n | x_i\neq t^\alpha x_j ~\text{for all}~\alpha\in \bZ, 
 i\neq j \}.
$$
Notice that $\tilde C_n(M)$ is the $\bZ^n$--covering of $C_n(M)$.

Suppose $i,j$ are distinct indices in $\{1,2,3\}$ and $\alpha\in \bZ$.
Let $\Co^i_j(\alpha)$ be the (closed) submanifold of $\tilde C_3({M})$ consisting of triples $(x_1,x_2,x_3)$ such that the $t^\alpha x_i$ is above $x_j$.
We define a canonical orientation on $\Co^i_j(\alpha)$. Let $\Delta$ be the diagonal of $\widetilde M_0\times \widetilde M_0$, and let $pr_{i,j}$ be the map from $\tilde C_3({M})$ to $\widetilde{M}_0\times \widetilde{M}_0$ given by the projections of the $i$th and $j$th points to $\widetilde M_0$.  Then the manifold $\Co^i_j(\alpha)$ is an open submanifold of the preimage of $\Delta$ under $pr_{i,j}$. 
Hence an orientation on $\Co^i_j(\alpha)$ is canonically induced by a choice of the orientation of $M_0$.

Fix an orientation on $S^5$ and let $\mu:S^5\to C_3(M)$ be a map. Let $\tilde \mu:S^5\to \tilde C_3(M)$ be a lifting of $\mu$. 
For $\alpha\in \mathbb{Z}$, let $\Co_i^j(\alpha,\mu)$ denote $\tilde \mu^{-1}(\Co^i_j(\alpha))\subset S^5$. After a generic perturbation of $\mu$, the set $\Co_i^j(\alpha,\mu)$ is an oriented submanifold of $S^5$ with codimension $3$.

If $A,B$ are disjoint, closed, oriented $2$--dimensional submanifolds of $S^5$, we use $\lk(A,B)$ to denote the linking number of $A$ and $B$. 
More generally, if $A_1,\dots,A_m$, $B_1,\dots,B_n$ are closed oriented $2$--dimensional submanifolds of $S^5$ such that $A_i\cap B_j=\emptyset$ for all $i,j$, and if $a_1,\dots,a_m,b_1,\dots,b_m$ are integers, we use 
$\lk(\sum_i a_i A_i ,\sum_j b_j B_j)$ to denote the value of $\sum_{i,j}a_ib_j\lk(A_i,B_j)$.

\subsection{The linking numbers}
\begin{Lemma}
\label{lem_Theta_well_defined}
	Suppose $\alpha, \beta\in \bZ$. Then for a generic $\mu$, we have 
	\begin{equation}
		\label{eqn_disjoint}
	\begin{split}
	&\Co^1_2(\alpha,\mu)\cap \Co^3_1(\beta-\alpha,\mu) = \emptyset, \\
	&\Co^1_2(\alpha,\mu)\cap  \Co^3_2(\beta,\mu) = \emptyset, \\
	&\Co^1_3( \alpha-\beta,\mu)\cap \Co^3_1(\beta-\alpha,\mu) = \emptyset, \\
	&\Co^1_3( \alpha-\beta,\mu)\cap \Co^3_2(\beta,\mu) = \emptyset.
	\end{split}
	\end{equation}
	Moreover, for generic $\mu$, the linking number  
	\begin{equation}
		\label{eqn_composite_linking}
    \lk\big(\Co^1_2(\alpha,\mu)-\Co^1_3( \alpha-\beta,\mu),\Co^3_1(\beta-\alpha,\mu)- \Co^3_2(\beta,\mu)\big)
	\end{equation} only depends on the homotopy class of $\mu$. 
\end{Lemma}

\begin{proof}
	The first part of the lemma follows from a standard transversality argument. Let $\mu(t)$ be a $1$--parameter family of maps $\mu$ parametrized by $t\in[0,1]$, then generically, $\cup_t \Co^i_j(-,\mu(t))$ are submanifolds of $[0,1]\times S^5$. After a further generic perturbation, we may assume that there are only finitely many $t$ such that \eqref{eqn_disjoint} does not hold with respect to $\mu(t)$, and that the intersections are all transverse in $[0,1]\times S^5$. Suppose $(s,x)\in[0,1]\times  S^5$ is an intersection point and suppose $\mu(s)(x) = (p_1,p_2,p_3)$, then one of the following holds:
	\begin{enumerate}
		\item $t^\alpha(p_1)$ is above $p_2$, and $t^{\beta-\alpha}(p_3)$ is above $p_1$.
		\item $t^{\alpha-\beta}(p_1)$ is above $p_3$, and $t^\beta (p_3)$ is above $p_2$. 
	\end{enumerate}
	In the first case, $x\in S^5$ is in the intersection of $\Co^1_2(\alpha,\mu(s))$, $\Co^3_1(\beta-\alpha,\mu(s))$, and $\Co^3_2(\beta,\mu(s))$. So it contributes to a change by $\pm 1$ in 
	$$
	\lk(\Co^1_2(\alpha,\mu), \Co^3_1(\beta-\alpha,\mu))
\quad \text{and} \quad 
	\lk(\Co^1_2(\alpha,\mu),  \Co^3_2(\beta,\mu)).
	$$
	By a straightforward calculation of orientations, the change in both linking numbers have the same sign.  So the intersection point does not change the linking number in \eqref{eqn_composite_linking}.  The same holds for the second case by a similar argument.
\end{proof}

Let $\alpha,\beta$ be integers. Let $\Theta_{\alpha,\beta}: \pi_5(C_3(M))\to \bZ$ be the map defined by \eqref{eqn_composite_linking}. 
It is clear from the definition that $\Theta_{\alpha,\beta}$ is a group homomorphism. Since $C_3\langle M,\partial\rangle$ is homotopy equivalent to $C_3(M)$, we define $\Theta_{\alpha,\beta}: \pi_5(C_3\langle M, \partial\rangle) \to \bZ$ to be the corresponding homomorphism. We also abuse notation and let $\Theta_{\alpha,\beta}: \pi_5(C_3'\langle M,\partial\rangle) \to \bZ$ be the composition of $\pi_5(C_3'\langle M,\partial\rangle) \to \pi_5(C_3\langle M,\partial\rangle)$ and $\Theta_{\alpha,\beta}$.

Similarly, suppose $\alpha\neq \beta\in \bZ$ and $i\neq j$. Then $\Co^i_j(\alpha)\cap \Co^i_j(\beta)=\emptyset$. So, for generic $\mu$, the linking number
$ \lk\big(\Co^i_j(\alpha,\mu),\Co^i_j(\beta)\big)$
only depends on the homotopy class of $\mu$. Define
\begin{equation}
\label{eqn_wij_linking}
\Theta^{i,j}_{\alpha,\beta}(\mu) = \lk\big(\Co^i_j(\alpha,\mu),\Co^i_j(\beta,\mu)\big).
\end{equation}

For each fixed $\mu$, there are only finitely many $\alpha,\beta$ such that $\Theta_{\alpha,\beta}(\mu)\neq 0$ or $\Theta^{i,j}_{\alpha,\beta}(\mu)\neq 0$ for some $i,j$. 
Let $\bZ^\infty$ denote the direct sum of countably infinite copies of $\bZ$. 

\begin{Definition}
Let $\Theta:\pi_5(C_3'\langle M,\partial\rangle) \to \bZ^\infty$ be the homomorphism whose coordinates are given by $\Theta_{\alpha,\beta}$ for $\alpha,\beta\in \bZ$ and $\Theta_{\alpha,\beta}^{i,j}$ for $i<j,\alpha<\beta$. We will abuse notation and also use $\Theta:\pi_5^\bQ(C_3'\langle M,\partial\rangle) \to \bQ^\infty$  to denote the same map in $\bQ$ coefficients. 
\end{Definition}

\begin{Lemma}
    If $\hat M = D^4$ and $M=S^1\times D^3$, then the map $\Theta: \pi_5^\bQ(C_3'\langle M,\partial\rangle) \to \bQ^\infty$ is an isomorphism.
\end{Lemma}

\begin{proof}
    By a straightforward computation, we see that the maps $\Theta_{\alpha,\beta}$ ($\alpha,\beta\in \bZ$) , $\Theta_{\alpha,\beta}^{i,j}$ $(\alpha<\beta, i<j)$ form a dual basis of $[t_1^\alpha w_{12}, t_2^\beta w_{23}]$, $[t_i^\alpha w_{ij}, t_i^\beta w_{ij}]$ $(\alpha< \beta, i<j)$ on $\pi_5^\bQ(C_3'\langle M,\partial\rangle)$. So the result is proved. 
\end{proof}

\begin{Lemma}\label{lem_image_theta_delta}
	Let $(\delta^i_2)_* :\pi_5^\bQ C_2'\langle M,\partial \rangle\to \pi_5^\bQ C_3'\langle M,\partial \rangle$ be the induced homomorphism by the coface map $\delta_2^i$. Then the image of $\Theta\circ (\delta^i_2)_*$ is contained in the image of the composition
	\begin{equation}
		\label{eqn_comp_map_htpy}
	\pi_5^\bQ C_2'\langle S^1\times D^3 ,\partial \rangle \to \pi_5^\bQ C_2'\langle M, \partial\rangle \xrightarrow{(\delta^i_2)_*} \pi_5^\bQ C_3'\langle M,\partial \rangle \xrightarrow{\Theta}\bQ^\infty
	\end{equation}
\end{Lemma}

\begin{proof}
Since $\pi_5^\bQ(S^3) =0$, we know that $\pi_5^\bQ C_2'\langle M,\partial \rangle$ is isomorphic to
$\pi_5^\bQ C_2\langle M,\partial \rangle$.
Note that by \eqref{eqn_hmty_grp_Cn'_split} and \eqref{eqn_hmty_grp_ST_split}, we have $ \pi_5^\bQ C_2\langle M,\partial \rangle \cong \pi_5^\bQ C_2(\inte(M))$
 is isomorphic to 
 \begin{equation}
 \label{eqn_pi5_C2_inte_decomp}
 \pi_5^\bQ(\inte(M))\oplus \pi_5^\bQ(\inte(M)\setminus\pt),
 \end{equation}
  The first component in \eqref{eqn_pi5_C2_inte_decomp} is realized in $C_2(\inte(M))$ by letting the first point be a fixed point $\tilde p$ and letting the second point move in the complement of a collar neighborhood $N(\partial M)$ of $\partial M$, where $\tilde p \in N(\partial M)$. The second factor is realized by letting the second point be fixed in $\inte(M)$ and letting the first point move in the complement of the second point.

If $\mu:S^5\to C_2(\inte(M))$ represents an element in the first component of \eqref{eqn_pi5_C2_inte_decomp}, then $\Theta\circ (\delta^i_2)_*([\mu]) = 0$. 

 If $\mu:S^5\to C_2(\inte(M))$ represents an element in the second component of \eqref{eqn_pi5_C2_inte_decomp}:
 Let $N_\epsilon(\partial (S^1\times D^3))$, $N_\epsilon(\partial M)$  be the $\epsilon$--neighborhoods of $\partial (S^1\times D^3)$, $\partial M$ in $(S^1\times D^3)$, $M$ respectively. 
 Choose $\epsilon>0$ sufficiently small so that $N_\epsilon(-)$ above are collar neighborhoods. 
  Let $M_\epsilon = M\setminus N_\epsilon(\partial M)$ and $S_\epsilon = S^1\times D^3\setminus N_\epsilon(\partial (S^1\times D^3)).$
  
 After homotopy, we may assume that under the map $\mu$, the trajectory of $x_2$ is equal to a fixed point $\tilde p\in \inte(S_\epsilon)$ and the trajectory of $x_1$ is induced by a map $\mu_1:S^5 \to M_\epsilon\setminus \{\tilde p\}$.

Let $\varphi: M_\epsilon \to S_\epsilon$ be a retraction that collapses $M_\epsilon \setminus S_\epsilon$ to a disk on the boundary of $S_\epsilon$. 
Let $\varphi_*\mu: S^5\to C_2(\inte(M))$ be the map such that the trajectory of $x_2$ is constant at $\tilde p$, and the trajectory of $x_1$ is given by $\varphi\circ \mu_1$. 

Take a non-vanishing unit tangent vector field that is parallel to the $M_0$ factor in the decomposition $M=I\times M_0$. Lift $\mu$ and $\varphi_*\mu$ to $C_2'(\inte(M))$ via this vector field. 
Then we have $\Co^u_v(\alpha,\delta^i_2(\mu)) = \Co^u_v(\alpha, \delta^i_2(\varphi_*\mu))$ as oriented submanifolds of $S^5$ for all $u,v,\alpha$. Therefore, $\Theta( \delta^i_2\circ \mu) = \Theta( \delta^i_2 \circ \varphi_* \mu)$.
\end{proof}

\begin{Lemma}
	\label{lem_domain_d_Y_times_I}
	The domain of the differential $d_{3,1}^2$ is equal to $\pi_4^\bQ C_1'\langle M,\partial \rangle$.
\end{Lemma}
\begin{proof}
We need to show that $(\delta^0_1)_*+(\delta^2_1)_*=(\delta^1_1)_*$ as maps from $\pi_4^\bQ C_1'\langle M,\partial \rangle$. 

Every element in $\pi_4 C_1'\langle M,\partial \rangle$ is represented by a map from $S^4$ whose image is contained in $C_1'(\pt\times M_0)$. Let $\eta$ be such a map.
Let $\tilde\eta$ be a translation of $\eta$ in the $I$ direction so that it is contained in a different horizontal slice of $M$. We require that $\tilde\eta$ is above $\eta$.
If the base points are placed suitably, the map representing $[\delta^0_1\circ \eta]+[\delta^2_1\circ \eta]$ can be homotoped to a map $S^4\to C_2'(M)$ such that the trajectory of $x_1$ is given by $\eta$ and the trajectory of $x_2$ is given by $\tilde \eta$ where the two points move simultaneously. This map represents $(\delta^1_1)_*([\mu])$.  
\end{proof}

\begin{Lemma}\label{lem_d2_theta}
	Suppose $\mu$ is in the image of the higher differential 
    \[
    d^2_{3,1}: \pi_4^\bQ C_1'\langle M,\partial \rangle \to \pi_5^\bQ C_3'\langle M,\partial \rangle /\sum_i 
	\Imm(\delta^i_2)_*.
    \]
    Then there exists a lift of $\mu$ in $\pi_5^\bQ C_3'\langle M,\partial \rangle$ such that
	$\Theta(\mu) = 0.$
\end{Lemma}

\begin{proof}We choose the boundary conditions for $\Emb(I,M)$ such that the starting point projects to $0\in I=[0,1]$ and the end point projects to $1\in I=[0,1]$. 

Without loss of generality, assume $\eta \in\pi_4 C_1'\langle M,\partial \rangle $ is represented by a map from $S^4$ to $C_1'(M)$ such that the associated unit vector always points upward. 
Let $\mu$ be the map representing $d^2_{3,1}([\eta])$ defined as in Section \ref{subsec_d_31^2_description}. In the construction of $\mu$, we need to choose a homotopy from the sum of 
$\delta^0_1\circ \eta$ and $\delta^2_1\circ \eta$ to $\delta^1_1\circ \eta$, and we use the homotopy constructed in Lemma \ref{lem_domain_d_Y_times_I}.

It is then straightforward to verify that after a small perturbation of $\mu$, the following holds on the image of $\mu$: the $I$--coordinate of $x_3$ is always greater than the $I$--coordinate of $x_2$, and the $I$--coordinate of  $x_2$ is always greater than the $I$--coordinate of $x_1$. Therefore 
$\Co^i_j(\alpha,\mu)$ is always empty when $i<j$, so $\Theta(\mu)=0$. 
\end{proof}

\subsection{The image in the diffeomorphism group}
In this subsection, we prove the following weaker form of Theorem \ref{Thm-main-2}.

\begin{Theorem}
\label{thm_YxI_smooth_ver}
Under the assumptions of Theorem \ref{Thm-main-2}, the image of the map
$$
\pi_0\Diff(S^1\times D^3,\partial)\to \pi_0\Diff(M,\partial)
$$
is of infinite rank.
\end{Theorem}

\begin{proof}
We will consider the spectral sequence associated with both $M= (S^1\times D^3)\natural \hat M$ and $M = S^1\times D^3$. To differentiate the two spectral sequences, we will use $E_{i,j}^k(M)$ to denote the corresponding term in the spectral sequence associated with $M$. When it is not clear from the context, we will also use $d_{i,j}^k(M)$ to denote the differentials in the spectral sequence associated with $M$.

Since $\pi_4(S^1\times D^3)=0$, we have 
\begin{equation}
\label{eqn_E3page_S1xD3}
E^3_{2,3}(S^1\times D^3)=E^2_{2,3}(S^1\times D^3)=\pi^\bQ_5 C_3'\langle S^1\times D^3,\partial\rangle/R
\end{equation}
 where 
$R=\Imm (d^1_{3,2}(S^1\times D^3))$. According to \cite[Proposition 3.4]{BG2019}, the quotient
\[
\Theta(\pi^\bQ_5 C_3'\langle S^1\times D^3,\partial \rangle)/\Theta(R)\cong \pi^\bQ_5 C_3'\langle S^1\times D^3,\partial \rangle/R
\]
has infinite rank. The embedding of $S^1\times D^3$
in $M$ induces a homomorphism between their spectral sequences, whose first pages give the following commutative diagram

\begin{equation*}
\begin{tikzcd}
\pi^\bQ_5C_2' \langle S^1\times D^3,\partial \rangle^{\oplus 4} \arrow[d] \arrow[r,"d^1_{3,2}"]&  \pi^\bQ_5C_3' \langle S^1\times D^3,\partial \rangle \arrow[d] \arrow[r,"\Theta"]
 & \bQ^\infty \arrow[d,"="]\\
\pi^\bQ_5 C_2'\langle M,\partial \rangle^{\oplus 4} \arrow[r,"d^1_{3,2}"] & \pi^\bQ_5 C_3'\langle M,\partial \rangle \arrow[r,"\Theta"]  & \bQ^\infty 
\end{tikzcd}
\end{equation*}

Recall that $E^3_{2,3}(M)$ is a quotient group of $\pi^\bQ_5 C_3'\langle M,\partial \rangle$.
According to Lemma \ref{lem_d2_theta}, the map $\Theta$ induces a map from $E^3_{2,3}(M)$ to  
\[
\Theta(\pi^\bQ_5 C_3'\langle M,\partial \rangle)/\Theta(\Imm d^1_{3,2}(M)),
\]
which we still denote by $\Theta$.

Consider the following commutative diagram:
\begin{equation*}
\begin{tikzcd}
\pi_0\Diff(S^1\times D^3, \partial) \arrow[r,"(\Psi_3)_\ast\circ \cS"]
\arrow[rd,"(\Psi_3)_\ast\circ \cS\circ E"]
&  E^3_{2,3}(S^1\times D^3) \arrow[d] \arrow[r,"\Theta"]
 &  \Theta(\pi^\bQ_5 C_3'\langle S^1\times D^3,\partial \rangle)/\Theta(\Imm d^1_{3,2}(S^1\times D^3))\arrow[d]\\
~  & E^3_{2,3}(M)  \arrow[r,"\Theta"]  &  \Theta(\pi^\bQ_5 C_3'\langle M,\partial \rangle)/\Theta(\Imm d^1_{3,2}(M)) ,
\end{tikzcd}
\end{equation*}
where $E:\pi_0\Diff(S^1\times D^3,\partial)\to \pi_0\Diff(M,\partial)$ denotes the map induced from the extension by identity. By \cite[Theorem 8.3]{BG2019}, the image of $\Diff(S^1\times D^3,\partial)$ in 
\[
\Theta(\pi^\bQ_5 C_3'\langle S^1\times D^3,\partial \rangle )/\Theta(\Imm d^1_{3,2}(S^1\times D^3)) \cong \pi^\bQ_5 C_3'\langle S^1\times D^3,\partial \rangle )/\Imm d^1_{3,2}(S^1\times D^3)
\]
is of infinite rank. Lemma \ref{lem_image_theta_delta} implies that the map
$$
\Theta(\pi^\bQ_5 C_3'\langle S^1\times D^3,\partial \rangle )/\Theta(\Imm d^1_{3,2}(S^1\times D^3))\to \Theta(\pi^\bQ_5 C_3'\langle M, \partial \rangle)/\Theta(\Imm d^1_{3,2}(M))
$$
is an injection. Hence the result is proved. 
\end{proof}

\subsection{The image in the homeomorphism group}
Now we prove Theorem \ref{Thm-main-2}, which is restated as follows.
\begin{repTheorem}{Thm-main-2}
Suppose $Y$ is a connected compact $3$--manifold with a non-empty boundary. Let $\hat{M}=I\times Y$ and let  $M = (S^1\times D^3)\natural \hat M$.  The image of the map
$$
\pi_0\Diff(S^1\times D^3,\partial)\to \pi_0\Homeo(M,\partial)
$$
is of infinite rank.
\end{repTheorem}

We need to introduce some notation.
Let $\Diff_0(M,\partial)$, $\Homeo_0(M,\partial)$ be the group of diffeomorphisms and homeomorphisms that are homotopic to the identity relative to $\partial M$.
Let $\iota(\xi_0) \in \Emb(I,M)$ be the base point defined in Section \ref{subsec_S}. 
Let $\Emb_0(I,M)\subset \Emb(I,M)$ be the connected component of $\Emb(I,M)$ containing $\iota(\xi_0)$. 
Let $\Emb^{top}(I,M)$ be the space of \emph{topological} embeddings of $I$ in $M$, endowed with the compact open topology; let $\Emb_0^{top}(I,M)$ be its connected component containing $\iota(\xi_0)$.

Then the scanning map construction in Section \ref{subsec_S} defines a set-theoretic map from $\pi_0 \Homeo_0(M,\partial)$ to $\pi_2\Emb^{top}_0(I,M)$, which we denote by $\cS^\tau$.

Let $x_0=\tilde p_0,x_{n+1}=\tilde p_1$ be two fixed points on the boundary of $\widetilde M$
which lift the boundary points $p_0,p_1$ used in the definition of $\Emb(I,M)$.
Following the strategy of \cite[Section 4]{BG2023}, we define
$$
C^\tau_n(M)=\{(x_1,\cdots,x_n)\in \widetilde M^n | x_i\neq t^\alpha x_j ~\text{for all}~\alpha\in \bZ\setminus \{0\}, 
 0\le i<j\le n+1 \}.
$$

Similar to $C_n'\langle M,\partial \rangle $, the sequence of spaces $C^\tau_n(M)$ can be made into a cosimplicial space where 
the coface map $\delta_n^i: C^\tau_n(M) \to C^\tau_{n+1}(M)$ ($0\le i \le n+1$) is defined by doubling the $i^{th}$
point
and the codegeneracy map $\sigma_n^i: 
C^\tau_n(M) \to C^\tau_{n-1}(M)$ ($1\le i \le n$) is defined by forgetting the $i^{th}$ point. Different from 
$C_n'\langle M,\partial \rangle $, this time
we do not need to use any compactification or add a tangent vector over each point in the configuration.
The coface maps are still cofibrations since the domain is a closed submanifold of the codomain.

As in Section \ref{sec_mapping_space}, we can also define the $\Delta^n$ structure on $C^\tau_n(M)$ using
its cosimplicial structure. Similarly, define the space
$\Map^\tau_n(M)$ to be the $\Delta^n$ structure-preserving maps from $\widetilde{C}_n'\langle I,\partial\rangle$ to $C^\tau_n(M)$.
The spectral sequence in Section \ref{sec_spec_seq} can then be constructed for the mapping space $\Map^\tau_3(M)$. We use $E^{k,\tau}_{i,j}$ to denote the terms in this spectral sequence and let $d_{i,j}^{k,\tau}$ be the corresponding differentials.

Define a map $\Psi_n^\tau: \Emb^{top}_0(I,M)\to \Map^\tau_n(M)$ as follows. For $f\in \Emb^{top}_0(I,M)$, let $\tilde{f}: I\to \widetilde{M}$ be the lifting of $f$ with respect to a fixed choice of $\tilde{f}(0)$, let $\Psi_n^\tau(f)$ be the composition $\widetilde{C}'_n\langle I,\partial \rangle \to I^n \to \widetilde{M}^n$, where the first map is the canonical projection and the second map is $\tilde{f}^n$.

We also define a map $\Phi_n: \Map_n(M)\to \Map_n^\tau(M)$ as follows. Suppose $f:\widetilde{C}'_n\langle I,\partial \rangle \to C_n'\langle M,\partial \rangle$ is a $\Delta^n$ structure-preserving map, let $pr:C_n'\langle M,\partial \rangle\to M^n$ be the canonical projection, define $\Phi_n(f)$ to be the lifting of $pr\circ f$ to $\widetilde{M}^n$ with respect to a fixed base point. 

Then we have a commutative diagram:
\begin{equation*}
\begin{tikzcd}
\Diff_0(M,\partial) \arrow[d] \arrow[r,"\cS"]& \pi_2\Emb_0(I,M)  \arrow[d] \arrow[r,"(\Psi_3)_*"]& \pi_2\Map_3(M)  \arrow[d,"(\Phi_3)_*"]\\
\Homeo_0(M,\partial) \arrow[r,"\cS^\tau"] & \pi_2\Emb_0^{top}(I,M) \arrow[r,"(\Psi_3^\tau)_*"]  & \pi_2\Map_3^\tau(M)  
\end{tikzcd}
\end{equation*}

The map $\Phi_3$ extends to maps from the fibration tower of $\Map_3(M)$ to the fibration tower of $\Map_3^\tau(M)$, and hence induces a homomorphism between the associated spectral sequences. We use $(\Phi_3)_*$ to denote this homomorphism.

 The term $E^1_{i,j}$ on the first page of the spectral sequence associated with $\Map_3(M)$ is given by direct sums of groups of the form $\pi_{i+j} C_j'\langle M,\partial\rangle $ with $i,j\ge 1$. 
 The term $E^{1,\tau}_{i,j}$ on the first page of the spectral sequence associated with $\Map_3^\tau (M)$ is given by direct sums of groups of the form $\pi_{i+j} C_j^\tau(M)$ with $i,j\ge 1$.

\begin{Lemma}
\label{lem_Phi_surj_fst_page}
    The spectral sequence homomorphism $(\Phi_3)_*$ satisfies 
    the following:
    \begin{itemize}
   \item[(a)] $(\Phi_3)_*:E^1_{3,2}\to E^{1,\tau}_{3,2}$ 
    is surjective;
    \item[(b)] 
     $(\Phi_3)_*: E^1_{3,1}\to E^{1,\tau}_{3,1}$ is an isomorphism; 

     \item[(c)] $(\Phi_3)_*: E^2_{3,1}\to E^{2,\tau}_{3,1}$
     is an isomorphism.
     \end{itemize}
\end{Lemma}

\begin{proof}
 Define
$$
\tilde C_n(\inte(M)):=\{(x_1,\cdots,x_n)\in \inte(\widetilde M)^n | x_i\neq t^\alpha x_j ~\text{for all}~\alpha\in \bZ, 
 i\neq j \}.
$$
Then $\tilde C_n(\inte(M))$ is a covering space of $C_n(\inte(M))$, so $\pi_k C_n\langle M,\partial \rangle \cong \pi_k \tilde C_n(\inte(M))$ for $k\ge 2$. 
Note that $\tilde C_n(\inte(M))$ is a subspace of 
$C^\tau_n(\inte M)\simeq C^\tau_n(M)$. 
Consider the commutative diagram of fibration sequences 
\begin{equation*}
\begin{tikzcd}
\inte \widetilde  M\setminus \{t^\alpha q|\alpha\in \bZ\} \arrow[r] \arrow[d,"h"] 
&  \tilde C_{2}(\inte M)  \arrow[r,"\sigma_2^2"] \arrow[d] &
\inte \widetilde   M \arrow[d]\\
\inte \widetilde M\setminus \{t^\alpha q|\alpha\in \bZ\setminus \{0\}\} \arrow[r] &   C^\tau_{2}(\inte M) \arrow[r,"\sigma_2^2"] &
\inte \widetilde  M,
\end{tikzcd}
\end{equation*}
where $q$ is a fixed point in $\inte(\widetilde{M})$. 
Up to homotopy equivalence, $\widetilde M\setminus \{t^\alpha q|\alpha\in \bZ\}$ can be 
viewed as the wedge sum of countable copies of $\hat M$ and $S^3$ and the map $h$ is obtained by collapsing
one copy of $S^3$. Thus $h$ induces an epimorphism on homotopy groups since its codomain is a retract of its domain.
The long exact sequences of homotopy groups for both fibrations split since both maps $\sigma_2^2$ have
global sections. Therefore the inclusion
$\tilde C_2(\inte(M)) \hookrightarrow C^\tau_2(\inte M)$ induces epimorphisms
on homotopy groups.

For $k\ge 2,j=1,2$, let $s_{k,j}$ be the homomorphism from 
$\pi^\bQ_{k} C_j'\langle M,\partial\rangle $ to $\pi^\bQ_{k} 
C_j^\tau (M)$ defined by the composition 
\begin{equation}
\label{eqn_Phi_n_component}
\pi^\bQ_{k} C_j'\langle M,\partial\rangle \to \pi^\bQ_{k} C_j\langle M,\partial\rangle \cong \pi^\bQ_{k} C_j(\inte(M))\cong \pi^\bQ_{k} \tilde C_j (\inte(M))\to \pi^\bQ_{k} C_j^\tau (M),
\end{equation}
where the first map is induced by the projection from  $C_j'\langle M,\partial\rangle$  to $ C_j\langle M,\partial\rangle $, the second map is induced by the homotopy equivalence between $C_j(\inte(M))$ and $C_j\langle M,\partial \rangle$, the third map is induced by the covering map from $\tilde C_j (\inte(M))$ to $C_j(\inte(M))$, and the last map is induced by the inclusion $\tilde C_j(\inte(M))\hookrightarrow C^\tau_j(M)$. By Lemma \ref{lem_pi_i_Cn'_split}, the first map in \eqref{eqn_Phi_n_component} is surjective.
Hence by the previous results, we know that $s_{k,2}$ is a surjection
and $s_{k,1}$ is an isomorphism for $k\neq 3$.

A straightforward diagram chasing shows that $s_{k,j}$ are the coordinate maps of $E^1_{k-j,j}\to E^{1,\tau}_{k-j,j}$. 
Taking $(k,j)=(5,2)$ and $(4,1)$, we obtain (a) and (b). 

By Lemma \ref{lem_domain_d_Y_times_I}, 
$E^2_{3,1}\cong \pi_4^\bQ C_1'\langle M,\partial \rangle\cong 
\pi_4^\bQ M$, 
which consists of the diagonal elements in 
$E^1_{3,1}\cong (\pi_4^\bQ M)^{\oplus 6}$.
Since the proof of Lemma \ref{lem_kerd^1_in_diagonal} only uses
the cosimplicial structure of $C_n'\langle M,\partial \rangle$,
it also works for the spectral sequence for
$C_n^\tau(M)$. Hence $E^{2,\tau}_{3,1}$ is a subspace
contained in the diagonal of 
$E^{2,\tau}_{3,1}\cong (\pi_4^\bQ M)^{\oplus 6}$. Now the natural 
homomorphism between the spectral sequences implies 
$E^{2,\tau}_{3,1}$ is also isomorphic to the diagonal and
(c) holds.
\end{proof}

\begin{proof}[Proof of Theorem \ref{Thm-main-2}]
    Note that when $\alpha\neq \beta$, $\alpha\neq 0$, $\beta\neq 0$, the definition of the terms in \eqref{eqn_composite_linking} and Lemma \ref{lem_Theta_well_defined} apply verbatim to maps from $S^5$ to $C_3^\tau(M)$. Hence \eqref{eqn_composite_linking} defines a homomorphism from $E^{1,\tau}_{2,3} \cong \pi_5 C_3^\tau(M)$ to $\bZ$. We denote this homomorphism by $\Theta_{\alpha,\beta}^\tau$. 
    Similarly, if $i<j$, and $\alpha<\beta$, $\alpha\neq 0$, $\beta\neq 0$, define $\Theta_{\alpha,\beta}^{\tau,i,j}$ from $\pi_5 C_3^\tau(M)$ to $\bZ$ using the linking number \eqref{eqn_wij_linking}.
    
    Then we have 
    $\Theta_{\alpha,\beta} = \Theta_{\alpha,\beta}^\tau \circ (\Phi_3)_*$ and $\Theta_{\alpha,\beta}^{i,j} = \Theta_{\alpha,\beta}^{\tau,i,j} \circ (\Phi_3)_*$ as maps from $E_{2,3}^1$.

    Define 
    \[
    \Theta' :  \pi_5^\bQ C_3'\langle M,\partial \rangle  \to \bQ^\infty,
    \]
    to be the map whose coordinates are given by $\Theta_{\alpha,\beta}$ where $\alpha\neq \beta$, $\alpha\neq 0$, $\beta\neq 0$, and $\Theta_{\alpha,\beta}^{i,j}$ where $\alpha< \beta$, $\alpha\neq 0$, $\beta\neq 0$, and $i<j$. Define 
    \[
    \Theta^{\tau} :  \pi_5^\bQ C_3^\tau(M)  \to \bQ^\infty,
    \]
    to be the maps whose coordinates are given by $\Theta^\tau_{\alpha,\beta}$ where $\alpha\neq \beta$, $\alpha\neq 0$, $\beta\neq 0$, and $\Theta_{\alpha,\beta}^{\tau,i,j}$ where $\alpha< \beta$, $\alpha\neq 0$, $\beta\neq 0$, and $i<j$. By 
    \cite[Corollary 8.2]{BG2019} and
    \cite[Section 4]{BG2023}, the image of $\pi_0\Diff(S^1\times D^3,\partial)$ in 
\begin{align*}
&\Theta'(\pi^\bQ_5(C_3\langle S^1\times D^3,\partial\rangle))/\Theta'(\Imm(d_{3,2}^1)+ \Imm(d_{3,1}^{2})) 
\\
= & \Theta'(\pi^\bQ_5(C_3\langle S^1\times D^3,\partial\rangle))/\Theta'(\Imm(d_{3,2}^1))
\end{align*}
is of infinite rank.

Therefore, by a straightforward diagram chasing similar to the proof of Theorem \ref{thm_YxI_smooth_ver}, we conclude that the image of $\pi_0\Diff_0(S^1\times D^3,\partial)$ in $\Theta^\tau(E_{2,3}^{1,\tau})/\Theta^\tau(\Imm(d_{3,2}^{1,\tau}) + \Imm(d_{3,1}^{2,\tau}))$ is of infinite rank. 
\end{proof}

%% file: appendix.tex

Let $M$ be a smooth manifold with boundary. This section studies the simplicial compactifications of configuration spaces of points on $M$. 
Assume that $M$ is smoothly and properly embedded in an ambient Euclidean space $\mathbb{R}^m$. We also assume that $M$ is contained in a larger manifold $M'\subset \bR^m$, such that $M'$ is a submanifold of $\bR^m$ without boundary, and that $M'$ intersects $\{0\}\times \bR^{m-1}$ transversely with $M = M'\cap ([0,+\infty)\times \bR^{m-1})$. 

We follow the notation in \cite{Sinha-cpt}.
For $n\in \mathbb{N}$, let $\underline{n}$ denote the set $\{1,\dots,n\}$. For each set $X$, let $C_n(X)$ denote the subset of $X^n$ such that the $n$ coordinates are mutually distinct. If $X$ is a manifold, then $C_n(X)$ is also a manifold. For $S$ a set, let $X^S$ be the product of copies of $X$ such that the factors are labeled by $S$. If $\{X_s\}_{s\in S}$ is a collection of sets labeled by the elements of $S$, we use $(X_s)^S$ to denote the product $\prod_{s\in S} X_s$.

Let $I=[0,+\infty]$ be the union of $[0,+\infty)$ and $\{+\infty\}$, and the smooth structure on $I$ is defined so that the map from $[0,1]$ to $I$ that sends $x$ to $x/(1-x)$ is a diffeomorphism.

Let $A_n[M]= M^{\underline{n}}\times (S^{m-1})^{C_2(\underline{n})}\times I^{C_3(\underline{n})}$, $A_n\langle M\rangle = M^{\underline{n}}\times (S^{m-1})^{C_2(\underline{n})}$. By definition, $C_n(M)\subset M^{\underline{n}}$ is the configuration space of $n$ ordered distinct points in $M$. 
Following \cite[Definition 1.3]{Sinha-cpt}, let $\alpha_n: C_n(M)\to A_n[M]$ be the map such that for $x=(x_1,\dots,x_n)$, the projection of $\alpha_n(x)$ to $M^{\underline{n}}$ is $x$, the projection of $\alpha_n(x)$ to the factor of $S^{m-1}$ labeled by $(i,j)\in C_2(\underline{n})$ equals  $x_{i}-x_{j}/|x_{i}-x_{j}|$, and the projection of $\alpha_n(x)$ to the factor of $I$ labeled by $(i,j,k)\in C_3(\underline{n})$ equals  $|x_i-x_j|/|x_i-x_k|$. Here, $x_i-x_j$ is defined in the ambient space $\bR^m$. Let $\beta_n:C_n(M)\to A_n\langle M \rangle$ be the composition of $\alpha_n$ with the coordinate projection from $A_n[M]$ to $A_n\langle M\rangle$.

Let $C_n[M]$ be the closure of $\alpha_n(C_n(M))$ in $A_n[M]$, let $C_n\langle M\rangle$ be the closure of $\beta_n(C_n(M))$ in $A_n\langle M \rangle$. Then the projection map from $A_n[M]$ to $A_n\langle M\rangle$ restricts to a quotient map from $C_n[M]$ to $C_n\langle M\rangle$, and we denote the quotient map by $Q_n$. 

Following \cite{Sinha-cosimplicial}, there is a variation of the compactification of configuration spaces when $\partial M\neq \emptyset$.

\begin{Definition}[{\cite[Definition 4.10]{Sinha-cosimplicial}}]
\label{def_compactify_w_boundary}
Suppose $\partial M \neq \emptyset$. Fix $p \neq q\in \partial M$. Let $C_n(M,\partial)$ be the subspace of $C_{n+2}(M)$ such that the first point equals $p$ and the last point equals $q$. Let $C_n[M,\partial]$ be the closure of $\alpha_{n+2}(C_n(M,\partial))$ in $A_{n+2}[M]$.  Let $C_n\langle M,\partial \rangle$ be the closure of $\beta_{n+2}(C_n(M,\partial))$ in $A_{n+2}\langle M\rangle$.
\end{Definition}

Let $\inte(M)$ denote the interior of $M$. 
Let $i: C_n(\inte(M)) \to C_n[M,\partial]$ be the inclusion map, and let $Q_n^\partial :C_n[ M,\partial] \to C_n\langle M,\partial \rangle$ be the quotient map. The main result of the appendix is the following.

\begin{Proposition}\label{prop_homotopy_equiv_Q}
	The maps $i$ and $Q_n^\partial $ are homotopy equivalences.
\end{Proposition}

This result is probably well-known to experts, but we were unable to find a reference for this particular statement, so we include a proof here for the sake of completeness. Although in this paper, we only need Proposition \ref{prop_homotopy_equiv_Q} when $M$ is compact, we prove Proposition \ref{prop_homotopy_equiv_Q} for general $M$.

\subsection{The inclusion map $i$}
\label{subsec_inclusion_map}
We first show that the inclusion $i: C_n(\inte(M)) \to C_n[M,\partial]$ is a homotopy equivalence. This is proved by showing that $C_n[M,\partial]$ is a smooth manifold with corners whose interior is the image of $i$.

Recall that $M$ is included in a larger manifold $M'$ such that $M'$ has no boundary, $M'$ is transverse to $\{0\}\times \bR^{m-1}$, and $M = M'\cap ([0,+\infty)\times \mathbb{R}^{m-1})$.

Note that $C_n[M,\partial]$ is a subset of $C_{n+2}[M']$, and both are subsets of $A_{n+2}[M']$. 
Following the convention of \cite{Sinha-cosimplicial}, when referring to a point in $C_n[M,\partial]$ as a configuration of $(n+2)$ points, we label the points by $x_0,\dots,x_{n+1}$.

Let $\iota_i:C_{n+2}(M')\to M'$ $(i=0,\dots,n+1)$ be the map that evaluates the position of the $i^{th}$ point. Let $\pi_{ij}:C_{n+2}(M')\to S^{m-1}$ $(i\neq j)$ be the map that evaluates the unit vector in $\bR^m$ pointing from the $i^{th}$ point to the $j^{th}$ point, let $\pi_{ij}^1$ be the first coordinate of $\pi_{ij}$ viewed as a map to $\bR^m$. 
Then the maps $\iota_i$ and $\pi_{ij}$ extend to smooth maps on $C_{n+2}[M']$.

Let $p,q$ be the fixed points on $\partial M$ as in Definition \ref{def_compactify_w_boundary}. 
The set $C_n[M,\partial]$ consists of elements $x\in C_{n+2}[M']$ such that
\begin{enumerate}
    \item $\iota_0(x)=p$, $\iota_{n+1}(x) = q$.
    \item $\pi_{ij}^1(x)\ge 0$ for $i=1,\dots,n$, $j=0,{n+1}$.
\end{enumerate}

By \cite[Theorem 4.4]{Sinha-cpt}, the space $C_{n+2}[M']$ is a smooth manifold with corners. The strata of $C_{n+2}[M']$ are labeled by f-trees (see \cite[Section 3]{Sinha-cpt}), and we use $C_T[M']$ to denote the closure of the stratum labeled by the tree $T$. By \cite[Theorem 4.15]{Sinha-cpt}, $C_T[M']$ is also a manifold with corners.

The point $(p,q)\in M'\times M'$ is a regular value of $\iota_0\times \iota_{n+1}$ on every $C_T[M']$. Therefore, $R:=(\iota_0\times \iota_{n+1})^{-1}(p,q)$ is a smooth manifold with corners. Let $R_T = R\cap C_T[M']$, we have $R_T$ is a smooth manifold with corners.  If $d$ is the dimension of $M$, then the codimension of $R_T$ in $C_T[M']$ is $2d$ when it is non-empty. 

For $i=1,\dots, n$, let 
\[
S_i = \{x\in R| \pi_{i0}^1(x)\ge 0, \pi_{i,n+1}^1(x)\ge 0\},
\]
\[
\partial S_i = \{x\in R| \pi_{i0}^1(x)= 0, \pi_{i,n+1}^1(x)= 0\}.
\]

Then for each $i$ and each tree $T$, the intersection $\partial S_i\cap R_T$ is a codimension--$1$ submanifold of $R_T$. Hence each $S_i$ is a submanifold with corners in $R$ with codimension zero. 

If $i,j,k$ are three leaf vertices of an f-tree $T$, we say that $(i,j),k$ satisfies the \emph{exclusion relation associated with $T$}, if there exists an internal vertex $v$ of $T$ such that $i,j$ are over $v$ but $k$ is not over $k$. 
For each f-tree $T$ and $i,j\in \{1,\dots,n\}$, we write $i\sim_T j$ if and only if $(i,j),0$ and $(i,j),n+1$ both satisfy the exclusion relation associated with $T$. Then $\sim_T$ is an equivalence relation on $\{1,\dots,n\}$. Let $\hat{V}(T)$ be the set of equivalence classes on $\{1,\dots,n\}$ with respect to $\sim_T$.
For each $i=1,\dots, n$, let $\hat{v}(i)\in \hat{V}(T)$ be the equivalence class represented by $i$.

For $i,j\in \{1,\dots,n\}$, we have $\partial S_i\cap R_T = \partial S_j \cap R_T$ if and only if $i \sim_T j$. Suppose $\hat{v}\in \hat{V}(T)$, let $S_{\hat{v},T}$ denote $S_i\cap R_T$ where $i$ is chosen so that $\hat{v}(i) = \hat{v}$, then the manifolds $\{S_{\hat{v},T}\}_{\hat{v}\in \hat{V}(T)}$ intersect transversely in $R_T$.  

Since $C_n\langle M,\partial \rangle = \cap_{i=1}^{n} S_i$, the above discussions imply that $C_n\langle M, \partial \rangle$ is a smooth manifold with corners. 
It is also clear that $C_n(\inte(M))$ is its interior. Therefore, the inclusion map $i: C_n(\inte(M)) \to C_n[M,\partial]$ is a homotopy equivalence.

\subsection{A general criterion for homotopy equivalence}
Now we prove that $Q_n^\partial$ is a homotopy equivalence. The first step is to establish the following result.

\begin{Proposition}
\label{prop_quotient_with_contractible_fiber}
    Suppose $q:X\to Y$ is a quotient map such that $X$ is compact and metrizable and $Y$ is Hausdorff. Suppose for each $y\in Y$, the pre-image $q^{-1}(\{y\})$ is homeomorphic to a closed ball. Let $L_i$ be the subset of $Y$ consisting of $y\in Y$ such that $q^{-1}(\{y\})$ is a ball of dimension at least $i$. Let $J_i = L_i\setminus L_{i+1}$. Assume that 
    \begin{enumerate}
        \item There exists $N$ such that $L_N = \emptyset$.
        \label{item_LN}
        \item Each $L_i$ is a closed subset of $Y$.
        \label{item_closed}
        \item The inclusion of $q^{-1}(L_i)$ in $X$ is a cofibration for each $i$.
        \label{item_cofib}
        \item Each $J_i$ is homeomorphic to a CW complex.
        \label{item_Ji_CW}
        \item The map $q|_{q^{-1}(J_i)}:q^{-1}(J_i)\to J_i$ is a topological fiber bundle for each $i$.
        \label{item_bundle}
    \end{enumerate}
    Then $q$ is a homotopy equivalence.
\end{Proposition}

\begin{proof}
    For each $i$, there is unique a factorization of $q$ as
    \[
    q: X \stackrel{r_i}{\to} X_i \stackrel{s_i}{\to} Y,
    \]
    such that both $r_i$ and $s_i$ are quotient maps, the map $s_i$ is a homeomorphism on $s_i^{-1}(L_i)$, and the map $r_i$ is a homeomorphism on $ q^{-1}(Y\setminus L_i)$. In other words, $r_i$ is the quotient map that collapses $q^{-1}(L_i)$ to a homeomorphic copy of $L_i$. For each $i$, there is a unique quotient map 
    \[q_i:X_{i+1}\to X_{i}\] 
    such that $r_{i} = q_i\circ r_{i+1}$. 

    Since the push-outs of cofibrations are cofibrations, by Assumption (\ref{item_cofib}), the inclusions of $s_k^{-1}(L_i)$ in $X_k$ are cofibrations, and the inclusions of $s_k^{-1}(L_i)$ in $s_k^{-1}(L_j)$ are cofibrations for $i>j$.

    \begin{Lemma}
    \label{lem_Xi_metrizable}
        Each $X_i$ is compact and metrizable. 
    \end{Lemma}
    \begin{proof}[Proof of Lemma \ref{lem_Xi_metrizable}]
    It is clear that $X_i$ is compact because it is a quotient space of a compact space.
    Let $\Delta_Y$ be the diagonal in $Y\times Y$, and let $\Delta_X$ be the diagonal in $X\times X$. 
    By Assumption (\ref{item_closed}) and the assumption that $Y$ is Hausdorff, the set $\Delta_Y\cap L_i\times L_i$ is closed in $Y\times Y$. Let $R = \{(x_1,x_2)\in X\times X|q_i(x_1)=q_i(x_2)\}$. Then 
    \[
    R = (q\times q)^{-1}(\Delta_Y\cap L_i\times L_i)\cup \Delta_X,
    \]
    therefore $R$ is a closed subset of $X\times X$. Since $X$ is compact Hausdorff and $R$ is closed, the quotient space $X_i$ is Hausdorff (\cite[Proposition 3.57]{lee2010introduction}).  If the quotient space of a compact metrizable space is Hausdorff, then the quotient space is metrizable (\cite[Theorem 4.2.13]{engelking1989general}). So $X_i$ is metrizable.
    \end{proof}
    
    We show that the quotient map $q_i:X_i\to X_{i-1}$ is a homotopy equivalence for each $i>0$. Note that by the definition of $q_i$, we have $s_{i+1}^{-1}(J_i)$ is homeomorphic to $q^{-1}(J_i)$ and $s_{i}^{-1}(J_i)$ is homeomorphic to $J_i$. By Assumption (\ref{item_bundle}), the map 
    \begin{equation}
    \label{eqn_qi_on_Ji}
        q_i|_{s_{i+1}^{-1}(J_i)}: s_{i+1}^{-1}(J_i)\to s_{i}^{-1}(J_i)
    \end{equation}
    is a topological fiber bundle with fiber homeomorphic to the closed ball $D^i$. By Assumption (\ref{item_Ji_CW}), the bundle has a section 
    \[
        \lambda_i: s_{i}^{-1}(J_i) \to s_{i+1}^{-1}(J_i).
    \]
    For each $p\in s_i^{-1}(L_{i+1})$, the pre-image $q_i^{-1}(\{p\})$ contains exactly one point.
    Extend $\lambda_i$ to a map $\bar{\lambda}_i: s_i^{-1}(L_i)\to s_{i+1}^{-1}(L_i)$ by 
    \[
        \bar{\lambda}_i(p) = 
        \begin{cases}
            \lambda_i(p) \quad \text{ if } p\in s_i^{-1}(J_i), \\
            q_i^{-1}(p) \quad \text{ if } p \in s_i^{-1}(L_{i+1}),
        \end{cases}
    \]
    where $q_i^{-1}(p)$ denotes the unique element in $q_i^{-1}(\{p\})$. 

\begin{Lemma}
\label{lem_lambda_bar_continuous}
      The map $\bar{\lambda}_i$ is continuous. 
\end{Lemma}
\begin{proof}[Proof of Lemma \ref{lem_lambda_bar_continuous}]
Suppose $\{p_k\}_{k\in \mathbb{N}}$ is a sequence in $s_i^{-1}(L_i)$ that converges to $p$, we show that $\{\bar{\lambda}_i(p_k)\}$ converges to $\bar{\lambda}_i(p)$ as $k\to\infty$. By Assumption (\ref{item_closed}), $s_i^{-1}(J_i)$ is an open subset of $s_i^{-1}(L_i)$. So, if $p\in s_i^{-1}(J_i)$, then the convergence follows from the continuity of $\lambda_i$. If $p\in s_i^{-1}(L_{i+1})$, we use an argument by contradiction. Assume  $\{\bar{\lambda}_i(p_k)\}$ does not converge to $\bar{\lambda}_i(p)$. Then since $X_{i+1}$ is a compact metrizable space, there is a subsequence of $\{\bar{\lambda}_i(p_k)\}$ that converges to a point $x\neq \bar{\lambda}_i(p)$ in $X_{i+1}$. We still use $\{\bar{\lambda}_i(p_k)\}$ to denote the subsequence. Since $q_i\circ \bar{\lambda}_i = \id$, we know that $q_i(x)$ equals the limit of $q_i\circ \bar{\lambda}_i(p_k)=p_k$ as $k\to \infty$. Since $p_k$ converges to $p$, we have $q_i(x) = p$. Since $p\in L_{i+1}$, we must have $x=\bar{\lambda}_i(p)$, which yields a contradiction.
\end{proof}

By Assumption (\ref{item_Ji_CW}) and standard obstruction theory, the space $s_{i+1}^{-1}(J_i)$ deformation retracts to the image of $\lambda_i$ by a deformation that preserves the fibers of \eqref{eqn_qi_on_Ji}. More precisely, there exists a homotopy
\[
H_i: s_{i+1}^{-1}(J_i)\times [0,1]\to s_{i+1}^{-1}(J_i),
\]
such that $H_i(p,0) = p$, $H_i(p,1) = \lambda_i(q_i(p))$, and $q_i\circ H_i(p,t) = q_i(p)$.
Extend $H_i$ to a map on $\bar{H}_i: s_{i+1}^{-1}(L_i)\times [0,1]\to s_{i+1}^{-1}(L_i)$ as follows:
\[
\bar{H}_i(p, t) = 
\begin{cases}
H_i(p,t) \quad \text{ if } p\in s_{i+1}^{-1}(J_i),\\
p \quad \quad \text{ if } p\in s_{i+1}^{-1}(L_{i+1}).
\end{cases}
\]

\begin{Lemma}
\label{lem_H_bar_continuous}
      The map $\bar{H}_i$ is continuous. 
\end{Lemma}
\begin{proof}[Proof of Lemma \ref{lem_H_bar_continuous}]
The argument is similar to Lemma \ref{lem_lambda_bar_continuous}. 
Suppose $\{(p_k,t_k)\}_{k\in \mathbb{N}}$ is a sequence in $s_{i+1}^{-1}(L_i)\times [0,1]$ that converges to $(p,t)$, we show that $\{\bar{H}_i(p_k,t_k)\}$ converges to $\bar{H}_i(p,t)$. By Assumption (\ref{item_closed}), $s_{i+1}^{-1}(J_i)\times [0,1]$ is an open subset of the domain of $\bar{H}_i$. So, if $p\in s_{i+1}^{-1}(J_i)$, then the convergence follows from the continuity of $H_i$. If $p\in s_{i+1}^{-1}(L_{i+1})$, we use an argument by contradiction. Assume  $\{\bar{H}_i(p_k,t_k)\}$ does not converge to $\bar{H}_i(p,t)$. Then since $s_{i+1}^{-1}(L_i)$ is a compact metrizable space, there is a subsequence of $\{\bar{H}_i(p_k,t_k)\}$ that converges to a point $x\neq \bar{H}_i(p,t)$ in $s_{i+1}^{-1}(L_i)$. We still use $\{\bar{H}_i(p_k,t_k)\}$ to denote the subsequence. 

Note that $q_i\circ H_i(p_k,t_k) = q_i(p_k)$. Hence we have 
\[
q_i(x) = \lim_{k\to \infty}q_i(\bar{H}_i(p_k,t_k)) = \lim_{k\to \infty}q_i(p_k) = q_i(p) \in s_i^{-1}(L_{i+1}).
\]
Therefore, we must have $x = p = \bar{H}_i(p,t)$, contradicting the assumptions.
\end{proof}

By Lemma \ref{lem_H_bar_continuous}, we know that $s_{i+1}^{-1}(L_i)$ deformation retracts to $\Imm(\bar{\lambda}_i)$. Consider the composition map
\[
\Imm(\bar{\lambda}_i) \hookrightarrow s_{i+1}^{-1}(L_i) \xrightarrow{q_i|_{s_{i+1}^{-1}(L_i) }} s_i^{-1}(L_i).
\]
We know that the inclusion map above is a homotopy equivalence, and the composition map is the identity map. Therefore, the map $q_i|_{s_{i+1}^{-1}(L_i) }$ is a homotopy equivalence.

The map $q_i$ is the push-out of $q_i|_{s_{i+1}^{-1}(L_i) }$ with the inclusion $s_{i+1}^{-1}(L_i)\hookrightarrow X_{i+1}$. Since the inclusion of $s_{i+1}^{-1}(L_i)$ in $X_{i+1}$ is a cofibration, and the push-out of a homotopy equivalence with a cofibration is also a homotopy equivalence, we conclude that the map $q_i:X_{i+1}\to X_i$ is a homotopy equivalence.

By Condition (\ref{item_LN}), the map $q$ is a finite composition of the maps $q_i$. So $q$ is a homotopy equivalence, and the proof of Proposition \ref{prop_quotient_with_contractible_fiber} is finished.
\end{proof}

\subsection{Configurations spaces on Euclidean spaces}
Let $d$ be a positive integer.  Later, we will take $d$ to be the dimension of $M$. 
In this subsection, we use Proposition \ref{prop_quotient_with_contractible_fiber} to show that $Q_n^\partial$ is a homotopy equivalence if $M$ is a compact codimension-zero submanifold in $\bR^d$ with an algebraic boundary.

We first recall the following notation from \cite{Sinha-cpt}.

\begin{Definition}
Define $\widetilde{A}_n[\mathbb{R}^d]= (S^{d-1})^{C_2(\underline{n})}\times I^{C_3(\underline{n})}$,  $\widetilde{A}_n\langle \mathbb{R}^d\rangle= (S^{d-1})^{C_2(\underline{n})}$. 
\end{Definition}

Recall that the maps $\alpha_n$ and $\beta_n$ were defined above Definition \ref{def_compactify_w_boundary}.
Let 
$\tilde{\alpha}_n: C_n(\bR^d) \to \widetilde{A}_n[\bR^d]$ be the composition of $\alpha_n$ with the projection map from $A_n[\bR^d]$ to $\widetilde{A}_n[\bR^d]$, and similarly, let $\tilde{\beta}_n: C_n(\bR^d) \to \widetilde{A}_n\langle\bR^d\rangle$ be the composition of $\beta_n$ with the projection map from $A_n\langle \bR^d\rangle$ to $\widetilde{A}_n\langle\bR^d\rangle$.

\begin{Definition}
Let
$\widetilde{C}_n[\mathbb{R}^d]$ be the closure of $\tilde{\alpha}_n(C_n(\bR^d))$ in $\widetilde{A}_n[\bR^d]$, let $\widetilde{C}_n\langle \bR^d\rangle$ be the closure of $\tilde{\beta}_n(C_n(\bR^d))$ in $\widetilde{A}\langle \bR^d\rangle$. 
\end{Definition}

\begin{Definition}
    Let $\widetilde{Q}_n: \widetilde{C}_n[\mathbb{R}^d] \to \widetilde{C}_n\langle \mathbb{R}^d\rangle$ be the restriction of the projection map from $\widetilde{A}_n[\mathbb{R}^d]$ to $\widetilde{A}_n\langle \bR^d \rangle$. 
\end{Definition}

Recall that $C_n[\bR^d]$ is a manifold with corners whose strata are labeled by f-trees (\cite[Section 4]{Sinha-cpt}). We follow the notation from \cite{Sinha-cpt} and denote the stratum labeled by the f-tree $T$ as $C_T(\bR^d)$, and denote its closure in  $\widetilde{C}_n[\bR^d]$ as $\widetilde{C}_T[\bR^d]$.  Following \cite{Sinha-cpt}, let $V(T)$ denote the set of internal vertices of $T$, let $v_0 \in V(T)$ denote the root vertex, let $V^i(T)=V(t)\setminus\{v_0\}$. For $v\in V(T)$,  let $\# v$ denote the number of edges such that $v$ is initial. 
We recall the following result.

\begin{Lemma}[{\cite[Theorem 4.15]{Sinha-cpt}}]
\label{lem_tilde_C_T_product}
The space $C_T[\bR^d]$ is diffeomorphic to 
\[
C_{\# v_0}[\bR^d] \times \big(\widetilde{C}_{\# v}[\bR^d]\big)^{V^i(T)}
\]
as manifolds with corners.
\end{Lemma}

We introduce the following notation.
\begin{Definition}
Suppose $n\ge 2$.
Let $\sigma=(i_1,\dots,i_n)$ be a permutation of $\underline{n}$, and let $\vec{v}\in S^{d-1}$. 
Define $\widetilde{Col}_{\sigma, \vec{v}}(\bR^d)$ to be the subset of $\widetilde{C}_n(\bR^d)$ consisting of configurations $x$ such that $\pi_{i_a,i_b} = \vec{v}$ for all $a<b$. (Recall that the maps $\pi_{ij}$ were defined at the beginning of Section \ref{subsec_inclusion_map}.) Let $\widetilde{Col}_n(\bR^d)\subset \widetilde{C}_n(\bR^d)$ denote the union of $\widetilde{Col}_{\sigma,\vec{v}}(\bR^d)$ for all $\sigma$, $\vec{v}$.
Let $\widetilde{Col}_{\sigma, v}[\bR^d]$ and $\widetilde{Col}_n[\bR^d]$ denote the closures of $\widetilde{Col}_{\sigma, v}(\bR^d)$ and $\widetilde{Col}_n(\bR^d)$ respectively in $\widetilde{C}_n[\bR^d]$.
\end{Definition}

\begin{remark}
\label{rem_col_fiber}
The space $\widetilde{Col}_{\sigma, v}[\bR^d]$ is diffeomorphic to an associahedron with dimension $n-2$. The space $\widetilde{Col}_n[\bR^d]$ has $n!$ connected components, and each component is diffeomorphic to the product of an associahedron with $S^{d-1}$. 
\end{remark}

\begin{Definition}
Suppose $T$ is an f-tree and $\Lambda\subset V(T)\setminus \{v_0\}$.
Define $S_{T,\Lambda}[\bR^d]$ to be the set of 
\[
(x_{v_0}, (x_v)_{v\in V^i(T)})\in C_{\# v_0}[\bR^d] \times \big(\widetilde{C}_{\# v}[\bR^d]\big)^{V^i(T)}
\] such that $x_v\in \widetilde{Col}_{\# v}[\bR^d]$ for all $v\in \Lambda$. We also abuse notation and let $S_{T,\Lambda}[\bR^d]$ denote the corresponding subspace of $C_T[\bR^d]$ via the canonical diffeomorphism given by Lemma \ref{lem_tilde_C_T_product}.
\end{Definition}

\begin{remark}
\label{rmk_generic_fiber}
For a generic element $x\in S_{T,\Lambda}[\bR^d]$, the set $Q_{n}^{-1}({Q}_n(x))$ is contained in ${S}_{T,\Lambda}[\bR^d]$ and is a product of associahedrons with dimensions $(\#v -2)$ for $v\in \Lambda.$ So the dimension of ${Q}_n^{-1}({Q}_n(x))$ for a generic $x\in {S}_{T,\Lambda}[\bR^d]$ is 
\begin{equation}
    \label{eqn_dim_fiber}
    d_{T,\Lambda} := \sum_{v\in \Lambda}(\#v -2 ).
\end{equation}
\end{remark}

In order to apply Proposition \ref{prop_quotient_with_contractible_fiber}, we need to show that certain inclusion maps are cofibrations. To facilitate the argument, we introduce the following terminology.

\begin{Definition}
    A subset $X$ of an Euclidean space $\bR^N$ is called a \emph{real affine algebraic variety with corners} if it has the form
    \[
    \bigcup_{i=1}^k \big(\bigcap_{j=1}^{n_i} S_{i,j}\big),
    \]
    such that each $S_{i,j}$ has the form $\{x\in \bR^N|f(x)=0\}$ or $\{x\in \bR^N|f(x)\ge 0\}$ with $f$ being a polynomial.
    A subset $Y$ of $X$ is called a \emph{subvariety with corners} if $Y$ is also a real affine algebraic variety with corners in $\bR^N$. 
\end{Definition}

\begin{remark}
    For simplicity, in the following, we will call a \emph{real affine algebraic variety with corners} as a \emph{variety}, and a \emph{subvariety with corners} as a \emph{subvariety}. 
\end{remark}

\begin{Definition}
    Suppose $X$ and $Y$ are varieties. We say that a map $f:X\to Y$ is \emph{algebraic}, if it is locally given by maps whose coordinates are rational functions. We say that $X$ and $Y$ are \emph{isomorphic}, if there exists an algebraic homeomorphism whose inverse is also algebraic.
\end{Definition}

By a straightforward extension of the arguments in \cite{whitney1992elementary}, we know that every variety $X$ is homeomorphic to a CW complex.  Moreover, if $Y$ is a subvariety of $X$, then the pair $(X,Y)$ is homeomorphic to a CW pair, and hence the inclusion of $Y$ in $X$ is a cofibration. By \cite[Theorems 4.1, 5.14]{Sinha-cpt}, $C_n[\bR^d]$ is a subvariety of $A_n[\bR^d]$, and $C_n\langle \bR^d\rangle$ is a subvariety of $A_n\langle \bR^d\rangle$. It is also clear that $C_T[\bR^d]$ is a subvariety of $C_n[\bR^d]$, and $S_{T,\Lambda}[\bR^d]$ is a subvariety of $C_T[\bR^d]$.

\begin{Lemma}
\label{lem_C_T_proj_subvariety}
    Let $T$ be an f-tree. Then ${Q}_n({C}_{T}[\bR^d])$ is a subvariety of ${C}_n\langle \bR^d\rangle$.
\end{Lemma}

\begin{proof}
In the tree $T$, we say that a vertex $v$ \emph{lies over} a vertex $w$ if the shortest path from the root $v_0$ to $v$ contains $w$. Note that by our convention, a vertex always lies over itself.

Consider the map
\[
\iota_T:  A_{\# v_0} \langle \bR^d\rangle \times \big(\widetilde{A}_{\# v}\langle \bR^d\rangle \big)^{V^i(T)} \to \widetilde{A}_n\langle \bR^d\rangle
\]
defined as follows. For each $(i,j)\in C_2(\underline{n})$, let $v$ be the highest vertex in $T$ such that both $i$ and $j$ lie over $v$. Let $v_i$ be the vertex over $v$ such that $i$ lies over $v_i$ and that $v_i$ is connected to $v$ by an edge in $T$, and define $v_j$ similarly. Then for $x\in A_{\# v_0} \langle \bR^d\rangle \times \big(\widetilde{A}_{\# v}\langle \bR^d\rangle \big)^{V^i(T)}$, the $(S^{d-1})$--coordinate labeled by $(i,j)\in C_2(\underline{n})$  of $\iota_T(x)$ is defined to be the $(S^{d-1})$--coordinate labeled by $(v_i,v_j)\in C_2(\underline{\# v})$ of the projection of $x$ to $\widetilde{A}_{\# v}\langle \bR^d\rangle$. 

For each $i\in \underline{n}$, let $w_i$ be the vertex such that $w_i$ is connected to $v_0$ by an edge and $i$ lies over $w_i$. The $I$--coordinate of $\iota_T(x)$ labeled by $i\in \underline{n}$ is defined to be the $I$--coordinate of $x$ labeled by $w_i$. 

Then the composition of 
\[
C_{\# v_0} [\bR^d] \times \big(\widetilde{C}_{\# v} [\bR^d] \big)^{V^i(T)}
\xrightarrow{\cong}
{C}_{T}[\bR^d] 
\xrightarrow{\widetilde{Q}_n}
{C}_n\langle \bR^d\rangle
\]
is equal to the composition of 
\[
C_{\# v_0} [\bR^d] \times \big(\widetilde{C}_{\# v} [\bR^d] \big)^{V^i(T)}
\xrightarrow{Q_{\# v_0}\times (\tilde{Q}_{\# v})^{V^i(T)}} 
C_{\# v_0} \langle \bR^d\rangle \times \big(\widetilde{C}_{\# v}\langle \bR^d\rangle \big)^{V^i(T)} 
\xrightarrow{\iota_T} 
{C}_n\langle \bR^d\rangle.
\]
Since $Q_{\# v_0}\times (\tilde{Q}_{\# v})^{V^i(T)}$ is surjective and it is clear that $\iota_T$ maps subvarieties of $A_{\# v_0} \langle \bR^d\rangle \times \big(\widetilde{A}_{\# v}\langle \bR^d\rangle \big)^{V^i(T)}$ to subvarieties of $\widetilde{A}_n\langle \bR^d\rangle$, the result is proved.
\end{proof}

\begin{Lemma}
\label{lem_S_T_Lambda_proj_subvariety}
    Let $T$ be an f-tree and $\Lambda\subset V(T)$. Then ${Q}_n({S}_{T,\Lambda}[\bR^d])$ is a subvariety of ${C}_n\langle \bR^d\rangle$.
\end{Lemma}

\begin{proof}
Let $\iota_T$ be as in the proof of Lemma \ref{lem_C_T_proj_subvariety}.
    Note that by definition, the image of ${S}_{T,\Lambda}[\bR^d]$ under the composition 
\[
{S}_{T,\Lambda}[\bR^d]
\hookrightarrow  
C_{\# v_0}[\bR^d]\times \big(\widetilde{C}_{\# v}[\bR^d]\big)^{V^i(T)}
\xrightarrow{(Q_{\# v_0}\times \tilde{Q}_{\# v})^{V^i(T)}} 
C_{\# v_0}\langle\bR^d\rangle \times \big(\widetilde{C}_{\# v}\langle \bR^d\rangle \big)^{V'(T)}
\]
is a subvariety, because it is a subset defined by colinearity conditions. Therefore, its image after composing with the map $\iota_T$ is a subvariety of ${C}_n\langle \bR^d\rangle$.
\end{proof}

Now suppose $M$ is a compact codimension-0 submanifold in $\bR^d$ whose boundary is the zero locus of a polynomial. We will later take $M$ to be a ball in $\bR^d$. Define $C_T[M]= C_n[M]\cap C_T[\bR^d]$, $S_{T,\Lambda}[M] = C_n[M]\cap C_T[\bR^d]$. Then $C_T[M]$ and $S_{T,\Lambda}[M]$ are subvarieties of $C_n[\bR^d]$. 

\begin{Lemma}
\label{lem_Qn_equiv_M_algebraic}
     Suppose $M$ is a compact codimension-0 submanifold in $\bR^d$ whose boundary is the zero locus of a polynomial. Then the quotient map $Q_n:C_n[M]\to C_n\langle M\rangle$ is a homotopy equivalence.
\end{Lemma}
\begin{proof}
We show that the map $Q_n$ satisfies all the assumptions of Proposition \ref{prop_quotient_with_contractible_fiber}. 

    By the assumptions on $M$, we know that $C_n[M]$ is a compact subvariety of $C_n[\bR^d]$, and hence it is compact and metrizable. It is also clear that for each $y\in C_n\langle M\rangle$, the preimage $Q_n^{-1}(\{y\})$ is a finite product of associahedra, and hence it is homeomorphic to a closed ball. 
    
    Let $L_i\subset C_n\langle M\rangle$ be the set of points $y$ such that $Q_n^{-1}(\{y\})$ has dimension at least $i$. Then $Q_n^{-1}(L_i)$ is the union of all $S_{T,\Lambda}[M]$ with $d_{T,\Lambda}\ge i$ (see Equation \eqref{eqn_dim_fiber}). By Lemma \ref{lem_S_T_Lambda_proj_subvariety}, $L_i$ is a finite union of subvarieties of $C_n\langle M\rangle$, so Assumptions (\ref{item_closed}), (\ref{item_cofib}) are verified. Assumption (\ref{item_LN}) is also clear by the definition of $d_{T,\Lambda}$. Lemma \ref{lem_S_T_Lambda_proj_subvariety} implies that $L_i\setminus L_{i+1}$ is homeomorphic to a CW complex. The set $J_i$ is the disjoint union of $S_{T,\Lambda}[M]\setminus L_{i+1}$ for all $d_{T,\Lambda} = i$, and each $S_{T,\Lambda}[M]\setminus L_{i+1}$ is an open subset of $J_i$. 
    By Remark \ref{rem_col_fiber}, it is clear that 
    \[
    Q_n|_{S_{T,\Lambda}[M]\setminus L_{i+1}}: S_{T,\Lambda}[M]\setminus L_{i+1} \to Q_n(S_{T,\Lambda}[M]\setminus L_{i+1})
    \] 
    is a fiber bundle. So the map $Q_n$ satisfies all the assumptions of Proposition \ref{prop_quotient_with_contractible_fiber}.
\end{proof}

Similarly, we have the following result.

\begin{Lemma}
\label{lem_Qn_partial_equiv_M_algebraic}
     Suppose $M$ is a compact codimension-0 submanifold in $\bR^d$ whose boundary is the zero locus of a polynomial.  Then the quotient map $Q_n^\partial :C_n[M,\partial]\to C_n\langle M,\partial\rangle$ is a homotopy equivalence.
\end{Lemma}

\begin{proof}
    The proof is essentially the same as Lemma \ref{lem_Qn_equiv_M_algebraic}. The only difference is that the configuration spaces are now stratified by f-trees with $n+2$ leaves. 
\end{proof}

Lemmas \ref{lem_Qn_equiv_M_algebraic} and \ref{lem_Qn_partial_equiv_M_algebraic} have the following corollaries.

\begin{Corollary}
\label{cor_M_R^d_Q_map}
Let $M$ be an open subset of $\bR^n$.
The quotient map $Q_n: C_n[M]\to C_n\langle M\rangle$ is a homotopy equivalence. 
\end{Corollary}

\begin{proof}
Let $\{M_i\}_{i\ge 1}$ be a sequence of compact polytopes in $M$ such that $M_i\subset \inte(M_{i+1})$ and $\cup_i M_i = M$. Then $C_n[M]$ is the union of $C_n[M_i]$, and hence it is homotopy equivalent to a CW complex. Similarly, $C_n\langle M\rangle$ is homotopy equivalent to a CW complex.

By Lemma \ref{lem_Qn_equiv_M_algebraic}, map $Q_n: C_n[M_i]\to C_n\langle M_i\rangle$ induces isomorphisms on all homotopy groups. Since the homotopy groups of $C_n[ M]$ are the direct limits of the homotopy groups of $C_n[ M_i]$, and the homotopy groups of $C_n\langle M\rangle$ are the direct limits of the homotopy groups of $C_n\langle M_i\rangle$, we know that the map $Q_n: C_n[M]\to C_n\langle M\rangle$ is a weak homotopy equivalence, and hence a homotopy equivalence. 
\end{proof}

Let $\bR_{\ge 0}^d = [0,+\infty)\times \bR^{d-1}$. The same argument proves the following corollary of Lemma \ref{lem_Qn_partial_equiv_M_algebraic}.
\begin{Corollary}
\label{cor_M_R^d_ge_Q_map}
Let $M$ be an open subset of $\bR^d_{\ge 0}$ that intersects $\{0\}\times \bR^{d-1}$ non-trivially. 
Then the quotient map $Q_n^\partial: C_n[M,\partial ]\to C_n\langle M,\partial\rangle$ is a homotopy equivalence. \qed
\end{Corollary}

\subsection{The map $Q_n^\partial$ for general manifolds}
    Now we show that the map $Q_n^\partial :C_n[ M,\partial] \to C_n\langle M,\partial \rangle$ is a homotopy equivalence for all $M$.
Let 
\[
\pi_{[M]}: A_n[ M,\partial] \to M^n,\quad \pi_{\langle M \rangle}: A_n\langle M,\partial \rangle \to M^n
\]
be the projections to the $M^n$ factor, then $\pi_{\langle M \rangle} \circ Q_n^\partial = \pi_{[M]}$ on $C_n[M,\partial]$.

Since $C_n[M,\partial]$ is a manifold with corners, it is homeomorphic to a CW complex.  Consider a locally finite covering of $M$ by $\{B_s\}_{s\in S}$ such that for each subset $I$ of $S$, the intersection $\cap_{s\in I} B_s$ is either empty or diffeomorphic to a convex polyhedron, and for $I_1\subset I_2\subset S$, the pair $(\cap_{s\in I_1} B_s, \cap_{s\in I_2} B_s)$ is diffeomorphic to the embedding of a convex polyhedron in a convex polyhedron. This gives a locally finite covering of $C_n\langle M, \partial \rangle$ indexed by $(S)^n$, such that for $\vec s = (s_1,\dots,s_n)\in (S)^n$, the corresponding subset $A_{(s_1,\dots,s_n)}$ of $C_n\langle M, \partial \rangle$ is the preimage of $B_{s_1}\times \dots\times B_{s_n}$ under $\pi_{[M]}$.   
Then for $I_1\subset I_2\subset (S)^n$, the pair $(\cap_{\vec s\in I_1} A_{\vec s}, \cap_{\vec s\in I_2} A_{\vec s})$ is homeomorphic to a variety and its subvariety. Hence  $C_n\langle M,\partial\rangle$ is homotopy equivalent to a CW complex. 

A map between CW complexes is a homotopy equivalence if and only if it induces an isomorphism on $\pi_1$ and induces isomorphisms on all homology groups with all local coefficient systems. Therefore, by the Seifert--van Kampen theorem and the Mayer--Vietoris sequence, we only need to show for every $x\in M^n$, there exists an open neighborhood $V_x$ of $x$, such that for every open subset $U\subset V_x$, the map
\begin{equation}
\label{eqn_restr_Q_to_Ux}
  Q_n^\partial|_{(\pi_{[M]})^{-1}(U)}: (\pi_{[M]})^{-1}(U) \to (\pi_{\langle M \rangle})^{-1}(U)
\end{equation}
is a weak homotopy equivalence. This is established in the following lemma.

\begin{Lemma}
For every point $x\in M^n$, there exists an open neighborhood $V_x\subset M^n$ of $x$ such that for every open subset $U\subset V_x$, the map
\eqref{eqn_restr_Q_to_Ux} is a weak homotopy equivalence. 
\end{Lemma}

\begin{proof}
We use induction on $n$. The case for $n=1$ is trivial. From now on, assume $n>1$, and assume that the statement is true for all $n'<n$.

If $x$ is not on the diagonal of $M^n$, we may choose a sufficiently small neighborhood $U_x$ of $x$ such that $(\pi_{[M]})^{-1}(U)$ is locally given by a product of configuration spaces with smaller numbers of points. So the desired result follows from the induction hypothesis and Lemma \ref{lem_prod_universal_hmty_equiv} below.

If $x =(x_0,\dots,x_0)$ and $x_0\in \inte(M)$, then there is a neighborhood $N$ of $x_0$ in $M$ such that $N$ is diffeomorphic to $\bR^d$. Let $U=N^n\subset M^n$, the desired result then follows from Corollary \ref{cor_M_R^d_Q_map}. 

If $x=(x_0,\dots,x_0)$ and $x_0\in \partial M$, then one can choose an open neighborhood $N$ of $x_0$ in $M$ such that $N$ is diffeomorphic to $\bR^d_{\ge 0}$, and such that $p,q\in N$, where $p,q$ are the  points in Definition \ref{def_compactify_w_boundary}.
Let $U = N^n\subset M^n$, the desired result then follows from Corollary \ref{cor_M_R^d_ge_Q_map}.
\end{proof}

Now we state and prove Lemma \ref{lem_prod_universal_hmty_equiv}.
\begin{Lemma}
\label{lem_prod_universal_hmty_equiv}
    Suppose $X_i$ and $Y_i$ are second countable. Suppose $f_i:X_i\to Y_i$ ($i=1,2$) are two continuous surjective maps such that for every open subset $U_i$ of $Y_i$, the map $f_i|_{f_i^{-1}(U_i)}: f_i^{-1}(U_i) \to U_i$ is a weak homotopy equivalence. Let $f = f_1\times f_2: X_1\times X_2\to Y_1\times Y_2$. Then for every open subset $U$ of $Y_1\times Y_2$, the map $(f)|_{f^{-1}(U)}: f^{-1}(U) \to U$ is a weak homotopy equivalence. 
\end{Lemma}
\begin{proof}
    Write $U$ as the union of a countable collection of open subsets $\{A_i\}_{i\in \mathbb{N}}$ such that each   $A_i$ is the product of open subsets of $Y_1$ and $Y_2$. For every finite subset $I\subset \mathbb{N}$, let $U_I = \cap_{i\in I} A_i$. By the assumptions, the map $f|_{f^{-1}(U_I)}: f^{-1}(U_I) \to U_I$ is a weak homotopy equivalence. Let $U_i = \cup_{j\le i} A_j$. Then by the Seifert--van Kampen theorem and the Mayer--Vietoris sequence, we know that the map $f|_{f^{-1}(U_i)}: f^{-1}(U_i) \to U_i$ is a weak homotopy equivalence. Since the homotopy groups of $U$ and $f^{-1}(U)$ are given by the direct limits of the homotopy groups of $U_i$ and $f^{-1}(U_i)$ respectively, the desired result is proved. 
\end{proof}

In conclusion, we have proved that $Q^\partial_n:C_n[M,\partial]\to C_n \langle M,\partial \rangle$ is a homotopy equivalence.

%% file: main.bbl
\begin{bibdiv}
\begin{biblist}

\bib{BG2019}{article}{
      author={Budney, Ryan},
      author={Gabai, David},
       title={Knotted 3-balls in ${S}^4$},
        date={2019},
     journal={arXiv preprint, arXiv:1912.09029},
}

\bib{BG2023}{article}{
      author={Budney, Ryan},
      author={Gabai, David},
       title={On the automorphism groups of hyperbolic manifolds},
        date={2023},
     journal={arXiv preprint, arXiv:2303.05010},
}

\bib{BottTu}{book}{
      author={Bott, Raoul},
      author={Tu, Loring~W.},
       title={Differential forms in algebraic topology},
      series={Graduate Texts in Mathematics},
   publisher={Springer-Verlag, New York-Berlin},
        date={1982},
      volume={82},
        ISBN={0-387-90613-4},
      review={\MR{658304}},
}

\bib{engelking1989general}{book}{
      author={Engelking, R},
       title={General topology},
      series={Sigma Series in Pure Mathematics},
   publisher={Heldermann Verlag Berlin},
        date={1989},
      volume={6},
}

\bib{GK}{article}{
      author={Goodwillie, Thomas~G.},
      author={Klein, John~R.},
       title={Multiple disjunction for spaces of smooth embeddings},
        date={2015},
        ISSN={1753-8416,1753-8424},
     journal={J. Topol.},
      volume={8},
      number={3},
       pages={651\ndash 674},
         url={https://doi.org/10.1112/jtopol/jtv008},
      review={\MR{3394312}},
}

\bib{GW}{article}{
      author={Goodwillie, Thomas~G.},
      author={Weiss, Michael},
       title={Embeddings from the point of view of immersion theory. {II}},
        date={1999},
        ISSN={1465-3060,1364-0380},
     journal={Geom. Topol.},
      volume={3},
       pages={103\ndash 118},
         url={https://doi.org/10.2140/gt.1999.3.103},
      review={\MR{1694808}},
}

\bib{lee2010introduction}{book}{
      author={Lee, John},
       title={Introduction to topological manifolds},
     edition={2},
      series={Graduate Texts in Mathematics},
   publisher={Springer-Verlag},
        date={2010},
      volume={202},
}

\bib{RW}{incollection}{
      author={Randal-Williams, Oscar},
       title={Diffeomorphisms of discs},
        date={2023},
   booktitle={I{CM}---{I}nternational {C}ongress of {M}athematicians. {V}ol. 4.
  {S}ections 5--8},
   publisher={EMS Press, Berlin},
       pages={2856\ndash 2878},
      review={\MR{4680344}},
}

\bib{Sinha-cpt}{article}{
      author={Sinha, Dev~P.},
       title={Manifold-theoretic compactifications of configuration spaces},
        date={2004},
        ISSN={1022-1824,1420-9020},
     journal={Selecta Math. (N.S.)},
      volume={10},
      number={3},
       pages={391\ndash 428},
         url={https://doi.org/10.1007/s00029-004-0381-7},
      review={\MR{2099074}},
}

\bib{Sinha-cosimplicial}{article}{
      author={Sinha, Dev~P.},
       title={The topology of spaces of knots: cosimplicial models},
        date={2009},
        ISSN={0002-9327,1080-6377},
     journal={Amer. J. Math.},
      volume={131},
      number={4},
       pages={945\ndash 980},
         url={https://doi.org/10.1353/ajm.0.0061},
      review={\MR{2543919}},
}

\bib{Watanabe2020}{article}{
      author={Watanabe, Tadayuki},
       title={Theta-graph and diffeomorphisms of some 4-manifolds},
        date={2020},
     journal={arXiv preprint, arXiv:2005.09545},
}

\bib{Weiss1}{article}{
      author={Weiss, Michael},
       title={Embeddings from the point of view of immersion theory. {I}},
        date={1999},
        ISSN={1465-3060,1364-0380},
     journal={Geom. Topol.},
      volume={3},
       pages={67\ndash 101},
         url={https://doi.org/10.2140/gt.1999.3.67},
      review={\MR{1694812}},
}

\bib{whitney1992elementary}{article}{
      author={Whitney, Hassler},
       title={Elementary structure of real algebraic varieties},
        date={1992},
     journal={Hassler Whitney Collected Papers},
       pages={456\ndash 467},
}

\end{biblist}
\end{bibdiv}
